\documentclass[a4paper,12pt]{amsart}
\usepackage{latexsym,amsmath,amssymb,amsthm,color}
\usepackage{constants}
\newtheorem{thm}{Theorem}[section]
\newtheorem{lemma}[thm]{Lemma}
\newtheorem{cor}[thm]{Corollary}

\newtheorem{define}[thm]{Definition}
\newtheorem{prop}[thm]{Proposition}

\numberwithin{equation}{section}
\newconstantfamily{eps}{symbol=\varepsilon}
\newconstantfamily{c}{symbol=c}
\newconstantfamily{al}{symbol=\alpha}
\newconstantfamily{be}{symbol=\beta}
\newconstantfamily{de}{symbol=\delta}
\newconstantfamily{La}{symbol=\Lambda}
\def\res{\hbox{ {\vrule height .3cm}{\leaders\hrule\hskip.3cm}}\hskip5.0\mu}

\input{amssym.def}

\title[the blow up method for Brakke flows]{the blow up method for \\
Brakke flows: networks near triple junctions}

\author[Y. Tonegawa]{Yoshihiro Tonegawa}

\author[N. Wickramasekera]{Neshan Wickramasekera}
\address{Department of Mathematics, Tokyo Institute of Technology,  2-12-1 Ookayama, Meguro-ku, 152-8551, Japan.}
\email{tonegawa@math.titech.ac.jp}
\address{DPMMS, University of 
Cambridge, Cambridge, CB3 0WB, United Kingdom.}
\email{N.Wickramasekera@dpmms.cam.ac.uk}
\date{}


\thanks{Y.T. is partially supported by JSPS Grant-in-aid for scientific research (A) $\#$25247008, (S) $\#$21224001
and challenging exploratory research $\#$23654057.}

\begin{document}
\setlength\parskip{5pt}
\begin{abstract}
We introduce a parabolic blow-up method to study the asymptotic behavior of a \emph{Brakke flow of planar networks} (i.e.\ a 1-dimensional Brakke flow in a two dimensional region) weakly close in a space-time region to a static multiplicity 1 triple junction $J$. 
We show that such a network flow is regular in a smaller space-time region, in the sense that it
consists of three curves coming smoothly together at a single point at 120 degree angles, staying smoothly close to $J$ and moving smoothly. 
Using this result and White's stratification theorem, we deduce that whenever a Brakke flow of networks in a space-time region ${\mathcal R}$ has no \emph{static} tangent flow with density $\geq2$, there exists a closed subset $\Sigma \subset {\mathcal R}$ 
of parabolic Hausdorff dimension at most 1 such that the flow is classical in ${\mathcal R} \setminus \Sigma$, i.e.\ near every point in ${\mathcal R} \setminus \Sigma$, the flow, if non-empty, consists of either an embedded curve moving smoothly or three embedded curves meeting smoothly at a single point at 120 degree angles and moving smoothly. In particular, such a flow is classical at all times except for a closed set of times of ordinary Hausdorff dimension at most $\frac{1}{2}$.
\end{abstract}

\maketitle
\makeatletter
\@addtoreset{equation}{section}                  
\renewcommand{\theequation}{\thesection.\@arabic\c@equation}
\makeatother

\tableofcontents

\section{Introduction}
Our goal in this paper is to introduce a general framework---a parabolic blow up method---to study the asymptotic nature of a multiplicity 1 Brakke flow near certain generic singularities of the flow. The theorems we prove here using this framework, which we shall describe shortly, concern the simplest non-trivial situation, namely, the asymptotic behavior near a static triple junction of a Brakke flow of planar \emph{networks}, i.e.\   a 1-parameter family of 
1-dimensional sets (corresponding to integral 1-varifolds) moving by generalized curvature in a domain 
$U\subset{\mathbb R}^2$ over a time interval.

In what might be called the classical setting, such a 1-dimensional flow consists of a locally finite 
union of smoothly embedded open curves moving smoothly 
with velocity equal to the 
curvature vector at each point and time, and such that at each of their boundary points (in $U$), three curves meet smoothly at 120 degree angles. Since globally in space the flow decreases the
total length of the curves in time, the curvature flow of networks
models the motion of grain boundaries driven by interfacial surface tension \cite{Gurtin}.
The 120 degree angle condition in this context is called the \emph{Herring condition}. 

While such classical solutions may stay classical time-globally in some special cases
\cite{Garcke,Ikota2005,Kinderlehrer,MMN,Mantegazza}, in general, various singularities may occur in finite time. For instance, physically, in the motion of grain boundaries, 
one observes that 
two or more triple junctions collide with each other and small grains are eliminated. 
This is a process of grain coarsening which should be an integral part of the mathematical modeling of the motion.
Motivated partly by such phenomena, in the pioneering work \cite{Brakke}, Brakke introduced 
a generalized notion of mean curvature flow (abbreviated MCF hereafter) 
using the notion of varifolds in Geometric Measure Theory, 
and studied existence and regularity of surfaces of any dimension and codimension moving by mean curvature.
Brakke's MCF (which we also call ``Brakke flow'', see Sec. \ref{defBF} for the definition)
naturally accommodates flows with
singularities, but allows the possibility of sudden loss of measure and non-uniqueness.

While the study of regularity of various classes of stationary varifolds (generalized minimal submanifolds), which are the equilibrium solutions of the Brakke flow, 
has seen several advances in the past 50 years or so, 
much less has been known concerning  Brakke flows apart from Brakke's own original work \cite{Brakke}.
Recently, the work of Kasai and the first author \cite{Kasai-Tonegawa} and 
of the first author \cite{Tonegawa} gave a new, streamlined  proof of a generalization of Brakke's local regularity theorem (\cite{Brakke})
which establishes a.e$.$ smoothness in time and space under the hypothesis that the moving surfaces have multiplicity 1 a.e. 
Roughly speaking, just like Allard's regularity theorem \cite{Allard}
for stationary varifolds, 
the results in \cite{Brakke,Kasai-Tonegawa,Tonegawa} show that a nearly flat part of a unit density Brakke flow is necessarily a smooth MCF, i.e. locally $C^{\infty}$ embedded submanifolds in space-time and moving smoothly
by the mean curvature.  \cite{Brakke,Kasai-Tonegawa,Tonegawa} however do not give any structural information about the flow in the vicinity 
of singularities including triple junctions. 
We note that 
there have been important results on Brakke's regularity theorem
when one is interested in special Brakke flows such as  those arising as weak 
limits of smooth MCFs (see for instance \cite{Ecker1,Ecker,White1}) or those produced by Ilmanen's elliptic regularization method (\cite{Ilmanen0}).

For both minimal submanifolds and mean curvature flows, as well as  for numerous other problems in geometric analysis and non-linear PDE,  describing the asymptotic behavior of the objects in question on approach to their singular sets, and understanding the structure of the singular sets themselves, remain largely open major challenges. For multiplicity 1 classes of minimal submanifolds,
the seminal work of Simon \cite{Simon2.1, Simon2, Simon3} established asymptotics near certain singularities, and also the structure results for  the singular sets in the full generality of varying tangent cone types and when there is no topological obstruction to perturbing singularities away. Earlier work of Allard--Almgren \cite{AA1}, Taylor \cite{Taylor2} and White \cite{White4} proved similar results in situations where the tangent cones satisfy more restrictive conditions in addition to the multiplicity 1 condition. 
Recent work of 
the second author \cite{Wick, Wick1} and of Krummel and the second author \cite{KW1, KW2} establish regularity results and asymptotics near singularities (branch points) for  certain classes of minimal and related submanifolds  for which the multiplicity 1 condition either fails or is not assumed a priori.  The work \cite{Taylor2, White4, Simon2, Simon3, Wick1, KW2} establish, for various classes of minimal submanifolds,  fine properties of the singular sets themselves,  such as smoothness or rectifiability.

Among the known results in this direction for MCF is the 
recent deep work of Colding--Minicozzi \cite{Colding, Colding1} (aided also by  the work of Colding--Ilmanen--Minicozzi \cite{Colding0}) which proves, for any flow of hypersurfaces,  uniqueness of the tangent flow whenever a multiplicity 1 shrinking cylinder occurs as one tangent flow, and provides strong structural information on  the singular sets of hypersurface flows whose tangent flows at singularities are all shrinking multiplicity 1 cylinders. In particular, these results apply to MCFs of mean-convex hypersurfaces, for which the earlier work of White \cite{White3, White2, White6} had established that all tangent flows at singularities are shrinking multiplicity 1 cylinders. 
The work of Schulze~\cite{Schulze} established such asymptotics near compact singularities of multiplicity 1 flows. 

Simon's work \cite{Simon2.1, Simon2, Simon3} mentioned above developed  two far reaching methods---one based on an infinite dimensional Lojasiewicz inequality and the other based on the so called blow-up method---for studying the singular sets of minimal submanifolds.  The work of Colding--Minicozzi \cite{Colding, Colding1} and of Schulze \cite{Schulze} referred to above study singularities of MCFs by establishing appropriate Lojasiewicz inequalities. In the present paper we  introduce a parabolic version of the blow-up method for studying fine properties of Brakke flow singularities, and implement it fully in the simplest non-trivial case with moving singularities, namely, for 1-dimensional Brakke flows in the vicinity of a static multiplicity 1 triple junction.

Our main result here (Theorem~\ref{mainreg} below) is a precise version of  the following:

\noindent
{\bf Theorem 1}. \emph{If a 1-dimensional Brakke flow in a planar region is weakly close in  a space-time neighborhood to a static multiplicity 1 triple junction $J$, then in a smaller space-time neighborhood, it is regular in the sense that it consists of three curves coming smoothly together at a single point, staying smoothly close to $J$ and moving by curvature.} 

Here, by \emph{weakly close} we mean that the flow has small space-time $L^{2}$ distance from $J$ and satisfies suitable mass hypotheses at the initial and final times. (See the statement of Theorem~\ref{mainreg}.) 
This is a natural, easily verifiable criterion in the analysis of singularities. For example, if a $1$-dimensional Brakke flow at a space-time singular point  has a tangent
flow equal to  a static multiplicity 1 triple junction, then our theorem is applicable near that point, and implies uniqueness of the tangent flow, and moreover, that in a space-time neighborhood of the point, the flow itself is a regular triple junction moving by curvature.
Using the above result and White's stratification theorem \cite{White0}, we deduce the following partial regularity result (Theorem~\ref{mainpa}) for 1-dimensional Brakke flows in a planar region:

\noindent
{\bf Theorem 2}. \emph{If a 1-dimensional Brakke  flow in a planar region has no static tangent flow consisting of more than 3 half-lines meeting at the origin (or, equivalently, if the density of every static tangent flow is $<2$), 
then the flow is a classical network flow away from a closed singular set of parabolic Hausdorff dimension at most 1; in particular, such a flow is a classical network flow at each instance of time except for a closed set of times of ordinary Hausdorff dimension at most 1/2}. 

The hypothesis  concerning tangent flows in Theorem 2 is motivated by the physics of motion of grain boundaries where the triple junction seems to be the unique stable junction; other types of junctions may form but seem to disappear instantly.
(See more discussion after the statement of Theorem~\ref{mainpa}.)
In this regard, 
Theorem 1 may also be understood as a result about local asymptotic stability of triple junction within a 
broad class of weak varifold solutions and with respect to a weak topology of measure. There have been 
related works which prove global-in-time asymptotic stability of triple junction, i.e. classical 
network curvature flow converges to a triple junction as $t\rightarrow\infty$ \cite{Garcke,Ikota2005,Kinderlehrer,MMN}. 

The above results are formulated and proved here for a class of 1-dimensional flows more general than Brakke flows, in which the ``velocity of motion'' is given by the curvature vector plus any given space-time dependent vector field satisfying an optimal integrability condition. 

The simplicity of the spatial 1-dimensionality of the problem considered here allows us to essentially isolate the difficulties arising from the presence of the time variable. Although some of our arguments here take advantage of the spatial 1-dimensionality,  the overall method introduced here appears to hold promise for much further development.
Indeed, many  of the estimates developed here either directly extend to or can easily be modified
to work for Brakke flows of general dimension and codimension weakly close to certain types of multiplicity 1 tangent flows, including higher dimensional static triple junctions. However, there are also a few ingredients for which the arguments needed in higher dimensions seem to be much more complicated. We shall address  such generalizations elsewhere. 

One may also naturally wonder what could go wrong if the triple junction $J$ is replaced by a 1-dimensional multiplicity 1 stationary junction $J_{N}$ with $N(>3)$ half lines meeting at the origin. 
In this case, the direct analogue of the conclusion of Theorem 1, namely that in the space-time interior the flow consists of $N$ embedded curves coming together smoothly at a single point, is false. For instance in case $N=4$, consider the static junction $J_{4}$ consisting of two intersecting lines at the origin with a 120 degree angle  between them. We may construct a static configuration arbitrarily close to $J_{4}$ with precisely two triple junction singularities by splitting apart $J_{4}$ at the origin into two pairs of half-lines each making a 120 degree angle and connecting their vertices by a short line segment, and imagine non static flows that remain close to this configuration.
From the point of view of our method here (see below for an outline), general uniform regularity estimates fail (as they must in view of the example just mentioned) without further hypotheses in case $N \geq 4$ because the flow need not have the property that 
the moving curve at time $t$ has a singular point of density $\geq N/2$ for a.e.\ $t$. A further complicating issue in this case is that even when the curve does have singularities with the right density, its tangent cones may contain  \emph{higher multiplicity} lines or half-lines. Neither of these issues arises in the case $N = 3.$ 

Smoothly embedded 1-dimensional flows on the other hand cannot stay close to a singular static junction for too long. In higher dimensions, an analogue  is the question of what one can say about minimal surfaces weakly close to a pair of transverse planes (say, in ${\mathbb R}^{3}$). In that case, the difficulties are illustrated by Scherk's surfaces which show that no uniform estimates can hold without further hypotheses. 

\noindent
{\bf An outline of the proof of Theorem 1:} Without loss of generality, let $B_{2} \times [0, 4]$ be the space-time region in Theorem 1, where $B_{2}$ is the open ball in ${\mathbb R}^{2}$ (space) with radius $2$ and center at the origin. Let $V_{t}$, $t \in [0, 4]$ denote the moving 1-varifold at time $t$, and let $J$ denote a fixed stationary triple junction with vertex at the origin. Thus $J$ consists of three half-lines meeting at 120 degree angles at the origin. By assumption, the flow is weakly close to $J$, which in particular means that the space-time $L^{2}$ distance (height excess) $\mu$ of the flow $\{V_{t}\}_{t \in [0, 4]}$ relative to $J$, defined by $$\mu = \left(\int_{0}^{4}\int_{B_{2}} {\rm dist}^{2} \, (x, J) \, d\|V_{t}\|(x) dt \right)^{1/2},$$  is small.  

As mentioned before, our proof of Theorem 1 is based on a parabolic version of the blow-up method. We first use the full 
strength of \cite{Kasai-Tonegawa} to obtain (in Proposition~\ref{smprop}) a graphical representation of the varifolds $V_{t}$, with an appropriate estimate, away from the center of $J$.  We  use this graphical representation to establish various a priori  space-time and time uniform $L^{2}$-estimates that control the behavior of the flow  in the region near the center of $J$. In particular, a key step is to show that $\mu$ does not concentrate near the center of $J$. 

Our approach to establishing this a priori non-concentration estimate is inspired by the basic strategy developed by Simon \cite{Simon2} for minimal submanifolds. A key ingredient in Simon's method is the monotonicity formula for minimal submanifolds, whose role here is played by (a certain  local estimate inspired by) the Huisken monotonicity formula.
In the present parabolic setting, there are  several interesting new aspects also. These stem from firstly the fact that all we have at our 
disposal is Brakke's inequality defining the flow---which a priori only tells us something about the rate of change of mass (length) and not much about the velocity of motion---and secondly the fact that we need a number of nontrivial preliminary estimates 
involving curvature, which in the case of minimal submanifolds are not needed (regardless of the dimension). 
A key such estimate (established in Proposition~\ref{hde}) gives an interior space-time $L^{2}$ bound for the generalized curvature $h = h(V_{t}, x)$ of $V_{t}$ (where $x \in {\rm spt} \, \|V_{t}||$) in terms of $\mu$ whenever $\mu$ is sufficiently small; said more precisely, 
\begin{equation*}
\int_{1}^{3}\int_{B_{3/2}}|h|^2\,
d\|V_t\|dt\leq c\mu^2
\end{equation*}
provided $\mu$ is sufficiently small, where $c$ is a fixed constant independent of the flow. We use this estimate and computations similar to those used in the derivation of the Huisken monotonicity formula \cite{Huisken} (see also \cite{Ilmanenp}) to  
establish (in Proposition~\ref{pade}), whenever $\mu$ is sufficiently small,  that 
\begin{equation*}
\int_{5/4}^{s}\int_{B_{1}}\left|h+\frac{x^{\perp}}{2(s-t)}\right|^2 \rho_{(0,s)}(x,t)\, d{\|V_{t}\|}(x) dt \leq c\mu^{2}
\end{equation*}
for any $s \in [3/2, 3]$ such that $h(V_{s}, \cdot) \in L^{2}_{\rm loc} \, (\|V_{s}\|)$ and $\Theta \, (\|V_{s}\|, 0) \geq \Theta \, (\|J\|, 0) = 3/2$, where $c$ is a fixed constant independent of the flow, $\rho_{(0,s)}(x,t) = (4\pi(s-t))^{-1/2}e^{\frac{-|x|^{2}}{4(s-t)}}$ ($-\infty < t < s < \infty$)  is the backwards heat kernel
with  pole at $(0,s)$ and $\Theta \, (\|V\|, Z)$ denotes the density of $V$ at $Z$. 
This bound is then used (in Proposition~\ref{tildest}) to obtain,  for any $s \in [3/2, 3]$ as above and any $\kappa \in (0, 1)$, the crucial estimate 
\begin{equation*}
\sup_{t \in [5/4, s)} (s-t)^{-\kappa}\int_{B_{3/4}}\rho_{(0, s)}(\cdot, t) {\rm dist}^{2} \, (\cdot, J) \, d\|V_{t}\|(x) \leq c_{0}\mu^{2},
\end{equation*}
again provided $\mu$ is sufficiently small, where $c_{0}$ depends only on $\kappa$. This says that the $L^2$ distance of $V_{t}$ from the triple junction $J$ weighted by the backwards heat kernel decays quickly in time; in particular,  this estimate implies that the contribution to $\mu^{2}$ coming from a small spatial neighborhood of the origin and a slightly smaller time interval is a small proportion of $\mu^{2}$.

Using these estimates, we  carry out a careful blow-up analysis in Sec.~\ref{blowup-analysis}. We emphasize that the term $(s-t)^{-\kappa}$ appearing in the preceding estimate, though not needed for the non-concentration conclusion just pointed out, plays an important role in the blow up analysis. Once the appropriate asymptotic decay for the blow-ups are established, we obtain a space-time excess improvement lemma (Lemma~\ref{blowprop7}) for the flow, the iteration of which leads to Theorem 1 in a fairly standard way.

\noindent
{\bf Organization of the paper:} In Sec.\,2 we fix notation and state our main results. In Sec.\,3
we use results of \cite{Kasai-Tonegawa} to give a graph representation of the moving curves away from the center of the triple junction. The 
main result in Sec.\,4 is Proposition~\ref{hde}, which gives a time-uniform estimate on the difference of length
between the moving curve and the triple junction in terms of the space-time $L^2$ distance of the flow to the triple junction. The 
same estimate gives an $L^2$ curvature estimate in terms of the $L^2$ distance. Sec.\,5 contains
the main non-concentration estimate, Proposition~\ref{tildest}, which shows that the $L^2$ distance  does not concentrate
around the junction point. This is used in Sec.\,6 to estimate the location of and the H\"{o}lder norm (in time) for the junction points in term of the $L^2$ distance. All of these estimates are used to carry out a blow-up argument in Sec.\,7 on each of the three rays of the triple 
junction and to show that the three pieces of the blow-up come together at a single point in a regular fashion. Sec.\,8 describes the iteration procedure giving a H\"{o}lder estimate
of the gradient up to the (moving) junction points, proving the main local regularity theorem 
(Theorem~\ref{mainreg}).
Sec.\,9 contains the proof of the partial regularity theorem (Theorem~\ref{mainpa}).  
Sec.\,10 contains a further result concerning the nature of the tangent flows at singular points of a flow satisfying the hypotheses of Theorem~\ref{mainpa}. 
\section{Notation, background  and the main theorems}
\subsection{Basic notation}
 Let ${\mathbb N}$ be the the set of natural numbers and let ${\mathbb R}^+=\{x\geq 0\}$. For $r\in (0,\infty)$
 and $a\in {\mathbb R}^2$, define $B_r(a)=\{x\in {\mathbb R}^2 : |x-a|<a\}$ and when $a=0$, define $B_r=
 B_r(0)$. 
 We write ${\mathcal L}^1$ for the Lebesgue measure on ${\mathbb R}$ and 
 ${\mathcal H}^1$ for the 1-dimensional Hausdorff measure on ${\mathbb R}^2$. The restriction
 of ${\mathcal H}^1$ to a set $A$ is denoted by ${\mathcal H}^1\res_A$. 
 For an open set $U\subset {\mathbb R}^2$ let $C_c(U)$ be the set of all compactly supported 
 continuous functions on $U$ and let $C_c(U;{\mathbb R}^2)$ be the the set of all compactly supported
 continuous vector fields. The index $k$ of $C^k$ indicates continuous $k$-th order differentiability. 
 $\nabla$ always indicates differentiation with respect to the space variable. For $g\in C_c^1(U;{\mathbb R}^2)$,
 $\nabla g(x)$ is a $2\times 2$ matrix-valued function. 
 
 For any Radon measure $\lambda$ on ${\mathbb R}^2$ and $\phi\in C_c({\mathbb R}^2),$ we 
 shall often write $\lambda(\phi)$ for $\int_{\mathbb R^2} \phi\, d\lambda$. We let ${\rm spt}\,\lambda$ 
 be the support of $\lambda,$ and $\Theta(\lambda,x)$ be the 1-dimensional density of $\lambda$
 at $x$, i.e., $\Theta(\lambda, x) = \lim_{r\rightarrow0^{+}}\lambda(B_r(x))/(2r)$, when the limit exists. For a $\lambda$ measurable
 function $f$, $f\in L^2(\lambda)$ means $\int_{{\mathbb R}^2} |f|^2\, d\lambda<\infty$ and $f \in L^2_{\rm loc}(\lambda)$ means $\int_{K}|f|^{2}  \, d\lambda< \infty$ for each compact set $K \subset {\mathbb R}^{2}$. 
 
 For $-\infty<t<s<\infty$ and $x,y\in {\mathbb R}^2$, define the 1-dimensional backwards heat kernel
$\rho_{(y, s)}$ by
 \begin{equation}
 \rho_{(y,s)}(x,t)=\frac{1}{\sqrt{4\pi(s-t)}}\exp\left(-\frac{|x-y|^2}{4(s-t)}\right).
 \label{defheat}
 \end{equation}
 Let ${\bf G}(2,1)$ be the space of 1-dimensional subspaces of ${\mathbb R}^2$. For $S\in {\bf G}(2,1)$,
 we identify $S$ with orthogonal projection (and the $2\times 2$ matrix associated with orthogonal projection) of ${\mathbb R}^{2}$ onto $S$ and we let $S^{\perp}\in {\bf G}(2,1)$ be the orthogonal complement of $S$ in ${\mathbb R}^{2}$.
 For $A,B\in {\rm Hom}({\mathbb R}^2;{\mathbb R}^2)$ let $A\cdot B={\rm trace}\,(A^{\star}\circ B)$, where $\circ$ denotes composition and $A^{\star}$ is the transpose of $A$. 
 Let $u\otimes  v\in {\rm Hom}({\mathbb R}^2;
 {\mathbb R}^2)$ be the tensor product of $u,v\in {\mathbb R}^2$. 
 \subsection{Varifolds}
 We next recall the notion of varifolds and some related definitions. For a detailed discussion on varifolds, see 
 \cite{Allard,Simon}. For an open set $U\subset {\mathbb R}^2$,
 define $G_1(U)=U\times {\bf G}(2,1)$. A {\it 1-varifold}  in $U$ is a Radon measure on $G_1(U).$
 The set of $1$-varifolds in $U$  is denoted by ${\bf V}_1(U)$. Varifold convergence is the usual measure convergence on $G_1(U)$. 
 In this paper we are only concerned with 1-varifolds and shall often just refer to them as varifolds subsequently. 
 For a varifold $V\in {\bf V}_1(U)$ let $\|V\|$ denote the weight measure associated to $V$, defined by 
 $\|V\|(\phi) =\int_{G_1(U)}
 \phi(x)\, dV(x,S)$ for $\phi\in C_c(U)$. Given an ${\mathcal H}^1$ measurable countably 1-rectifiable
 set $M\subset U$ with locally finite ${\mathcal H}^1$ measure, there is a natural varifold denoted by $|M|$ and defined by
 \begin{equation*}
 |M|(\phi)=\int_M \phi(x,{\rm Tan}_x\, M)\, d{\mathcal H}^1(x), \ \ \forall \phi\in C_c(G_1(U)),
 \end{equation*}
 where ${\rm Tan}_x\, M\in {\bf G}(2,1)$ is the approximate tangent space of $M$ at
 $x$ which exists ${\mathcal H}^1$ a.e$.$ on $M$.  $V\in {\bf V}_1(U)$ is called {\it integral} if
 \begin{equation*}
 V(\phi)=\int_M \phi(x,{\rm Tan}_x\, M)\theta(x)\, d{\mathcal H}^1(x), \ \ \forall \phi\in C_c(G_1(U)),
 \end{equation*}
 for some ${\mathcal H}^1$ measurable countably 1-rectifiable set $M\subset U$ and ${\mathcal H}^1$
 a.e$.$ integer-valued locally ${\mathcal H}^1$ integrable function $\theta$ defined on $M$. 
 The function $\theta$ is called the {\it multiplicity} of $V$. 
 Set of all integral 1-varifolds in $U$ is denoted by ${\bf IV}_1(U)$. $V$ is called a {\it unit density} 
 varifold if $V$ is integral with $\theta=1$ a.e., that is, if $V=|M|$ with $M$ as above. For $V \in {\bf IV}_{1}(U)$ and a given $\|V\|$-measurable vector field $g$ on $U$, we shall often write $\int_U (g(x))^{\perp}\, d\|V\|(x)$
 for $\int_{G_1(U)} S^{\perp}(g(x))\, dV(x,S)$. 
 
 For $V\in {\bf V}_1(U),$ let $\delta V$ denote the first variation of $V$, defined by
 \begin{equation*}
 \delta V(g) =\int_{G_1(U)} \nabla g(x)\cdot S\, dV(x,S) \;\;\; \ \ \forall g\in C^1_c(U;{\mathbb R}^2).
 \end{equation*}
 Let $\|\delta V\|$ be the total variation measure of $\delta \, V$ when it exists (which is the case precisely when $\delta \, V$ is locally bounded). If $\|\delta V\|$ is absolutely 
 continuous with respect to $\|V\|$, we have for some $\|V\|$ measurable vector field $h(V,\cdot),$
 \begin{equation}
 \delta V(g)=-\int_U g(x)\cdot h(V,x)\, d\|V\|(x), \ \ \forall g\in C^1_c(U;{\mathbb R}^2).
 \label{fvf}
 \end{equation}
By $h(V,\cdot)$, we always mean the vector field satisfying \eqref{fvf}. We simply call $h(V,\cdot)$ 
the (generalized) curvature of $V$. 
If $V=|M|$ and $M$ is a $C^2$ curve, then $h(V,\cdot)$ is the curvature of $M$  
times the unit normal vector. 

For any $V\in {\bf IV}_1(U)$ with locally bounded first variation, we note that Brakke's perpendicularity
theorem \cite[Ch.\,5]{Brakke} says that
\begin{equation}
\int_U (g(x))^{\perp} \cdot h(V,x)\, d\|V\|(x)=\int_U g(x)\cdot h(V,x)\, d\|V\|(x), \ \ \forall g\in C_c(U;{\mathbb R}^2).
\label{perpthm}
\end{equation}
\subsection{Definition of Brakke flow}
\label{defBF}
To motivate the definition of Brakke flow, suppose that we have a family of classical planar networks 
$\{M_t\}_{t\in I}$ for some open interval $I\subset {\mathbb R}$, 
i.e., each $M_t$ consists of a locally finite union of embedded $C^{\infty}$
open curves contained in some fixed open set $U\subset {\mathbb R}^2$ 
which meet smoothly at  triple junctions of 120 degree angles and which are moving smoothly in time. 
Let $v$ be the normal velocity vector of $M_t$. For any test function $\phi=\phi(x,t)$ with support ($\phi(\cdot,t)$)
a compact subset of $U$, computation of the
first variation of length gives
\begin{equation}
\label{firstlen}
\frac{d}{dt}\int_{M_t}\phi(x,t)\, d\mathcal H^1(x) = \int_{M_t} (-\phi h+\nabla\phi)\cdot v+\frac{\partial \phi}{\partial t}\, d\mathcal H^1(x). 
\end{equation}
Here, $h$ is the curvature vector of $M_t$ of open curves.  It is important to note that the 120 degree angle condition at triple 
junctions helps to cancel out the contribution of the first variations at boundary points of open curves. Conversely, if some normal 
vector field $\tilde v$ on $M_t$ satisfies the following {\em inequality}
\begin{equation}
\label{firstlen2}
\frac{d}{dt}\int_{M_t}\phi(x,t)\, d\mathcal H^1(x) \leq  \int_{M_t} (-\phi h+\nabla\phi)\cdot {\tilde v}+\frac{\partial \phi}{\partial t}\, d\mathcal H^1(x)
\end{equation}
for all {\em non-negative} test function $\phi=\phi(x,t)$ with support ($\phi(\cdot,t)$) a compact subset of $U$
then one can prove that $\tilde v$ is necessarily equal to 
the normal velocity $v$. 
In particular, \eqref{firstlen2} with $\tilde v=h$
is satisfied if and only if $\{M_t\}$ is a curvature flow. 
Since quantities in \eqref{firstlen2} naturally extend to the varifold setting, 
this characterization of normal velocity allows us to formulate a weak notion of curvature flow  of varifolds, which we call the Brakke flow. 
\begin{define}
Let $U\subseteq \mathbb R^2$ be open and $I\subseteq\mathbb R$ be an interval (i.e.\ a connected open 
subset of $\mathbb R$). 
A family of 1-varifold $\{V_t\}_{t\in I}$ is a Brakke flow if 
\begin{itemize}
\item[(1)] $V_t\in {\bf IV}_1(U)$ for a.e.\ $t\in I$,
\item[(2)] $h(V_t,\cdot)$ exists and belongs to $L^2_{loc}(\|V_t\|)$ for a.e.\ $t\in I$ and 
\item[(3)]
for each $\phi\in C^1(U\times I; \mathbb R^+)$ with $\phi(\cdot, t)\in C^1_c(U)$ for $t\in I$, and for any $t_1,t_2\in I$ with $t_1<t_2$,
\begin{equation}
\|V_t\|(\phi(\cdot,t))\Big|_{t=t_1}^{t_2}\leq \int_{t_1}^{t_2}\int_{U} (-\phi h(V_t,\cdot)+\nabla\phi)\cdot h(V_t,\cdot)
+\frac{\partial\phi}{\partial t}\, d\| V_t\|dt.
\label{firstlen3}
\end{equation}
\end{itemize}
\end{define}
The correspondence between $M_t$ above and $V_t$ is that $V_t=|M_t|$, $\|V_t\|=\mathcal H^1\res_{M_t}$, and $h(V_t,\cdot)=h$. 
The inequality \eqref{firstlen3} is an integrated formulation of \eqref{firstlen2} in time. The requirement that $h(V_t,\cdot)\in L^2_{loc}(\|V_t\|)$ 
exists for a.e.\ $t$ presupposes $\|\delta V_t\|$ is  locally bounded and absolutely continuous
with respect to $\|V_t\|$. Thus, if there is a smooth junction of three curves, they need to meet at 120 degree angles. 
\subsection{The right-hand side of the curvature flow equation}
To define a flow more general than the Brakke flow, for
$V\in {\bf V}_1(U)$, a vector field $u\in L^2_{loc}(\|V\|)$ and $\phi\in C^1_c(U;{\mathbb R}^+)$, define ${\mathcal B}(V, u, \phi)$ as follows:
\begin{equation}
{\mathcal B}(V,u,\phi) =\int_U (-\phi(x)h(V,x)+\nabla\phi(x))\cdot(h(V,x)+(u(x))^{\perp})\, d\|V\|(x)
\label{defB}
\end{equation}
if $V\in {\bf IV}_1(U)$, $\|\delta V\|$ exists  and is absolutely continuous with respect to
$\|V\|$ and $h(V,\cdot)\in L^2_{loc}(\|V\|)$; otherwise ${\mathcal B}(V,u,\phi) =-\infty$. 
If $\{M_t\}$ is a family of smoothly embedded curves moving with normal velocity $v=h+u^{\perp}$, where $u$ is a 
given smooth ambient vector field, then one can prove just as Sec. \ref{defBF} that 
for all $\phi\in C^1_c(U\times I;\mathbb R^+)$,
\begin{equation}
\frac{d}{dt}\int_{M_t}\phi(x,t)\, d\mathcal H^1(x)=\mathcal B(|M_t|,u(\cdot,t),\phi(\cdot,t))+\int_{M_t}\frac{\partial\phi}{\partial t}(\cdot,t)\, d\mathcal H^1(x).
\label{defB2}
\end{equation}
Conversely, having the property \eqref{defB2} with $\leq$ in place of equality for smooth $M_{t}$ implies
that $v=h+u^{\perp}$ on $M_t$, so we may use \eqref{defB2} with 
$\leq$ in place of equality as a weak formulation of the condition $v=h+u^{\perp}$. (See  Hypothesis 
(A4) below.) 
\label{defBFG}
\subsection{The triple junction $J$, some test functions and norms}
Let $J\subset{\mathbb R}^2$ be defined by
\begin{equation}
J =\left\{(s,0)\,:\, s\geq 0\right\}\cup \left\{(-s,\sqrt{3} \,s)\,:\, s\geq 0\right\}\cup\left\{(-s,-\sqrt{3}\, s)\,:\, s\geq 0\right\}.
\label{regtridef}
\end{equation}
Because of the 120 degree angles of intersection, $|J|$ is a stationary varifold. Being stationary, $|J|$ is also 
a Brakke flow which is static in time. We often think of $J$ as a set in $\mathbb R^2$ but also identify
(with a slight abuse of notation) $J$ with a static multiplicity 1 
Brakke flow in space-time, as in the statement of Theorem 1.

For $\theta\in{\mathbb R}$ denote by ${\bf R}_{\theta}:{\mathbb R}^2\rightarrow {\mathbb R}^2$ the map 
corresponding to the counterclockwise orthogonal rotation by angle $\theta$. 
Define a similarity class of $J$ by 
\begin{equation}
{\mathcal J}=\left\{{\bf R}_{\theta}(J)+\xi\, :\, \theta\in \left(-\frac{\pi}{3},\frac{\pi}{3}\right)\mbox{ and }\xi\in {\mathbb R}^2\right\},
\label{defj}
\end{equation}
and for $R\in (0,\infty)$,
define $d_R:{\mathcal J}\times{\mathcal J}\rightarrow{\mathbb R}$ by
\begin{equation}
d_R (J_1,J_2)=\max\left\{R^{-1}|\xi_1-\xi_2|, |\theta_1-\theta_2|\right\}
\label{hje1}
\end{equation}
for $J_1={\bf R}_{\theta_1}(J)+\xi_1$ and $J_2={\bf R}_{\theta_2}(J)+\xi_2$. It is clear that $d_R$ defines a metric on ${\mathcal J}$. We set
\begin{equation}
d(J_1,J_2)=d_1(J_1,J_2).
\label{hje1s}
\end{equation}
Let $\hat\phi:{\mathbb R}^2\rightarrow{\mathbb R}$ 
be a smooth radially symmetric non-negative function such that $\hat\phi(x)=1$ for
$|x|\leq \frac14$, $\hat\phi(x) =0$ for $|x|\geq \frac12$, $0\leq \hat\phi\leq 1$ and $|\nabla\hat \phi|\leq 8$. Set 
\begin{equation}
\Cl[c]{c-p} =\int_{\mathbb R}\hat\phi(s,0)\, ds.
\label{c-pdef}
\end{equation}
Define
\begin{equation}
\phi_j(x)=\hat\phi\left({\bf R}_{-\frac{2(j-1)\pi}{3}}(x)-(1,0)\right)
\label{phij}
\end{equation}
for $j=1,2,3$ and $x\in {\mathbb R}^2$.
Note that $\phi_1$ is just $\hat\phi$ composed with translation by $(1,0)$, and $\phi_2$ and $\phi_3$ 
are $\phi_{1}$ composed with orthogonal rotation by
$\frac{2\pi}{3}$ and $\frac{4\pi}{3}$ respectively.  We have $\int_J \phi_j\, d{\mathcal H}^1
=\Cr{c-p}$ for each $j=1,2,3$. 
For $J'={\bf R}_{\theta}(J)+\xi \in {\mathcal J}$, $R\in (0,\infty)$ and $j=1,2,3$, define
\begin{equation}
\phi_{j,J',R}(x)=\phi_j({\bf R}_{-\theta}(R^{-1}(x-\xi))), \ \ \phi_{j,J'}(x) =\phi_{j,J',1}(x)
\label{hje2}
\end{equation}
for $x\in {\mathbb R}^2$, where $\phi_j$ is as in \eqref{phij}. Note that $\phi_{j,J}=\phi_j$. 

For $R > 0$, let $Q_{R} =\{(s,t)\in (-R,R)\times(-R^2,R^2)\}$. We shall use the following norm for functions $f \,: \, Q_{R} \to {\mathbb R}$:
\begin{equation*}
\begin{split}
\|f\|_{C^{1,\zeta}(Q_{R})}=&\sup_{(s,t)\in Q_{R}}(R^{-1}|f(s,t)| +|\nabla f(s,t)|)\\
&+\sup_{(s_1,t_2),(s_2,t_2)\in Q_{R}, (s_1,t_2) \neq (s_2,t_2)}
\frac{R^{\zeta}|\nabla f (s_1,t_1)-\nabla f (s_2,t_2)|}{\max\{|s_1-s_2|,|t_1-t_2|^{\frac{1}{2}}\}^{\zeta}} \\
&+\sup_{(s,t_1),(s,t_2)\in Q_{R}, t_{1} \neq t_{2}} \frac{R^{\zeta}|f(s,t_1)-f(s,t_2)|}{|t_1-t_2|^{\frac{1+\zeta}{2}}}.
\end{split}
\end{equation*}
 Note that this norm is 
invariant under the parabolic change of variables in the sense that if $\tilde f(\tilde s, \tilde t)=R^{-1}f(s,t)$ where $\tilde s=R^{-1} s$ and $\tilde t=R^{-2}t$, then $\|\tilde f\|_{C^{1, \zeta}(Q_{1})} = \|f\|_{C^{1, \zeta}(Q_{R})}$. We shall denote by $C^{1,\zeta}(Q_{R})$ the space of functions $f \, : \, Q_{R} \to {\mathbb R}$ with   $\|f\|_{C^{1, \zeta}} (Q_{R}) < \infty$. 
\subsection{Hypotheses and the main theorems}\label{hypotheses}
Let $U\subseteq {\mathbb R}^2$ be open and let $I \subseteq {\mathbb R}$  be an interval (i.e.\ a connected open subset of ${\mathbb R}$). Assume

\begin{itemize}
\item[(A0)] $p \in [2, \infty)$ and $q \in (2, \infty)$ are fixed numbers such that 
\begin{equation}
\zeta \equiv 1-\frac{1}{p}-\frac{2}{q}>0.
\label{power}
\end{equation}
\end{itemize}
For each $t \in I$, let $V_{t}$ be a 1-varifold in $U$ and  
$u(\cdot, t) \, : \, U \to {\mathbb R}^{2}$ a $\|V_{t}\|$ measurable vector field such that: 

\begin{itemize}
\item[(A1)]  $V_t\in {\bf IV}_1(U)$ for a.e.\ $t\in I$;
\item[(A2)] there exists $E_{1} \in [1, \infty)$ such that for each $B_r(x)\subset U$ and each $t\in I$,
\begin{equation}
\|V_t\|(B_r(x))\leq 2r E_1;
\label{ddd}
\end{equation}

\item[(A3)] $u$ satisfies
\begin{equation}
\left(\int_{I}\left(\int_U |u(x,t)|^p\, d\|V_t\|(x)\right)^{\frac{q}{p}}\, dt\right)^{\frac{1}{q}}<\infty;
\label{powerfin}
\end{equation}
\item[(A4)] for each $\phi\in C^1(U\times I;{\mathbb R}^+)$ with $\phi(\cdot,t)\in C_c^1(U)$ for $t \in I,$ and for any $t_{1}, t_{2} \in I$ with $t_1<t_2$, 
\begin{equation}
\label{meq}
\|V_{t}\|(\phi(\cdot,t))\Big|_{t=t_1}^{t_2}\leq \int_{t_1}^{t_2}{\mathcal B}(V_t,u(\cdot,t),\phi(\cdot,t))\, dt
+\int_{t_1}^{t_2}\int_U \frac{\partial\phi}{\partial t}(\cdot,t)\, d\|V_t\|dt.
\end{equation}
\end{itemize}
{\bf Remark:} 
If $u=0$, \eqref{meq} reduces to \eqref{firstlen3} and (A3) is irrelevant. \eqref{meq} also requires ${\mathcal B}
(V_t,0,\phi(\cdot,t))>-\infty$ for a.e.\ $t\in I$ hence $h(V_t,\cdot)$ exists and in $L^2_{loc}(\|V_t\|)$. Thus in this case, (A1), (A3) and (A4)
are equivalent to $\{V_t\}$ being a Brakke flow.
Moreover, Huisken's monotonicity formula which can be derived from \eqref{firstlen3} 
shows (see e.g.\ \cite{Ilmanenp}) that
(A2) is locally satisfied for some $E_1$ so that Brakke flows always satisfy (A1)-(A4) locally. 
(A4) is a weak formulation of the condition $v=h+u^{\perp}$, 
as discussed in Sec. \ref{defBFG}.

The following is our $\varepsilon$-regularity theorem, whose proof takes up a 
major part of the paper:

\begin{thm}\label{mainreg}
Corresponding to $p$, $q$ as in $({\rm A}0)$, $E_{1} \in [1, \infty)$ and $\nu\in (0,1)$, there exist $\Cl[eps]{e-8}\in (0,1)$ and $\Cl[c]{c-14}\in (1,\infty)$ such that
the following holds: For $R\in (0,\infty)$ and $U=B_{4R}$, let $\{V_t\}_{t\in [-2R^2,2R^2]}$ and
$\{u(\cdot,t)\}_{t\in [-2R^2,2R^2]}$ satisfy (A1)-(A4) with $I = [-2R^{2}, 2R^{2}]$. Suppose 
\begin{equation}
\mu \equiv \Big( R^{-5}\int_{-2R^2}^{2R^2} \int_{B_{4R}} {\rm dist}\,(\cdot,J)^2 \, d\|V_t\|dt\Big)^{\frac12}<\Cr{e-8};
\label{mr1}
\end{equation}
there exist $j_{1}, j_{2} \in \{1, 2, 3\}$ such that 
\begin{equation}
R^{-1}\|V_{-2R^2}\|(\phi_{j_1,J,R})\leq (2-\nu)\Cr{c-p}, \ \ R^{-1}\|V_{2R^2}\|(\phi_{j_2,J,R})\geq \nu\Cr{c-p}
\label{mr2}
\end{equation}
where $J$, $\Cr{c-p}$ and $\phi_{j,J,R}$ are as defined in \eqref{regtridef}, \eqref{c-pdef} and \eqref{hje2} respectively, and 
\begin{equation}
\|u\| \equiv R^{\zeta} \Big(\int_{-2R^2}^{2R^2}\big(\int_{B_{4R}} |u|^p\, d\|V_t\|\big)^{\frac{q}{p}}\Big)^{\frac{1}{q}}
<\Cr{e-8}.
\label{mr3}
\end{equation}
Then there exists $\hat a\in C^{\frac{1+\zeta}{2}}([-R^2,R^2];B_R)$
and,  letting 
\begin{equation}
l_j(t)=\mbox{first coordinate of }{\bf R}_{-\frac{2\pi(j-1)}{3}}
\left(\hat a(t)\right)\ \ and \ \ 
D_j =\cup_{t\in [-R^2,R^2]} \big([l_j(t),R]\times\{t\}\big),
\label{mr4}
\end{equation}
there exist functions $f_j\in C^{1,\zeta}(D_j),$ $j=1,2,3,$
such that for all $t\in [-R^2,R^2]$, we have
\begin{equation}
\label{mr6}
{\bf R}_{-\frac{2\pi(j-1)}{3}}({\rm spt}\,\|V_t\|)\cap \big([l_j(t),R]\times
[-R,R]\big)=\{(x,f_j(x,t)): x\in [l_j(t),R]\}
\end{equation}
for $j=1,2,3$ and 
\begin{equation}
\frac{\partial f_1}{\partial x}(l_1(t),t)=\frac{\partial f_2}{\partial x}(l_2(t),t)
=\frac{\partial f_3}{\partial x}(l_3(t),t).
\label{mr6.5}
\end{equation}
Furthermore, we have
\begin{equation}
\label{mr7}
\|\hat a\|_{C^{\frac{1+\zeta}{2}}([-R^2,R^2];B_R)}+\sum_{j=1}^3
\|f_j\|_{C^{1,\zeta}(D_j)}\leq \Cr{c-14}\max\{\mu,\|u\|\}.
\end{equation}
\end{thm}

\noindent
{\bf Remark:} In case $u\in C^{\beta}(B_{4R}\times [-2R^2,2R^2];{\mathbb R}^2)$ 
(where H\"{o}lder 
continuity is in the usual parabolic sense),   the result of
\cite{Tonegawa} shows that $f_j$ are $C^{2,\beta}$ away from 
the junction point and that the flow satisfies $v=h+u^{\perp}$ in the classical
sense. In fact in this case, up-to-the-junction-point $C^{2,\beta}$ regularity of $f_j$ as well as 
$C^{1,\frac{\beta}{2}}$ regularity of $\hat a$  also hold, and can be proved 
by the well-know reflection technique of \cite{KNS} combined with the
regularity theory of linear parabolic systems \cite{Solo}. If $u$ is smooth, 
or zero in particular, then $\hat a$ is smooth, and $f_j$ are smooth up to the junction point. Note that
for the reflection technique of \cite{KNS}, having $C^{1,\zeta}$ regularity given by our theorem  
provides  the crucial starting hypothesis. 

\vspace{.2in}

We note that all quantities are scale invariant under parabolic change of variables
so we may and we shall, without loss of generality,  set $R=1$ in the proof of Theorem~\ref{mainreg}. The inequality \eqref{mr1} provides a closeness
to $J$ of $\|V_t\|$ in the $L^2$ distance, and \eqref{mr2} requires some closeness to $J$
in terms of measure at the initial and final times. The latter also prevents complete loss
of measure $\|V_t\|$ during this time interval. Assumption \eqref{mr3} ensures smallness
of the perturbation from the $u=0$ case, and obviously is not needed if $u=0$. 
The conclusion is that each time slice of the flow in a smaller space-time domain consists of  three embedded curves meeting precisely at one common junction point, that this junction point $\hat{a}(t)$ at time $t$ is a
H\"{o}lder continuous function of $t$, and that the three curves are represented as
$C^{1,\zeta}$ graphs up to the junction point as in \eqref{mr4} and \eqref{mr6}, satisfying  the estimate \eqref{mr7}. Moreover, 
at each time, the three curves meet at 120 degree angles, as expressed by \eqref{mr6.5}. 
If more regularity is assumed on $u$, then we have, as stated in
the Remark above, better regularity up to the junction point.

By combining Theorem \ref{mainreg} with a stratification theorem of White \cite{White0}, we obtain 
the following partial regularity theorem. We make an assumption on the density of static
tangent flows, which is equivalent to assuming that any static tangent flow is either a unit
density line or a unit density triple junction. For 
the precise definition of tangent flow, 
see Sec.\,\ref{secpart}. Tangent flows are analogous to tangent cones for minimal submanifolds:
at each point in space-time, at least one tangent flow is obtained by passing to a subsequential varifold limit of parabolic rescalings of the flow, and tangent flows enjoy a nice
homogeneity property called backwards-cone-like (see Sec.\,\ref{secpart} (b)). 
\begin{thm}
In addition to the assumptions ($A1$)-($A4$), assume that 
\begin{itemize}
\item[($A5$)] at each point in space-time, whenever a tangent flow to $\{V_{t}\}_{t \in I}$ is static, the density at the
origin of the tangent flow is strictly less than 2. 
\end{itemize}
Then there exists a closed set $\Sigma_1\subset U\times I$ with the parabolic Hausdorff 
dimension at most 1 such that the flow in $U \times I \setminus \Sigma_{1}$ is classical in the sense that for any $(x,t)\in U\times I \setminus \Sigma_1$, there exists 
a space-time neighborhood $U_{x,t}$ containing $(x,t)$ such that $U_{x,t}\cap \cup_{t'}({\rm spt}\,\|V_{t'}\|\times\{t'\})$ is
either empty, a $C^{1,\zeta}$ graph over a line segment or a $C^{1,\zeta}$ triple junction as described in Theorem~\ref{mainreg}.
\label{mainpa}
\end{thm}

\noindent
{\bf Remarks:} {\bf (1)} Under hypotheses (A1)-(A4), we may completely classify all non-trivial static tangent flows. 
They are time-independent stationary integral
varifolds whose supports are unions of half-lines emanating from the origin.
Thus hypothesis (A5) requires that any static tangent flow is a single line or
a triple junction, either one with unit density, and nothing else. 
This hypothesis is motivated by the fact that any static tangent flow with density greater than or equal to 2 should be unstable for various physical models. For the motion of grain boundaries, one
observes that junctions with more than 3 edges appear and break up instantaneously. 
Mathematically, any junction (including lines) with multiplicity strictly greater than 1 is not
mass minimizing in the sense that one can always set the multiplicity equal to 1 and reduce
the mass. Any unit density junction with more than 3 edges may be mapped by a 
suitable Lipschitz function so that the image of the map has locally 
less ${\mathcal H}^1$ measure {\it as a set}.
(Note that the usual varifold push-forward counts multiplicities of the image and the mapping
here is different from it.) It is called {\it reduced mass model} according to Brakke \cite[p.57]{Brakke}.

\noindent
{\bf (2)} Other than the singularities coming from collisions of triple junctions, 
we may also have some curve disappearing suddenly. Such singularities are included
in the closed set $\Sigma_1$. 

\noindent
{\bf (3)} The parabolic Hausdorff dimension counts the time variable as 2. 
Thus, $\Sigma_{1}$ having parabolic dimension at most 1 implies that the times at which singularities can occur form a closed subset of $I$ of usual 
Hausdorff dimension at most $\frac{1}{2}.$  

\noindent
{\bf (4)} The short-time existence of classical network flows (i.e.\ those consisting of curves meeting smoothly and only at a locally finite number of triple junctions) was established by Bronsard--Reitich \cite{Bronsard} when the initial network itself is classical, and has recently been extended to more general initial networks satisfying certain regularity and non-degeneracy assumptions by Ilmanen--Neves--Schulze \cite{INS}. (The work \cite{INS} also gives a result that says that a flow weakly close to a triple junction $J$ is $C^{1, \alpha}$ close to $J$ in the interior (see \cite{INS}, Theorem~1.3 and remark (iii)) in the special case when the flow is a priori assumed to be regular, using methods limited to such an priori regularity hypothesis.) It remains an interesting open problem to prove a general existence theorem for curvature flows
satisfying (A1)-(A5). 

\noindent
{\bf (5)} See Sec.\,\ref{difpat} for a more detailed characterization of $\Sigma_1$
in terms of tangent flows. 

\section{A graph representation away from the singularity of $J$}
 We apply results from \cite{Kasai-Tonegawa}
to show that the supports of the moving varifolds, in the region outside a small neighborhood of the singularity of $J$, are represented as a $C^{1,\zeta}$ graphs, with the $C^{1,\zeta}$
norm bounded in terms of the $L^2$ distance of the flow to $J$. This result will be used frequently in the rest of the paper. 
\begin{prop}
Corresponding to $\tau\in (0,\frac12),\nu\in (0,1),E_1 \in [1, \infty)$ and $p \in [2, \infty)$, $q \in (2, \infty)$ satisfying (\ref{power}), there exist $\Cl[eps]{e-p}\in (0,1)$ and 
$\Cl[c]{c-p-1}\in (1,\infty)$ such that the following holds: Suppose that $\{V_t\}_{t\in [0,4]}$ and $\{u(\cdot,t)\}_{t\in [0,4]}$ satisfy (A1)-(A4) with $U = B_{2}$ and $I = [0,4],$ and that
\begin{equation}
\mu \equiv \left(\int_0^{4}\int_{B_{2}}{\rm dist}(\cdot,J)^2\, d\|V_t\|dt\right)^{\frac12}\leq \Cr{e-p},
\label{smep1}
\end{equation}
\begin{equation}
\|u\| \equiv \left(\int_0^{4}\left(\int_{B_{2}} |u|^p\, d\|V_t\|\right)^{\frac{q}{p}}dt\right)^{\frac{1}{q}}\leq \Cr{e-p},
\label{smep2}
\end{equation}
\begin{equation}
\|V_0\|(\phi_{j_1})\leq (2-\nu)\Cr{c-p} \ \ \mbox{for some} \ \ j_{1} \in \{1, 2, 3\} \ \ \mbox{and}, 
\label{smep3}
\end{equation}
\begin{equation}
 \|V_{4}\|(\phi_{j_2})\geq \nu \Cr{c-p} \ \ \mbox{for some} \ \  j_2\in \{1,2,3\}.
\label{smep4}
\end{equation}
Then 
\begin{equation}
B_{2-\frac{\tau}{5}}\cap\left\{{\rm dist}(\cdot,J)>\frac{\tau}{5}\right\}\cap {\rm spt}\,\|V_t\|=\emptyset \ \ \ \mbox{for each} \ \  t\in [\tau,4]
\label{smep5}
\end{equation}
and there exist $f_j: [\tau ,2-\tau]\times [\tau,4-\tau]\rightarrow {\mathbb R}$, $j=1,2,3$, such that
$f_1,f_2,f_3$ are  $C^{1,\zeta}$ in space and $C^{\frac{1+\zeta}{2}}$ in time with 
\begin{equation}
\begin{split}
\|f_j\|_{C^{1,\zeta}(Q)}
\leq \Cr{c-p-1}\max\{\mu,\|u\|\}
\end{split}
\label{smep6}
\end{equation}
where $Q=[\tau,2-\tau]\times [\tau,4-\tau]$, 
and
\begin{equation}
([\tau ,2-\tau]\times [-\tau,\tau])\cap \big({\bf R}_{-\frac{2(j-1)\pi}{3}}({\rm spt}\,\|V_t\|)\big)=
\{(s,f_j(s,t))\,:\, s\in [\tau ,2-\tau]\},
\label{smep7}
\end{equation}
for $j=1,2,3$ and for all $t\in [\tau,4-\tau]$. 
Furthermore, for a.e$.$ $t\in [\tau,4-\tau]$, we have that
\begin{equation}
\label{exden}
\Theta(\|V_t\|,x) \in \{1, 3/2\} \;\; \mbox{for each} \;\; x\in B_{2-\tau}\cap {\rm spt}\,\|V_t\|.
\end{equation}
\label{smprop}
\end{prop}
{\it Proof}. 
The claim \eqref{smep5} may be deduced by applying \cite[Prop.\,6.4 \& Cor.\,6.3]{Kasai-Tonegawa}.
We also see that for fixed $\tau \in (0, 1/2)$, ${\rm spt}\|V_t\|\cap B_{2-\tau}$ approaches $J$ as $\mu,\,\|u\|\rightarrow 0$ for all $t\in [\tau,4]$. 

To prove the existence of a graph representation
and the $C^{1,\zeta}$ estimate as asserted, assume for fixed $E_{1}$, $\nu$, $\tau$ that
the claim is false. Then for each $m\in {\mathbb N}$
there exist $\{V_t^{(m)}\}_{t\in [0,4]}$ and $\{u^{(m)}(\cdot,t)\}_{t\in [0,4]}$ satisfying (A1)-(A4) and \eqref{smep1}-\eqref{smep4} with $U = B_{2}$, $I=[0,4]$, $\Cr{e-p} = m^{-1}$ and with $V_{t}^{(m)}, u^{(m)}$ in place of $V_{t}, u$, but there are no 
functions $f_{1}$, $f_{2}$, $f_{3}$ with the stated regularity satisfying \eqref{smep6} and \eqref{smep7} with $\Cr{c-p-1}=m,$ $u^{(m)}$ in place of $u$ and $\mu^{(m)} = \left(\int_{0}^{4}\int_{B_{2}} {\rm dist}\, (\cdot, J)^{2} \, d\|V_{t}^{(m)}\| dt \right)^{\frac{1}{2}}$ in place of $\mu$. 
To obtain a contradiction, we will use \cite[Th.\,8.7]{Kasai-Tonegawa} which shows the existence of such a graph 
representation as in the asserted conclusion under a set of hypotheses. In order to check that the hypotheses of \cite[Th.\,8.7]{Kasai-Tonegawa}, with $V^{(m)}_{t},$ $u^{(m)}$ in place of $V_{t}$, $u$, are satisfied for sufficiently large $m$, we first prove that
\begin{equation}
\|V_t^{(m)}\|\rightarrow {\mathcal H}^1\res_{J}
\label{smep8}
\end{equation}
as $m\rightarrow\infty$ on $B_{2}$ for all $t\in (0,4)$. To see this, take any $\varphi\in C_c^2(B_{2};{\mathbb R}^+)$ and use
(A4) to obtain for any $t_1,t_2\in [0,4]$ with $t_1<t_2$ that
\begin{equation}
\|V^{(m)}_{t_2}\|(\varphi)-\|V_{t_1}^{(m)}\|(\varphi)\leq \int_{t_1}^{t_2}\int -\frac{|h^{(m)}|^2}{2}\varphi
+\frac{|\nabla\varphi|^2}{2\varphi}+|u^{(m)}|^2\varphi+|\nabla\varphi||u^{(m)}|\, d\|V_t^{(m)}\|dt
\label{smep9}
\end{equation}
where $h^{(m)}=h(V_t^{(m)},x)$. Define $\Phi^{(m)}(t)=\int_0^t\int\frac{|\nabla\varphi|^2}{2\varphi}+|u^{(m)}|^2\varphi+|\nabla\varphi||u^{(m)}|\, d\|V_t^{(m)}\|dt$. 
Since $\sup |\nabla\varphi|^2/\varphi\leq \sup 2|\nabla^2\varphi|$, we see by H\"{o}lder's inequality combined with the fact that $\|u^{(m)}\| \leq m^{-1}$ and $\|V_{t}^{(m)}\|({\rm spt} \, \varphi) \leq 4E_{1}$, that $\{\Phi^{(m)}\}_{m\in {\mathbb N}}$
is bounded uniformly in the $\frac{q-2}{q}$-H\"{o}lder norm on $[0,4]$. In particular, we may choose a uniformly convergent subsequence of $\{\Phi^{(m)}\}$. 
We also see from \eqref{smep9} that $\|V_t^{(m)}\|(\varphi)-\Phi^{(m)}(t)$ is monotone decreasing in $t$. Because of this, we may 
extract a subsequence $\{m_j\}_{j\in {\mathbb N}}$ such that, for each $t$ except for a countable number of $t$, 
$\{\|V_t^{(m_j)}\|(\varphi)\}_{j\in {\mathbb N}}$ is a Cauchy sequence. Next choose a countable
set $\{\varphi_l\}_{l\in {\mathbb N}}\subset C_c^2(B_{2};{\mathbb R}^+)$ such that it is dense with respect to the $C^0$ topology in $C_c^2(B_{2};
{\mathbb R}^+)$. By a diagonal argument, we may choose a subsequence so that $\{\|V_t^{(m_j)}\|(\varphi_l)\}_{j\in {\mathbb N}}$ is a Cauchy 
sequence except for a countably many $t$ and for all $l\in {\mathbb N}$. 
Since $\{\|V_t^{(m_j)}\|\}_{j\in {\mathbb N}}$ is a set of uniformly bounded measures, this shows
that $\|V_t^{(m_j)}\|$ converges to a Radon measure, say, $\lambda_t$ for all $t\in [0,4]$ except for countably many $t$. By the weak 
compactness of Radon measures, we may extend $\lambda_{t}$ to all $t \in [0, 4]$ as Radon measures such that passing to a further subsequence, $\|V_t^{(m_j)}\|$ converges to $\lambda_t$ 
for all $t\in [0,4]$. By the first part of the proof, we know that ${\rm spt}\, \lambda_t\subset J$ for all $t\in (0,4]$. By \eqref{smep3} and
\eqref{smep4}, we also have that
\begin{equation}
\lambda_0(\phi_{j_1})\leq (2-\nu) \Cr{c-p}, \ \lambda_4(\phi_{j_2})\geq \nu \Cr{c-p}.
\label{smep10}
\end{equation}
for some $j_1,j_2\in \{1,2,3\}$.

Next note that by \eqref{smep9}, for any fixed $\varphi\in C^2_c(B_{2};{\mathbb R}^+)$, we have
\begin{equation}
\int_0^4 \int\frac{|h(V_t^{(m)},\cdot)|^2}{2}\varphi\, d\|V_t^{(m)}\|dt
\leq \|V_0^{(m)}\|(\varphi)+\Phi^{(m)}(4)
\label{smep11}
\end{equation}
where the right-hand side is uniformly bounded. By Fatou's lemma applied to \eqref{smep11},
for a.e$.$ $t\in[0,4]$, there exists a (time-dependent) subsequence such that 
\begin{equation}
\int |h(V_t^{(m_{j_k})},\cdot)|^2\varphi\, d\|V_t^{(m_{j_k})}\|\leq C(t,\varphi, E_1)
\label{smep12} 
\end{equation}
for all $k\in {\mathbb N}$. By (A1), for a.e$.$ $t\in [0,4]$, $V_t^{(m_{j_k})}$ is integral. 
By the compactness theorem for integral varifolds \cite{Allard}, there exists a further subsequence 
which converges to a limit varifold $\tilde V _t$ which is integral with locally
$L^2$ generalized curvature. Since ${\rm spt}\, \|\tilde V_{t}\|\subset J$ and $h(\tilde V_t,\cdot)\in L^2(\|\tilde V_t\|)$, the density function
of $\|\tilde V_t\|$ must be constant on each line segment of $J \cap B_{2}$. In fact, the density must be
 constant on all of $J \cap B_{2}$. 
Since $\|\tilde V_t\|=\lambda_t$, we conclude that for a.e.\ $t \in [0, 4]$, $\lambda_t=\theta(t){\mathcal H}^{1}\res_{J \cap B_{2}}$ for some 
$\theta(t)\in \{0\}\cup {\mathbb N}$. On the other hand, using $\varphi = \phi_j$ in \eqref{smep9}
and taking a limit, we obtain for any $0\leq t_1<t_2\leq 4$
\begin{equation}
\lambda_{t_2}(\phi_j)-\lambda_{t_1}(\phi_j)\leq (t_2-t_1)E_1 \sup\frac{|\nabla\phi_j|^2}{2\phi_j}.
\label{smep13}
\end{equation}
In general, we have $\lambda_t(\phi_j)=\theta(t)\Cr{c-p}$, but by \eqref{smep10} and \eqref{smep13}, 
we see that $\lambda_t(\phi_j)=\Cr{c-p}$ for a.e$.$ $t\in [0,4]$. 
Thus we have shown that $\theta(t)=1$ for a.e.\  $t\in [0,4]$. Since \eqref{smep13} holds also for any function in
$C_c^2(B_{2};{\mathbb R}^+)$ in place of $\phi_{j}$, we see that $\lambda_t={\mathcal H}^1\res_{(J\cap B_{2})}$ for all $t\in (0,4)$.
Since the limit measure is uniquely determined, the whole sequence converges, and this establishes 
\eqref{smep8}. 

We are now ready to use \cite[Th.\,8.7]{Kasai-Tonegawa}. For any $(x_0,t_0)\in
J\cap(B_{2}\setminus B_{\tau})\times (\tau,4-\tau)$, consider a small domain containing $(x_{0}, t_{0})$ which is a positive
distance away from the origin. In view of \eqref{smep8}, we have (8.85) and (8.86) (with $\nu=1/2$, for
example) of \cite[Th.\,8.7]{Kasai-Tonegawa} satisfied for all sufficiently large $m$, 
and so are the other assumptions (8.83) and (8.84). Note that \cite[Th.\,8.7]{Kasai-Tonegawa}
assumes the varifold has unit density a.e.\, but one can indeed prove using a variant of Huisken's monotonicity formula
(see \cite[Prop.\,6.2]{Kasai-Tonegawa})
and \eqref{smep8} that there cannot be a point $ (x, t) \in B_{2-\tau}\times(\tau,4-\tau)$ with $\Theta(\|V_{t}^{(m)}\|,x)\geq 2$
for all sufficiently large $m$. 
Thus, near $(x_0,t_0)$, ${\rm spt}\, \|V^{(m)}_t\|$ is represented as a graph a function of the desired regularity
satisfying the estimate \eqref{smep6} on this domain. By covering $Q$ with a finite number of such
small domains, we obtain \eqref{smep6}  with a suitable constant $\Cr{c-p-1}$ which
depends only on $\tau,\nu,E_1, p,q$. 

To see \eqref{exden}, first note that we have $h(V_t,\cdot)\in 
L^2_{loc}(\|V_t\|)$ and that $V_t$ is integral for a.e$.$ $t$; thus, for such $t$,  
$\Theta(\|V_t\|,x)$ exists and is greater than or equal to $1$ for all $x\in {\rm spt}\,\|V_t\|\cap B_2.$ Moreover, at each $x\in {\rm spt}\,\|V_t\|\cap B_2$, there exists a tangent cone which is
 a stationary 1-dimensional integral varifold, and thus we may conclude that $\Theta(\|V_t\|,x)$ is an
 integer multiple of $1/2$. Again using a variant of Huisken's monotonicity formula and \eqref{smep8},
 we conclude that $\Theta(\|V_t\|,x)\leq 3/2$ for
 $x\in B_{2-\tau}$ for sufficiently small $\Cr{e-p}$. 
 This proves \eqref{exden}. 
\hfill{$\Box$}
\section{A priori estimates I: the space-time $L^2$-curvature estimate}
\label{curv-est}
The estimates in this section are analogous to those in \cite[Sec.\,5]{Kasai-Tonegawa}
where, roughly speaking, one obtains a time-uniform estimate for the difference of the $k$-dimensional area of the moving
varifolds and that of a flat plane. There are a few subtle differences however. First, unlike in \cite[Prop.\,4.6]{Kasai-Tonegawa}, we do not have a useful
Lipschitz graph approximation for the varifolds near the triple junction. Here, in Proposition 
\ref{firprop}, we rely on  the fact that $L^2$ control 
of curvature gives good $C^{1,\frac12}$ norm control since the varifolds are 1-dimensional,  with the end result being similar to \cite[Prop.\,5.2]{Kasai-Tonegawa}. Lemma \ref{ode} is 
similar in spirit to \cite[Lem.\,5.5]{Kasai-Tonegawa}, although since we do not have, unlike in \cite[Sec.\,6]{Kasai-Tonegawa}, an $L^{\infty}$
estimate for ${\rm dist}\,(\cdot,J)$ on ${\rm spt} \, \|V_{t}\|$ in terms of the $L^2$-distance, we need to employ a different estimate in Lemma~\ref{ode}. 
Once this lemma is established, we obtain \eqref{thap9} and
\eqref{thap10} just as we obtain the corresponding estimates in \cite[Th.\,5.7]{Kasai-Tonegawa}.
\begin{define}
Let $\phi_{\rm rad}:{\mathbb R}^2 \rightarrow{\mathbb R}^+$ be a non-negative radially symmetric
function such that $\phi_{\rm rad}=1$ on $B_1$, $|\nabla\phi_{\rm rad}|\leq 4$ and $
\phi_{\rm rad}\in C^{\infty}_{c}(B_{\frac32})$. Define
\begin{equation}
{\bf c}=\int_{J}\phi_{\rm rad}^2\, d{\mathcal H}^1.
\label{cfi1}
\end{equation}
\end{define}
Proposition \ref{firprop} and Corollary \ref{fircor} below do not involve the time variable. 
\begin{prop}
Corresponding to $E_1\in [1,\infty)$ there exist constants $\Cl[c]{c-1}\in (1,\infty)$, $\Cl[al]{alpha-1},\, \Cl[be]{beta-1}\in (0,1)$ such that
the following holds. For $V\in {\bf IV}_1(B_{2})$ with $h(V,\cdot)\in L^2(\|V\|)$, 
define
\begin{equation}
\hat\alpha=\left(\int_{B_{2}}|h(V,x)|^2\phi_{\rm rad}(x)^2\, d\|V\|(x)\right)^{\frac12},
\label{fir3.2}
\end{equation}
and
\begin{equation}
\hat\mu=\left(\int_{B_{2}}{\rm dist}\, (x,J)^2\, d\|V\|(x)\right)^{\frac12}.
\label{fir4}
\end{equation}
Assume the following \eqref{firm1}-\eqref{fir3}: 
\begin{equation}
\|V\|(B_{2})\leq 4E_1;
\label{firm1}
\end{equation}
\begin{equation}
\Theta(\|V\|,x) \in \left\{1, \frac32\right\} \,\,\mbox{for each} \,\, x\in {\rm spt}\, \|V\|;
\label{fir0}
\end{equation} 
\begin{equation}
\left(\left[\frac12, \frac32\right]\times \left[-\frac12,\frac12\right]\right)\cap
 \left({\bf R}_{-\frac{2(j-1)\pi}{3}}\big({\rm spt}\, \|V\|\big)\right)=\left\{(s, f_j(s))\,:\, |s-1|\leq\frac12\right\}
\label{fir1}
\end{equation}
for $j=1,2,3$ and for some $f_j\in C^{1}(\{s\in {\mathbb R}\,:\, |s-1|\leq \frac12\})$ with
\begin{equation}
\hat\beta \equiv \max_{j\in \{1,2,3\}}\|f_j\|_{C^{1}(\{|s-1|\leq \frac12\})}\leq \Cr{beta-1};
\label{fir2}
\end{equation}
\begin{equation}
{\rm spt}\, \|V\|\cap B_{\frac{19}{10}}\cap \left\{x\in{\mathbb R}^2\,:\, {\rm dist}\, (x,J)\geq \frac{1}{10} \right\}=\emptyset
\label{fir6}
\end{equation}
and
\begin{equation}
\hat\alpha\leq \Cr{alpha-1}.
\label{fir3}
\end{equation}
Then there exist three $C^{1,\frac12}$ curves $l_1,\, l_2,\, l_3$ 
having one common end point near the origin in $B_{1}$ meeting at
$120$ degree angles such that ${\rm spt} \, \|V\| \cap B_{1}=\cup_{j=1}^3 l_j\cap B_{1}$. 
Moreover, we have
\begin{equation}
\big|{\mathcal H}^1({\rm spt}\, \|V\|\cap B_1)-3\big|\leq \Cr{c-1}(\hat\alpha\hat\mu+\hat\beta^2)
\label{fir5}
\end{equation}
and
\begin{equation}
\big|\|V\|(\phi_{\rm rad}^2) -{\bf c}\big| \leq \Cr{c-1}(\hat\alpha\hat\mu+\hat\beta^2).
\label{fir5.5}
\end{equation}
\label{firprop}
\end{prop}
{\it Proof}. If $\Theta(\|V\|,x)=1$, the Allard regularity theorem \cite{Allard} combined with the fact that $h(V,\cdot)\in
L^2(\|V\|)$ shows that ${\rm spt}\, \|V\|$
is an embedded $C^{1,\frac12}$ curve in some neighborhood of $x$. A standard argument
shows that the set of points with $\Theta(\|V\|,x)=\frac32$ is discrete. Specifically,
assume for a contradiction that this set has an accumulation point $a\in B_{2}$ and let $a_{1}, a_{2}, \ldots$ be such that $a_{j} \neq a$, $\Theta \, (\|V\|,a_{j}) = \frac32$ for each $j=1, 2, \ldots$ and  
$a_{j} \rightarrow a$ as $j \to \infty$. By (\ref{fir0}), $\Theta \, (\|V\|, a) = \frac32$.  Consider a subsequential limit ${\mathbf C}$ of rescalings 
of $V$ about $a$ by the scale factors $|a_i-a|^{-1}$. This limit ${\mathbf C}$ is a stationary integral cone which has density $=\frac32$ at the origin and at 
$b = \lim_{i\rightarrow\infty}\frac{a_i-a}{|a_i-a|} \in {\mathbf S}^{1}$, hence along the whole ray determined by $b$, producing a positive measure portion of ${\rm spt} \, \|{\mathbf C}\|$ on which the density is equal to $\frac32$, a
contradiction to  integrality of ${\mathbf C}$. In particular this shows that
there are only finitely many points $a_1,\cdots, a_N$ in $B_1$ with $\Theta \, (\|V\|, a_{j}) = \frac32$. Away from these points, one may
parametrize the curves by the arc length parameters. One can prove that the weak second derivative 
is precisely $h$ a.e., thus the $\frac12$-H\"{o}lder norms of the unit tangent vectors along the curves 
can be estimated by $L^2$ norm of $h$. Because of this fact, at each $a_1,\cdots,a_N$, 
${\rm spt}\,\|V\|$ consists of three emanating $C^{1,\frac12}$ curves, and due to the stationarity of
tangent cone, the meeting angles of these three curves are all $120$ degrees. Let us call these points
``junction point''.
On $B_1\setminus\{a_1,\cdots,a_N\}$, ${\rm spt}\, \|V\|$ consists of 
$C^{1,\frac12}$ curves, with the end points on $\partial B_1$ or at junction point
or without any end points (i.e. closed curve). Since any closed curve $l\subset B_1$ have $\int_{l}|h|\geq 2\pi$,
H\"{o}lder's inequality and \eqref{firm1} show $\hat\alpha\geq 2\pi/ \sqrt{4E_1}$. 
Thus for small $\Cr{alpha-1}$, there cannot be 
any closed curve in $B_1$. By \eqref{fir1} and \eqref{fir6}, there are at least three curves $l_1,\, l_2,\, l_3$ 
which are the extensions of curves represented as graphs of $f_j,\, f_2,\, f_3$, respectively. As noted 
already, small $\Cr{alpha-1}$
implies small change of unit tangent vector along the curves. Thus choosing sufficiently small $\Cr{alpha-1}$ and
$\Cr{beta-1}$, we may assume that $l_1,\, l_2,\, l_3$ are very close to straight lines and close to $J$ in $C^1$ norm. They hit one of 
junction point, or otherwise \eqref{fir6} would be violated. We note that there cannot be more than one
triple junction. Suppose otherwise. Suppose one follows the parametrization of $l_1$ from the right and one hits
 the first junction point $a_1$. Then one follows the curve
emanating from $a_1$ by turning $60$ degrees ``to the right''. Suppose that one hits another junction point $a_2$
along this curve. Then one follows the curve just like before by turning $60$ degrees. By choosing $\beta_1$ and
$\alpha_1$ small, we can make sure that the latter curve has the tangent vector whose 
direction is only at most, say, $1$ degree 
different from that of $(1,\sqrt{3})$. Note that $l_1$ is very close to $x$-axis, so after turning $60$ degrees twice,
the curve should be almost parallel to such vector. Unless one hits the next junction point $a_3$ along this 
curve, we would have a contradiction to \eqref{fir6}. In this manner, one can argue that there would have to be 
infinitely many junction points in $B_1$, noting that the sum of total variations of tangent vector for each 
curve can be made arbitrarily small. This contradicts the finiteness of the number of triple junctions. 
Thus we may conclude that $l_1,\, l_2,\, l_3$ meet
at a unique junction point close to the origin under appropriate restrictions on $\Cr{alpha-1}$ and $\Cr{beta-1}$. 
This proves the first part of the claim. For the rest of the proof, continue to denote the graph representations of
$l_1,\, {\bf R}_{-\frac{2\pi}{3}} l_2,\, {\bf R}_{-\frac{4\pi}{3}} l_3$ by (for $j=1,2,3$, respectively)
\begin{equation}
\{(s,f_j(s))\in {\mathbb R}^2\,:\, s\in [s_j,1]\},
\label{fir7}
\end{equation}
where we note that 
\begin{equation}
(s_1,f_1(s_1))={\bf R}_{\frac{2\pi}{3}}(s_2,f_2(s_2))={\bf R}_{\frac{4\pi}{3}}(s_3,f_3(s_3))
\label{fir8}
\end{equation}
is the meeting point of triple junction of three curves. 
It then follows from \eqref{fir8} and ${\bf R}_0+{\bf R}_{-\frac{2\pi}{3}}
+{\bf R}_{-\frac{4\pi}{3}}={\bf 0}$ that
\begin{equation}
s_1+s_2+s_3=0
\label{fir9}
\end{equation}
and
\begin{equation}
 f_1(s_1)+f_2(s_2)+f_3(s_3)=0.
\label{fir9.5}
\end{equation}
The fact that the curves meet at $120$ degrees implies that
\begin{equation}
f'_1(s_1)=f'_2(s_2)=f'_3(s_3),
\label{fir10}
\end{equation}
where $f_j'$ here means the right derivative of $f_j$. 
Recalling $|h|=|f''_j|/(1+|f'_j|^2)^{\frac32}$, \eqref{fir3.2} and \eqref{fir2}, we have
\begin{equation}
\sup_{s\in [s_j,1]}|f'_j(s)|\leq |f'_j(1)|+\int_{s_j}^1|f''_j(s)|\, ds\leq \hat\beta+2\hat\alpha
\label{fir10.5}
\end{equation}
for all sufficiently small $\Cr{alpha-1}$ and $\Cr{beta-1}$.
For each $j=1,2,3$, we next compute ${\mathcal H}^1(l_j\cap B_1)$.
\begin{equation}
\big|{\mathcal H}^1(l_j\cap B_1)-(1-s_j)\big|\leq \hat\beta^2 +\int_{s_j}^1(\sqrt{1+|f'_j(s)|^2}-1)\, ds
\leq \hat\beta^2+\frac12\int_{s_j}^1|f'_j(s)|^2\, ds,
\label{fir11}
\end{equation}
where the first inequality is due to the error estimate using \eqref{fir2} near $\partial B_1$ 
and the second inequality is by $\sqrt{1+t^2}-1\leq t^2/2$.
Summing over $j$ and using \eqref{fir9}, we obtain
\begin{equation}
\big|\sum_{j=1}^3 {\mathcal H}^1(l_j\cap B_1)-3\big|\leq 
 3\hat\beta^2+\frac12\sum_{j=1}^3\int_{s_j}^1|f'_j(s)|^2\, ds.
\label{fir12}
\end{equation}
By integration by parts, and using \eqref{fir9.5}, \eqref{fir10} and \eqref{fir2}, we have
\begin{equation}
\sum_{j=1}^3\int_{s_j}^{1}|f'_j|^2=\sum_{j=1}^3\big( f_j(1)f'_j(1)-\int_{s_j}^{1}
f_j f''_j\big)\leq 3\hat\beta^2+2\hat\alpha\big(\sum_{j=1}^3\int_{s_j}^{1} |f_j|^2\big)^{\frac12}.
\label{fir13}
\end{equation}
For each $j$, we note that $|f_j(s)|={\rm dist}\,((s,f_j(s)),J)$ away from the origin and close to
$J$. More precisely, the equality holds when $(s,f_j(s))$ and positive $x$-axis 
has angle $\leq \frac{\pi}{3}$. When $(s,f_j(s))$ and negative $x$-axis has angle
$\leq \frac{\pi}{6}$, then we have $|f_j(s)|\leq {\rm dist}\, ((s,f_j(s)),J)$. 
For $s\in (\tilde{s}_1,\tilde{s}_2)$, suppose $(s, f_j(s))$ lies in the 
sector $\{tv\, :\, t>0,\, |v|=1,\, \cos\frac{\pi}{3}>v\cdot (1,0)>\cos\frac{5\pi}{6}\}$, where
we have $|f_j(s)|>{\rm dist}\, ((s,f_j(s)),J)$. Denote $\tilde{\delta}=
\tilde{s}_2-\tilde{s}_1$ and note that $\tilde{\delta}$ is small when $\Cr{alpha-1}$ and $\Cr{beta-1}$ are 
small. Considering the geometry of graph, for sufficiently small $\sup |f'_j|$, we have
\begin{equation}
\tilde{\delta}\sup_{s\in (\tilde{s}_1,\tilde{s}_2)}|f_j(s)|^2\leq 2\tilde{\delta}\inf_{s\in (\tilde{s}_2,\tilde{s}_2+\tilde{\delta})}|f_j(s)|^2
\leq 2\int_{l_j\cap B_1}{\rm dist}\, (x,J)^2\, d\|V\|(x).
\label{fir13.1}
\end{equation}
Outside of $(\tilde{s}_1,\tilde{s}_2)$, as stated, we have $|f_j(s)|\leq {\rm dist}\,((s,f_j(s)),J)$, thus we have
by \eqref{fir13.1}
\begin{equation}
\int_{s_j}^1|f_j|^2\leq 3\int_{l_j\cap B_1}{\rm dist}\, (x,J)^2\, d\|V\|(x).
\label{fir14}
\end{equation}
Combining \eqref{fir12}, \eqref{fir13} and \eqref{fir14}, we obtain \eqref{fir5}. To obtain \eqref{fir5.5}, 
note that $\phi_{\rm rad}=1$ on $B_1$ and thus we need to be concerned with region of integration 
over $B_{\frac32}\setminus B_1$ of $\|V\|$. But we have \eqref{fir2}, thus the difference of
integrations in this region can be estimated by $c\hat\beta^2$. In the estimate one uses the radial 
symmetry of $\phi_{\rm rad}$ to obtain the quadratic estimate. Thus we obtain \eqref{fir5.5} with some 
suitable constant $\Cr{c-1}$. 
\hfill{$\Box$}
\begin{cor}
For a given $E_1\in [1,\infty)$, let $\Cr{c-1},\Cr{alpha-1},\Cr{beta-1}$ be the
corresponding constants obtained in Proposition 
\ref{firprop}. For $V\in {\bf IV}_1(B_{2})$ with $h(V,\cdot)\in L^2(\|V\|)$, define
$\hat\alpha,\hat\beta, \hat\mu$ as \eqref{fir3.2}, \eqref{fir2} and \eqref{fir4}. 
Assume \eqref{firm1}-\eqref{fir6}.
Define
\begin{equation}
\hat{E}=\|V\|(\phi_{\rm rad}^2)-{\bf c},
\label{hatE}
\end{equation}
and assume that
\begin{equation}
3\Cr{c-1}\hat\beta^2\leq |\hat{E}|.
\label{cor1}
\end{equation}
Then we have
\begin{equation}
\hat\alpha^2\geq \min\{\Cr{alpha-1}^2,(2\Cr{c-1}\hat \mu)^{-2}|\hat{E}|^2\}.
\label{cor2}
\end{equation}
\label{fircor}
\end{cor}
{\it Proof}. If $\hat\alpha\geq \Cr{alpha-1}$ holds, then \eqref{cor2} holds and there is nothing further to prove. 
Thus consider the case $\hat\alpha< \Cr{alpha-1}$. Since \eqref{firm1}-\eqref{fir6} are assumed, 
we fulfill all the assumptions of Proposition \ref{firprop}, thus we have \eqref{fir5.5}. Using the
notation of \eqref{hatE}, this implies that we have either $|\hat{E}|\leq 2\Cr{c-1} \hat\alpha\hat\mu$, 
or $|\hat{E}|\leq 2\Cr{c-1}\hat\beta^2$. The last possibility is excluded by \eqref{cor1}.
Thus we have $\hat\alpha^2\geq (2\Cr{c-1}\hat\mu)^{-2}|\hat{E}|^2$. 
Thus we have either $\hat\alpha^2\geq \Cr{alpha-1}^2$ or the last possibility. 
This proves \eqref{cor2}.
\hfill{$\Box$}

The next ODE lemma connects Corollary \ref{fircor} to (A4). 
\begin{lemma}
Corresponding to $P,T\in (0,\infty)$ there exist $\Cl[c]{ha},\Cl[c]{ha2}\in (0,\infty)$ 
such that the following holds: Given a non-negative function $g \in L^2([0,T])$ and a monotone
decreasing function $\Phi\,:\, [0,T]\rightarrow{\mathbb R}$, define
$f\,:\, [0,T]\rightarrow {\mathbb R}^+$ by
\begin{equation}
f(t)=P\min\left\{1,g(t)^{-2} |\Phi(t)|^2\right\}
\label{le1}
\end{equation}
when $g(t)>0$ and $f(t)=P$ when $g(t)=0,$ 
and suppose that  
\begin{equation}
\Phi(t_2)-\Phi(t_1)\leq -\int_{t_1}^{t_2}f(t)\, dt,\hspace{.2cm}
0\leq \forall t_1<\forall t_2\leq T.
\label{le2}
\end{equation}
Then 
\begin{itemize}
\item[(1)] if $\Phi(0)\leq \Cr{ha}$, then 
\begin{equation}
\Phi(T)\leq \Cr{ha2} \|g\|_{L^{2}([0, T])}^{2}.
\label{le3}
\end{equation}
\item[(2)] if $\Phi(T)\geq -\Cr{ha}$, then
\begin{equation}
\Phi(0)\geq - \Cr{ha2}\|g\|_{L^{2}([0, T])}^{2}.
\label{le4}
\end{equation}
\end{itemize}
\label{ode}
\end{lemma}
{\it Proof}.
We prove (1) first. Set
\begin{equation}
\Cr{ha}=\frac{PT}{8},\ \ \ \Cr{ha2}=\frac{8}{PT^2}.
\label{kys}
\end{equation}
We may assume $\Phi(t)>0$ for 
all $t\in [0,T]$ since $\Phi(T)\leq 0$ otherwise and \eqref{le3} is trivially true. 
Assume for a contradiction that \eqref{le3} were false. 
Set 
\begin{equation}
c=\|g\|_{L^2}\sqrt{2/T}
\label{ky0}
\end{equation}
and define $A_1=\{t\in [0,T]\, :\, g(t)\geq c\}$
and $A_2=[0,T]\setminus A_1$. It is easy to check that ${\mathcal L}^1(A_1)
\leq T/2$, and thus 
\begin{equation}
{\mathcal L}^1 (A_2)\geq T/2.
\label{ky1}
\end{equation}
We next define 
$A_{2,a}=\{t\in A_2\, :\, g(t)^{-2}|\Phi(t)|^2\geq 1\}$ and $A_{2,b}=A_2\setminus A_{2,a}$.
By \eqref{ky1}, we have either ${\mathcal L}^1(A_{2,a})\geq T/4$ or ${\mathcal L}^1(A_{2,b})\geq T/4$.
In the first case, since $f(t)=P$ on $A_{2,a}$, we have 
\begin{equation}
\Phi(T)-\Phi(0)\leq -\int_{A_{2,a}} f(t)\, dt=-P{\mathcal L}^1(A_{2,a})\leq -\frac{PT}{4}.
\label{ky2}
\end{equation}
We have $\Phi(T)-\Phi(0)\geq -\Cr{ha}$ and \eqref{ky2} gives a contradiction
to \eqref{kys}. 
In the second case, using $f(t)=P g(t)^{-2}\Phi(t)^2$ on $A_{2,b}$ and integrating $(-\Phi(t)^{-1})'
\leq -Pg(t)^{-2}$ we have
\begin{equation}
-\Phi(T)^{-1}+\Phi(0)^{-1}\leq -P\int_{A_{2,b}}\frac{dt}{g(t)^2}
\leq -P c^{-2} {\mathcal L}^1 (A_{2,b})\leq -\frac{PT^2}{8} \|g\|_{L^2}^{-2}.
\label{ky2.5}
\end{equation}
We used $g(t)\leq c$ on $A_{2,b}\subset A_2$ and \eqref{ky0}. 
Since we are assuming \eqref{le3} is false, by \eqref{kys}, we have
\begin{equation}
-\Phi(T)^{-1}> -\Cr{ha2}^{-1} \|g\|_{L^2}^{-2}=-\frac{PT^2}{8} \|g\|_{L^2}^{-2}.
\label{ky3}
\end{equation}
Two inequalities \eqref{ky2.5} and \eqref{ky3} give a contradiction. This 
proves (1). 
For (2), replace $\Phi(\cdot)$ by $-\Phi(T-\cdot)$ and $f(\cdot)$ by $f(T-\cdot)$,
and then apply the previous argument. If $-\Phi(T)\leq \Cr{ha}$, then one concludes that
$-\Phi(0)\leq \Cr{ha2}\|g\|_{L^2}^2$. Thus we obtain the result of (2).
\hfill{$\Box$}
\begin{prop}
Corresponding to $\nu,E_1,p,q$ there exist $\Cl[eps]{e-1}\in (0,1)$
and $\Cl[c]{c-2}\in (1,\infty)$ with the following property. 
Suppose $\{V_t\}_{t\in [0,4]}$ and $\{u(\cdot, t)\}_{t\in [0,4]}$
satisfy (A1)-(A4) on $B_{2}\times [0,4]$. Assume \eqref{smep1} and \eqref{smep2}
with $\Cr{e-p}$ there replaced by $\Cr{e-1}$, \eqref{smep3} and \eqref{smep4}.
Then we have
\begin{equation}
\sup_{t\in [1,3]}\big|\|V_t\|(\phi_{\rm rad}^2)-{\bf c}\big|
\leq \Cr{c-2}\max\{\mu,\|u\|\}^2
\label{thap9}
\end{equation}
and 
\begin{equation}
\int_{1}^{3}\int_{B_2}|h(V_t,\cdot)|^2\phi_{\rm rad}^2\,
d\|V_t\|dt\leq \Cr{c-2}\max\{\mu,\|u\|\}^2.
\label{thap10}
\end{equation}
\label{hde}
\end{prop}
{\it Proof}. 
We first use Proposition \ref{smprop} with $\tau=\frac14$ to obtain $\Cr{e-p}$ and $\Cr{c-p-1}$ so that we have
\eqref{smep5}-\eqref{exden} with $\tau=\frac14$ there. We also use Proposition \ref{firprop} to obtain $\Cr{c-1},
\Cr{alpha-1},\Cr{beta-1}$ corresponding to $E_1$. Then, for a.e$.$ $t\in [\frac{1}{2},\frac{7}{2}]$, 
we have conditions \eqref{firm1}-\eqref{fir1} and \eqref{fir6} satisfied for $V=V_t$ there. If we further assume that 
\begin{equation}
\Cr{c-p-1}\max\{\mu,\|u\|\}\leq \Cr{beta-1},
\label{bb1}
\end{equation}
then \eqref{fir2} is also satisfied due to \eqref{smep6}. 
We restrict $\Cr{e-1}$ so that
\eqref{bb1} holds by the following: 
\begin{equation}
\Cr{e-1}\leq \min\{\Cr{e-p}, \Cr{beta-1} \Cr{c-p-1}^{-1}\}.
\label{asu1}
\end{equation}
Next fix $P$ and $T$ as 
\begin{equation}
P=\frac{1}{16}\min\{\Cr{alpha-1}^2,(2\Cr{c-1})^{-2}\}, \ \ T=\frac12.
\label{th11.1}
\end{equation}
With these choices of $P$ and $T$, we obtain $\Cr{ha}$ and
$\Cr{ha2}$ by Lemma \ref{ode}.
With $\Cr{ha}$ fixed, we choose a small $\tau$ and then restrict $\Cr{e-1}$ so that, by using Proposition \ref{smprop}, we have
\begin{equation}
\|V_{\frac12}\|(\phi_{\rm rad}^2)\leq {\bf c}+\frac{\Cr{ha}}{2}
\label{va1}
\end{equation}
and
\begin{equation}
\|V_{\frac72}\|(\phi_{\rm rad}^2)\geq {\bf c}-\frac{\Cr{ha}}{2}.
\label{va2}
\end{equation}
We will fix $\Cr{c-2}$ later.
We also set 
\begin{equation}
\beta_*=\max_{j=1,2,3} \sup_{t\in [\frac{1}{2},\frac{7}{2}]} \|f_j(\cdot, t)\|_{C^1(\{|s-1|\leq \frac12\})}(\leq \Cr{c-p-1}\max\{\mu,\|u\|\}\leq \Cr{beta-1}),
\label{thap4}
\end{equation}
and define $C(u)$ and estimate it by H\"{o}lder's inequality as follows: 
\begin{equation}
C(u)=\int_0^4\int_{B_2}|u|^2\, d\|V_t\|dt\leq \Cl[c]{cuc}(p,q,E_1)\|u\|^2.
\label{bb2}
\end{equation}
Define for $t\in [\frac{1}{2}, \frac{7}{2}]$
\begin{equation}
E(t)=\|V_t\|(\phi_{\rm rad}^2)-{\bf c}-\int_{\frac{1}{2}}^t\int_{B_2} |u|^2\phi_{\rm rad}^2\, d\|V_s\|ds-\Cl[c]{c-3} \beta_*^2 (t-\frac{1}{2}),
\label{th1}
\end{equation}
where $\Cr{c-3}$ will be fixed later.
We first prove that 
\begin{equation}
E(t_2)-E(t_1)\leq -\frac14 \int_{t_1}^{t_2}\int_{B_2} |h(V_t,\cdot)|^2\phi_{\rm rad}^2\, d\|V_t\|dt,
\hspace{.2cm}\frac{1}{2}\leq \forall t_1<\forall t_2\leq \frac{7}{2}.
\label{th2}
\end{equation}
By (A1), (A2) and (A4), for a.e$.$ $t$, we have $V_t\in {\bf IV}_1(B_2)$, $h(V_t,\cdot)\in L^2(\|V_t\|)$, $u(\cdot,t)\in L^2(\|V_t\|)$. At such time $t$, using the perpendicularity of 
mean curvature \eqref{perpthm}, (omitting $t$ dependence for simplicity)
\begin{equation}
{\mathcal B}(V,u,\phi_{\rm rad}^2)\leq \int_{B_2}-|h|^2\phi_{\rm rad}^2+\phi_{\rm rad}^2|h||u|
+|u^{\perp}\cdot\nabla\phi_{\rm rad}^2|+(\nabla\phi_{\rm rad}^2)^{\perp}\cdot h\, d\|V\|.
\label{th3}
\end{equation}
The last term of \eqref{th3} may be computed as 
\begin{equation}
\int_{G_1(B_2)} 2\phi_{\rm rad}S^{\perp}(\nabla\phi_{\rm rad})\cdot h\, dV
\leq \frac14 \int_{B_2}|h|^2\phi_{\rm rad}^2\,d\|V\|+4\int_{G_1(B_2)}|S^{\perp}(\nabla\phi_{\rm rad})|^2\, dV.
\label{th4}
\end{equation}
Note that $\nabla\phi_{\rm rad}$ is $0$ outside of $B_{\frac32}\setminus B_1$, and ${\rm spt}\, \|V\|$
is represented as the union of three graphs of $C^1$ functions by \eqref{smep5} and \eqref{smep7} in $B_{\frac32}
\setminus B_1$. 
Consider the neighborhood of $(1,0)$ in which ${\rm spt}\, \|V\|$ is represented by $f_1$. 
Since $\phi_{\rm rad}$ is radially symmetric function, $\nabla\phi_{\rm rad}$ at
$(s,f_1(s))$ is parallel to $(s,f_1(s))$ and $|\nabla\phi_{\rm rad}|\leq 4$. 
On the other hand, the projection matrix $S^{\perp}$ at the same point is easily seen to be
\begin{equation}
S^{\perp}=(1+(f_1')^2)^{-1}\left(\begin{array}{ll} (f_1')^2 & -f_1' \\ -f_1' & 1\end{array}\right)
\label{th5}
\end{equation}
which is obtained by computing $I-\hat\nu\otimes\hat\nu$ with $\hat\nu=(1+(f_1')^2)^{-\frac12}(1,f_1')$, $I$ being the
identity $2\times 2$ matrix. Thus, we have
\begin{equation}
|S^{\perp}(\nabla\phi_{\rm rad})|\leq \frac{4}{
\sqrt{s^2+f_1(s)^2}}\big|S^{\perp}\left(\begin{array}{l} s\\ f_1(s)\end{array}
\right)\big|\leq c\sqrt{(f_1')^2+f_1^2}\leq c\beta_*
\label{th6}
\end{equation}
by \eqref{thap4}, where $c$ is an absolute constant. 
We have similar computations for $f_2$ and $f_3$. 
Thus by \eqref{th4} and \eqref{th6}, the last term of \eqref{th3} may be
estimate by
\begin{equation}
\int_{B_2}(\nabla\phi_{\rm rad}^2)^{\perp}\cdot h\, d\|V\|
\leq \frac14\int_{B_2}|h|\phi_{\rm rad}^2\, d\|V\|+c\beta_*^2.
\label{th7}
\end{equation}
The same computations show that the third term of \eqref{th3} may be estimated by
\begin{equation}
\int_{B_2}|u^{\perp}\cdot\nabla\phi_{\rm rad}^2|\, d\|V\|
\leq \frac12 \int_{B_2}|u|^2\phi_{\rm rad}^2\, d\|V\|+c\beta_*^2.
\label{th8}
\end{equation}
The second term of \eqref{th3} may be estimated by 
\begin{equation}
\int_{B_2}\phi_{\rm rad}^2|h||u|\, d\|V\|\leq \frac12\int_{B_2}\phi_{\rm rad}^2|h|^2\,d\|V\|
+\frac12\int_{B_2}\phi_{\rm rad}^2|u|^2\, d\|V\|.
\label{th9}
\end{equation}
Combining \eqref{th3}, \eqref{th7}-\eqref{th9}, we obtain (by recovering the notation for $t$ dependence)
\begin{equation}
{\mathcal B}(V_t,u(\cdot,t),\phi_{\rm rad}^2)\leq -\frac14\int_{B_2}|h(V_t,\cdot)|^2\phi_{\rm rad}^2
\, d\|V_t\|+\int_{B_2}|u(\cdot,t)|^2\phi_{\rm rad}^2\, d\|V_t\|+\Cr{c-3}\beta_*^2,
\label{th10}
\end{equation}
where $\Cr{c-3}$ is an absolute constant, and this holds for a.e$.$ $t
\in [\frac{1}{2},\frac{7}{2}]$. Due to (A4) and \eqref{th10}, now it is 
clear that the inequality \eqref{th2} holds if we define $E(t)$ as in
\eqref{th1}. 
We restrict $\Cr{e-1}$ further by
\begin{equation}
\Cr{e-1}^2\leq \min\left\{\frac{\Cr{ha}}{4\Cr{cuc}}, \frac{\Cr{ha}}{16\Cr{c-3}\Cr{c-p-1}^{2}}\right\}.
\label{th12.1}
\end{equation}
We proceed to prove \eqref{thap9}. 
For a.e$.$ $t\in [\frac12,\frac72]$, we have \eqref{firm1}-\eqref{fir6}. 
Denoting
\begin{equation}
\hat{E}(t)=\|V_t\|(\phi_{\rm rad}^2)-{\bf c},
\label{th13}
\end{equation}
\begin{equation}
\alpha(t)=\left(\int_{B_2}|h(V_t,\cdot)|^2\phi_{\rm rad}^2\, d\|V_t\|\right)^{\frac12},\,\,
\mu(t)=\left(\int_{B_2}{\rm dist}\, (\cdot, J)^2\, d\|V_t\|\right)^{\frac12},
\label{th14}
\end{equation}
Corollary \ref{fircor} and the definition of $\beta_*$ in \eqref{thap4} show that
\begin{equation}
3\Cr{c-1}\beta_*^2\leq |\hat{E}(t)|\, \Longrightarrow\, \alpha(t)^2 
\geq \min\{\Cr{alpha-1}^2,
(2\Cr{c-1}\mu(t))^{-2}|\hat{E}(t)|^2\}.
\label{th15}
\end{equation}
Fix any ${\hat s}\in [1,3]$ and we first prove the following upper bound,
\begin{equation}
{\hat E}({\hat s})\leq \Cr{ha2}\mu^2+ (3\Cr{c-1}+4\Cr{c-3})\beta_*^2+2C(u).
\label{th16}
\end{equation}
{\it Proof of \eqref{th16}}.
\newline
{\bf (i)} Suppose that there exists some $t_0\in [\frac12, 1]$
such that 
\begin{equation}
{\hat E}(t_0)<(3\Cr{c-1}+\Cr{c-3})\beta_*^2+C(u).
\label{th17}
\end{equation}
By the monotone decreasing property of $E(\cdot)$, \eqref{th1} and \eqref{th17}, we then have
\begin{equation}
E({\hat s})\leq E(t_0)\leq {\hat E}(t_0)<(3\Cr{c-1}+\Cr{c-3})\beta_*^2+C(u).
\label{th18}
\end{equation}
But then again by \eqref{th1}, \eqref{th18} and \eqref{thap4}, we have
\begin{equation}
{\hat E}({\hat s})\leq E({\hat s})+C(u)+3\Cr{c-3}\beta_*^2
<(3\Cr{c-1}+4\Cr{c-3})\beta_*^2+2C(u).
\label{th19}
\end{equation}
With \eqref{th19} we proved \eqref{th16} under the assumption of
(i). Now consider the complementary situation.
\newline
{\bf (ii)} Suppose that for all $t\in [\frac12,1]$, we have
\begin{equation}
{\hat E}(t)\geq (3\Cr{c-1}+\Cr{c-3})\beta_*^2+C(u).
\label{th20}
\end{equation}
This in particular means $|{\hat E}(t)|\geq 3\Cr{c-1}\beta_*^2$, 
thus \eqref{th20} and \eqref{th15} show
\begin{equation}
\alpha(t)^2\geq \min\{\Cr{alpha-1}^2,(2\Cr{c-1}\mu(t))^{-2}|{\hat E}(t)|^2\}
\label{th21}
\end{equation}
for a.e$.$ $t\in [\frac12,1]$. 
By \eqref{th1} and \eqref{th20},
\begin{equation}
{\hat E}(t)\geq E(t)\geq {\hat E}(t)-\Cr{c-3}\beta_*^2-C(u)
\geq 0.
\label{th22}
\end{equation}
\eqref{th22} shows $|{\hat E}(t)|\geq |E(t)|$ in particular and by \eqref{th21} 
and \eqref{th11.1} we have
\begin{equation}
\alpha(t)^2  
\geq 4 P\min\{1,\mu(t)^{-2}|E(t)|^2\}
\label{th23}
\end{equation}
for a.e$.$ $t\in [\frac12,1]$.
Now we are in the position to apply Lemma \ref{ode}. 
Define
\begin{equation}
\Phi(t)=E(t+\frac12),\,\, f(t)=P\min\{1,\mu(t+\frac12)^{-2}
|\Phi(t)|^2\},\,\,\,t\in [0,\frac12]. 
\label{th24}
\end{equation}
By \eqref{th2} \eqref{th23} and \eqref{th24}, we have for $0\leq \forall t_1<\forall t_2\leq \frac12$
\begin{equation}
\Phi(t_2)-\Phi(t_1)\leq -\frac14 \int_{t_1+\frac12}^{t_2+\frac12}
\alpha(t)^2\, dt\leq -\int_{t_1}^{t_2}f(t)\, dt.
\label{th25}
\end{equation}
We also have 
\begin{equation}
\Phi(0)=E(\frac12)= \|V_{\frac12}\|(\phi_{\rm rad}^2)-{\bf c}\leq \frac{\Cr{ha}}{2}
\label{th26}
\end{equation}
by \eqref{th1} and \eqref{va1}. Hence the assumptions of Lemma \ref{ode}
are all satisfied with $g(t)=\mu(t+\frac12)$, and noticing that $\|g\|_{L^2}^2\leq \mu^2$, 
we conclude
\begin{equation}
 (E(1)=)\,\Phi(\frac12)\leq \Cr{ha2}\mu^2.
\label{th27}
\end{equation}
For any $\hat s \in [1,3]$, by \eqref{th1}, \eqref{th27} and the monotone decreasing property of $E$, we have
\begin{equation}
{\hat E}({\hat s})\leq E({\hat s})+C(u)+3\Cr{c-3}\beta_*^2 
\leq \Cr{ha2}\mu^2+C(u)+3\Cr{c-3}\beta_*^2.
\label{th28}
\end{equation}
Thus \eqref{th28} shows that 
\eqref{th16} holds under the assumption of (ii). This concludes the proof of \eqref{th16}. 
\newline
Fix any ${\hat s}\in [1,3]$ and we next prove the following lower
bound,
\begin{equation}
-\Cr{ha2}\mu^2-(3\Cr{c-1}+8\Cr{c-3})\beta_*^2-2C(u)\leq {\hat E}({\hat s}).
\label{th29}
\end{equation}
The idea is similar to the upper bound estimate with a few differences, 
but we present the proof for the completeness. 
\newline
{\it Proof of \eqref{th29}}.
\newline
{\bf (i)} Suppose that there exists some $t_0\in [3,\frac72]$
such that
\begin{equation}
{\hat E}(t_0)>-(3\Cr{c-1}+4\Cr{c-3})\beta_*^2 -C(u).
\label{th30}
\end{equation}
By \eqref{th1} and \eqref{th30},
\begin{equation}
E(t_0)\geq {\hat E}(t_0)-C(u)-4\Cr{c-3}\beta_*^2
>-(3\Cr{c-1}+8\Cr{c-3})\beta_*^2-2C(u).
\label{th31}
\end{equation}
By the monotone decreasing property of $E$, we have
$E({\hat s})\geq E(t_0)$ while ${\hat E}({\hat s})\geq E({\hat s})$ by
\eqref{th1}. Thus \eqref{th31} proves \eqref{th29} in case of (i).
\newline
{\bf (ii)} Suppose that for all $t\in [3,\frac72]$, we have
\begin{equation}
{\hat E}(t)\leq -(3\Cr{c-1}+4\Cr{c-3})\beta_*^2-C(u).
\label{th32}
\end{equation}
This means $|{\hat E}(t)|\geq 3\Cr{c-1}\beta_*^2$, thus by \eqref{th15}, we have
\eqref{th21} for a.e$.$ $t\in [3,\frac72]$. 
We need to change ${\hat E}$ in \eqref{th21} to $E$. To do so, observe that
\begin{equation}
|E(t)|\leq |{\hat E}(t)|+C(u)+4\Cr{c-3}\beta_*^2
\leq 2|{\hat E}(t)|,
\label{th33}
\end{equation}
the last inequality of \eqref{th33} coming from \eqref{th32}. 
Thus, \eqref{th21} with \eqref{th33} (as well as recalling \eqref{th11.1}) shows
\eqref{th23} for a.e$.$ $t\in[3,\frac72]$. Again we apply
Lemma \ref{ode}. Set
\begin{equation}
\Phi(t)=E(t+3),\,\, f(t)=P\min\{1,\mu(t+3)^{-2}|\Phi(t)|^2\},\,\, t\in [0,\frac12].
\label{th34}
\end{equation}
By having \eqref{th23}, we have \eqref{th25} and 
\begin{equation}
\Phi(\frac12)=E(\frac72)\geq \|V_{\frac72}\|(\phi_{\rm rad}^2)-{\bf c}-C(u)-4\Cr{c-3}\beta_*^2
\geq -\frac{\Cr{ha}}{2}-\Cr{cuc}\Cr{e-1}^2-4\Cr{c-3}\Cr{c-p-1}^2\Cr{e-1}^2\geq -\Cr{ha}
\label{th35}
\end{equation}
by \eqref{th34}, \eqref{th1}, \eqref{va2}, \eqref{bb2}, \eqref{thap4} and the last inequality due to \eqref{th12.1}. 
The assumptions of Lemma \ref{ode} (for case (2))
are thus satisfied, and we obtain
\begin{equation}
-\Cr{ha2}\mu^2\leq \Phi(0)(=E(3)).
\label{th36}
\end{equation}
Since $E$ is decreasing, for any $\hat s\in [1,3]$, we have
\begin{equation}
\hat E (\hat s)\geq E(\hat s) \geq E(3).
\label{th36s}
\end{equation}
Hence under the assumption of (ii), \eqref{th36} and \eqref{th36s} show \eqref{th29}.

Since ${\hat s}\in [1,3]$ is arbitrary, \eqref{th16} and
\eqref{th29} combined with \eqref{thap4} and \eqref{bb2} prove the first claim \eqref{thap9}
with a suitable constant $\Cr{c-2}$. 

To prove \eqref{thap10}, observe that \eqref{th2} with $t_2=3$ and
$t_1=1$ shows (recalling \eqref{th1})
\begin{equation}
\int_{1}^{3}\int_{B_2}|h|^2\phi_{\rm rad}^2\, d\|V_t\|dt
\leq 4(E(1)-E(3)) 
\leq 4(|{\hat E}(1)|+
|{\hat E}(3)|+C(u)+2\Cr{c-3}\beta_*^2).
\label{th37}
\end{equation}
Then using \eqref{thap9} to \eqref{th37}, we obtain \eqref{thap10}, again with a suitable choice of $\Cr{c-2}$.
\hfill{$\Box$}
\vspace{.2in}

In the next two sections we derive further a priori estimates for the flow $\{V_{t}\}$ whenever it is  weakly close, in space-time at scale one, to the static triple junction $J$. These estimates provide enough control of the behavior of the moving curves near the singularity of $J$ for us to establish (in Section~\ref{blowupsection}) decay, by a fixed factor at a fixed smaller scale,  of the space-time $L^{2}$ distance of the flow to $J$, and consequently (by iterating this decay result) Theorem~\ref{mainreg}. These estimates are in the spirit of those proved first by L.~Simon (\cite{Simon2}), for a similar purpose,  for the case of multiplicity 1 minimal submanifolds weakly close to certain cylindrical minimal cones (in arbitrary dimension and codimension). However, in the present parabolic setting, their statements are often different and proofs require new ideas.

\section{A priori estimates II: non-concentration of the $L^2$-distance near the singularity of $J$}
The main result in this section is the estimate (\ref{pa21}) of Proposition~\ref{tildest}. This estimate in  full strength plays an important role in Section~\ref{blowupsection} where we establish asymptotics for the blow-ups of sequences of flows converging weakly to the static triple junction $J$. It also implies that the space-time $L^{2}$ distance $\mu$ of the flow from $J$ does not concentrate near the singularity of $J$, a fact that is indispensable in the proof of the key decay result (Proposition~\ref{blowprop7}) for $\mu$. 

An essential ingredient in the proof of Proposition~\ref{tildest} is Proposition~\ref{pade} below, which is based on  the results of Section~\ref{curv-est} and (the main idea behind) Huisken's monotonicity formula.  
\begin{prop} Corresponding to $\nu \in (0, 1)$, $E_1 \in [1, \infty)$ and $p$, $q$ as in $({\rm A}0),$
 there exist $\Cl[eps]{e-2}\in (0,\Cr{e-1}]$ (where $\Cr{e-1} = \Cr{e-1}(\nu, E_{1}, p, q)$ is as in Proposition~\ref{hde}) and $\Cl[c]{c-4}\in (1,\infty)$ with the following property: 
If $\{V_t\}_{t\in [0,4]}$ and $\{u(\cdot, t)\}_{t\in [0,4]}$
satisfy (A1)-(A4) with $U = B_{2}$ and $I =[0,4]$,  if $\eqref{smep1}, \eqref{smep2}, \eqref{smep3}, \eqref{smep4}$ hold with $\Cr{e-2}$ in place of $\Cr{e-1}$ and if
\begin{equation}
V_{t_0}\in {\bf IV}_1(B_2),\ \ h(V_{t_0},\cdot)\in L^2(\|V_{t_0}\|)\ \mbox{ and }\ 
\Theta(\|V_{t_0}\|,0)=\frac32
\label{pa1}
\end{equation}
for some $t_0\in [\frac32, 3]$, then 
\begin{equation}
\int_{5/4}^{t_0}\int_{B_1} \left|h+\frac{x^{\perp}}{2(t_0-t)}\right|^2\rho_{(0,t_0)}(x,t)
\, d\|V_t\|dt\leq \Cr{c-4}\max\{\mu,\|u\|\}^2.
\label{pa2}
\end{equation}
 \label{pade}
\end{prop}
{\it Proof}. By \eqref{pa1}, there exists a tangent cone to $V_{t_0}$ at $x=0$ which is  
$|{\bf R}_{\theta}(J)|$ for some $\theta\in [0,2\pi)$. From this fact, it follows that 
\begin{equation}
\lim_{\epsilon\searrow 0} \int_{B_2}\phi_{\rm rad}^2(x)\rho_{(0,t_0+\epsilon)}(x,t_0)\, d\|V_{t_0}\|(x)=\frac32.
\label{pa3}
\end{equation}
In the following, we fix $\epsilon>0$ arbitrarily close to 0. We choose
$t_1\in [1,\frac54]$ so that
\begin{equation}
\int_{B_2}|h(V_{t_1},\cdot)|^2\phi_{\rm rad}^2\, d\|V_{t_1}\|\leq
8\Cr{c-2}\max\{\mu,\|u\|\}^2
\label{pa4}
\end{equation}
where $\Cr{c-2} = \Cr{c-2}(\nu, E_{1}, p, q)$ is as in Proposition~\ref{hde}. This is possible in view of the estimate \eqref{thap10}.
Arguing as in Proposition~\ref{firprop}, for $\Cr{e-2}$ suitably small (so that
$8\Cr{c-2}\max\{\mu,\|u\|\}^2\leq \Cr{alpha-1}^2$ where $\alpha_{1} = \alpha_{1}(E_{1})$ is as in Proposition~\ref{firprop}), we may 
conclude that ${\rm spt}\,\|V_{t_1}\| \cap B_{1}$ consists of three $C^{1,\frac12}$ curves $l_1,l_2,l_3$ meeting
at a common point $p$ near the origin, with associated numbers $s_{1}, s_{2}, s_{3} \in (-1/2, 1/2)$ and functions $f_{j} \in C^{1, 1/2}([s_{j}, 1])$, $j=1, 2, 3,$ such that 
\begin{equation}
{\bf R}_{-\frac{2\pi (j-1)}{3}} l_{j} = \{(s,f_j(s))\in {\mathbb R}^2\,:\, s\in [s_j,1]\}
\end{equation}
for $j=1, 2, 3;$ furthermore, using the estimate \eqref{fir10.5}, we see that $\sup \, \{{\rm dist} \, (x, J) \, : \, x \in l_{j}\}$ for $j=1,2,3,$ and hence also $|s_j|$, are all $\leq 2\beta_*+2\sqrt{8\Cr{c-2}}\max\{\mu,\|u\|\}$, where $\beta_*$ is as in \eqref{thap4}. 
These estimates and radial symmetry of $\phi_{\rm rad}$ and $\rho_{(0,t_0)}(\cdot,t_1)$ imply, 
for a suitable choice of $\Cr{c-4}$ depending only on $p,q,\nu,E_1,$ that
\begin{equation}
\left|\int_J \phi_{\rm rad}^2\rho_{(0,t_0)}(\cdot,t_1)\, d{\mathcal H}^1
-\int_{B_2} \phi_{\rm rad}^2\rho_{(0,t_0)}(\cdot,t_1)\, d\|V_{t_1}\|\right|\leq\Cr{c-4}
\max\{\mu,\|u\|\}^2.
\label{pa5}
\end{equation}
Here it is important that $t_1\in [1,\frac54]$ so that $t_0-t_1\geq \frac14$, allowing the choice
of $\Cr{c-4}$ to be independent of $t_0$ and $t_1$. 

We next use $\rho_{(0,t_0+\epsilon)}(\cdot,t)\phi_{\rm rad}^2$ as a 
test function in \eqref{meq} with $t\in [t_1, t_0]$. For simplicity of notation write $\rho$ for 
$\rho_{(0, t_{0} + \epsilon)}$ and define ${\hat \rho}(x, t) =\rho_{(0,
t_0+\epsilon)}(x,t)\phi_{\rm rad}(x)^2$.  By direct computation, 
\begin{equation}
\begin{split}
{\mathcal B}(V_t,u(\cdot,t),{\hat \rho}(\cdot,t))&=\int_{B_2}(-h{\hat \rho}+\nabla{\hat\rho})\cdot(h+u^{\perp})\, d\|V_t\|
\\ & =\int_{B_2}-|h|^2{\hat\rho}+2\nabla{\hat\rho}\cdot h+u^{\perp}\cdot(-h{\hat\rho}+\nabla{\hat\rho})
-\nabla{\hat\rho}\cdot h\, d\|V_t\|\\
&=\int_{G_1(B_2)}-{\hat\rho}\big|h-\frac{(\nabla{\hat\rho})^{\perp}}{\hat\rho}\big|^2
+\frac{|(\nabla{\hat\rho})^{\perp}|^2}{\hat\rho}+u\cdot(-h{\hat\rho}+(\nabla{\hat\rho})^{\perp})+S\cdot\nabla^2
{\hat\rho}\, d V_t \\
&\leq \int_{G_1(B_2)}-\frac{\hat\rho}{2}\big|h-\frac{(\nabla{\hat\rho})^{\perp}}{\hat\rho}\big|^2
+\frac{|(\nabla{\hat\rho})^{\perp}|^2}{\hat\rho}+\frac{|u|^2{\hat\rho}}{2}+S\cdot\nabla^2{\hat\rho}\,
dV_t,
\end{split}
\label{pa7}
\end{equation}
where we have used the fact that by \eqref{perpthm}, for $\|V_t\|$ a.e$.$, $h\cdot \nabla{\hat\rho}=
h\cdot(\nabla{\hat\rho})^{\perp}$. 
We now need to carefully evaluate the terms involving $\nabla\phi_{\rm rad}^2$ in the above. For the second term on the right hand side of \eqref{pa7}, we have
\begin{equation}
\frac{|(\nabla{\hat\rho})^{\perp}|^2}{\hat\rho}=\phi_{\rm rad}^2\frac{|(\nabla\rho)^{\perp}|^2}{\rho}
+2(\nabla\rho)^{\perp}\cdot(\nabla\phi_{\rm rad}^2)^{\perp}
+\rho\frac{|(\nabla\phi_{\rm rad}^2)^{\perp}|^2}{\phi_{\rm rad}^2}
\leq \phi_{\rm rad}^2\frac{|(\nabla\rho)^{\perp}|^2}{\rho}+c\beta_*^2
\label{pa8}
\end{equation}
where the last inequality follows from the estimate \eqref{th6} which holds also with $\rho$ in place of $\phi_{\rm rad}$.
The constant $c$ is an absolute constant which may differ from line to line.
By combining \eqref{pa7} and \eqref{pa8}, we obtain
\begin{equation}
\begin{split}
{\mathcal B}(V_t,u(\cdot,t),{\hat\rho}(\cdot,t))\leq& \int_{G_1(B_2)}
-\frac{\hat\rho}{2}\big|h-\frac{(\nabla{\hat\rho})^{\perp}}{\hat\rho}\big|^2
+\phi_{\rm rad}^2\big(\frac{|(\nabla\rho)^{\perp}|^2}{\rho}+S\cdot\nabla^2\rho\big)\\
&\hspace{.5in}+\frac{|u|^2{\hat\rho}}{2}+c\beta_*^2 +2S\cdot(\nabla\rho\otimes
\nabla\phi_{\rm rad}^2)+\rho S\cdot\nabla^2\phi_{\rm rad}^2\, dV_t.
\end{split}
\label{pa9}
\end{equation}
By \eqref{pa9}, (A4) and the identity 
\begin{equation}\label{idenbh}
\frac{\partial \rho}{\partial t}+S\cdot \nabla^2\rho+\frac{|S^{\perp}
(\nabla\rho)|^2}{\rho}=0,
\end{equation}
we conclude that
\begin{equation}
\begin{split}
\left.\int_{B_2}{\hat\rho}(\cdot,t)\, d\|V_t\|\right|_{t=t_1}^{t_0}\leq& \int_{t_1}^{t_0}\int_{G_1(B_2)}-\frac{\hat\rho}{2}
\big|h-\frac{(\nabla{\hat\rho})^{\perp}}{\hat\rho}\big|^2 
+\frac{|u|^2{\hat\rho}}{2}+c\beta_*^2 \\ & \hspace{1in}+2S\cdot(\nabla\rho\otimes
\nabla\phi_{\rm rad}^2)+\rho S\cdot\nabla^2\phi_{\rm rad}^2\, dV_t dt.
\end{split}
\label{pa10}
\end{equation}
We next proceed to estimate the last two terms on the right hand side of \eqref{pa10}. 
Note that the integrands in these two terms are zero outside $B_{\frac32}\setminus
B_{\frac12}$. (So in particular the derivatives of $\rho$ appearing there are bounded uniformly.) 
Since by Proposition~\ref{smprop} integration with respect to $\|V_{t}\|$ in $B_{\frac 32} \setminus B_{\frac 12}$ is along three $C^{1,\frac12}$ curves $l_{1}^{(t)}, l_{2}^{(t)}, l_{3}^{(t)}$ represented as graphs of  functions $f_1(\cdot, t),f_2(\cdot, t),f_3(\cdot, t)$ respectively as in \eqref{smep7}, one can compute them explicitly. For instance for  the curve $l_{1}^{(t)}$, 
by explicit calculation (suppressing the $t$ dependence of the functions involved), 
\begin{equation}
\begin{split}
\int_{l_{1}^{(t)}} & 2S\cdot(\nabla\rho\otimes
\nabla\phi_{\rm rad}^2)+\rho S\cdot\nabla^2\phi_{\rm rad}^2\\
 & \begin{split}= \int_{\frac12}^{\frac32}&
(1+(f_1')^2)^{-\frac12}\left(\begin{array}{ll} 1 & f_1' \\ f_1' & (f_1')^2\end{array}\right)\cdot
\left\{2r^{-2}(x\otimes x)\frac{d\rho}{dr}\frac{d\phi_{\rm rad}^2}{dr}\right. \\ &+\left.\rho r^{-2}
(x\otimes x)(\frac{d^2}{dr^2}-r^{-1}\frac{d}{dr} )\phi_{\rm rad}^2+r^{-1}\rho I\frac{d\phi_{\rm rad}^2}{dr}\right\}ds,
\end{split}
\end{split}
\label{pa11}
\end{equation}
where $f_1=f_1(s,t)$, $r=\sqrt{s^2+(f_1)^2}$, $x=(s,f_1(s,t))$, $I$ is the identity $2\times 2$
matrix and $d/dr$ is the differentiation with respect to the radial direction. Since 
$|f_1|,\,|f_1'|\leq \beta_*$ by \eqref{thap4}, we may estimate terms
on the right hand side of \eqref{pa11} up to errors of order $\beta_*^2$ to obtain 
\begin{equation}
\begin{split}
\int_{l_{1}^{(t)}} & 2S\cdot(\nabla\rho\otimes
\nabla\phi_{\rm rad}^2)+\rho S\cdot\nabla^2\phi_{\rm rad}^2\\ &\leq \int_{\frac12}^{\frac32}\left(2\frac{d\rho}{dr}\frac{d\phi_{\rm rad}^2}{dr}+
\rho\left(\frac{d^2}{dr^2}-s^{-1}\frac{d}{dr}\right)\phi_{\rm rad}^2+s^{-1}\rho\frac{d\phi_{\rm rad}^2}{dr}\right)\, ds
+c\beta_*^2\\
&=\int_{\frac12}^{\frac32}\left(2\frac{d\rho}{dr}\frac{d\phi_{\rm rad}^2}{dr}+
\rho\frac{d^2\phi_{\rm rad}^2}{dr^2}\right)\, ds+c\beta_*^2.
\end{split}
\label{pa12}
\end{equation}
Since the functions appearing in the integrand on the right hand side of the above are radially symmetric, their values at 
$(s,f_1(s,t))$ and those at $(s,0)$ differ by at most $c\beta_*^2$. Thus we have 
\begin{equation}
\begin{split}
\int_{l_{1}^{(t)}} & 2S\cdot(\nabla\rho\otimes
\nabla\phi_{\rm rad}^2)+\rho S\cdot\nabla^2\phi_{\rm rad}^2\\ & \leq \int_{\frac12}^{\frac32}\left(2\frac{\partial \rho}{\partial x}\frac{\partial \phi_{\rm rad}^2}{\partial x}
+\rho\frac{\partial^2 \phi_{\rm rad}^2}{\partial x^2}\right)(s,0)\, ds+c\beta_*^2
=\int_{\frac12}^{\frac32} \frac{\partial \rho}{\partial x}\frac{\partial \phi_{\rm rad}^2}{\partial x}(s,0)\, ds
+c\beta_*^2
\end{split}
\label{pa13}
\end{equation}
where we integrated by parts and used the property that $\frac{\partial\phi_{\rm rad}^2}{\partial x}(s,0)=0$
at $s=1\pm \frac12$. The same computation holds after rotation for the other two curves $l_{2}^{(t)}, l_{3}^{(t)}$, so by \eqref{pa10} we  deduce 
\begin{equation}
\begin{split}
\left.\int_{B_2}{\hat\rho}(\cdot,t)\, d\|V_t\|\right|_{t=t_1}^{t_0}\leq  &
\int_{t_1}^{t_0}\int_{B_2}-\frac{\hat\rho}{2}
\big|h-\frac{(\nabla{\hat\rho})^{\perp}}{\hat\rho}\big|^2 
+\frac{|u|^2{\hat\rho}}{2}\, d\|V_t\|dt \\
& +3 \int_{t_1}^{t_0}\int_{\{(x,0)\in{\mathbb R}^2\,:\, |x-1|<\frac12\}} \frac{\partial \rho}{\partial x}\frac{\partial \phi_{\rm rad}^2}{\partial x}\, d{\mathcal H}^1 dt+c\beta_*^2.
\end{split}
\label{pa14}
\end{equation}
We note that 
\begin{equation*}
\left.\int_{J}{\hat\rho}(\cdot,t)\, d{\mathcal H}^1\right|_{t=t_1}^{t_0}=\int_J\int_{t_1}^{t_0}
\frac{\partial{\hat\rho}}{\partial t}(\cdot,t)\, dtd{\mathcal H}^1=3\int_{t_1}^{t_0}\int_{\{(x,0)\in
{\mathbb R}^2\,:\, x\geq 0\}}\frac{\partial{\hat\rho}}{\partial t}(\cdot,t)d{\mathcal H}^1 dt
\label{pa15}
\end{equation*}
by radial symmetry of $\hat\rho$, and thus, since $\frac{\partial{\hat\rho}}{\partial
t}=-\phi_{\rm rad}^2 \frac{\partial^2 \rho}{\partial x^2}$ on the $x$-axis, 
\begin{equation}
\begin{split}
\left.\int_{J}{\hat\rho}(\cdot,t)\, d{\mathcal H}^1\right|_{t=t_1}^{t_0}&=3\int_{t_1}^{t_0}\int_{\{(x,0)\in
{\mathbb R}^2\,:\, x\geq 0\}} \frac{\partial \rho}{\partial x}\frac{\partial \phi_{\rm rad}^2}{\partial x}d{\mathcal H}^1 dt\\
&=3\int_{t_1}^{t_0}\int_{\{(x,0)\in
{\mathbb R}^2\,:\, |x-1|<\frac12\}} \frac{\partial \rho}{\partial x}\frac{\partial \phi_{\rm rad}^2}{\partial x}d{\mathcal H}^1 dt
\end{split}
\label{pa16}
\end{equation}
where we used integration by parts and the fact that $\frac{\partial\rho}{\partial x}=0$ at $x=0$ and $\phi_{\rm rad}^2=0$ at $x=\frac32$.
Substituting \eqref{pa16} into \eqref{pa14}, we obtain
\begin{equation}
\left.\big(\int_{B_2}{\hat\rho}(\cdot,t)\, d\|V_t\|-\int_J{\hat\rho}(\cdot,t)\, d{\mathcal H}^1\big)\right|_{t=t_1}^{t_0}
\leq \int_{t_1}^{t_0}\int_{B_2}-\frac{\hat\rho}{2}\big|h-\frac{(\nabla{\hat\rho})^{\perp}}{\hat\rho}\big|^2+{\hat\rho}
\frac{|u|^2}{2}\, d\|V_t\|dt+c\beta_*^2.
\label{pa17}
\end{equation}
We now let $\epsilon\rightarrow 0$ in \eqref{pa17}. Since $\int_J{\hat\rho}_{(0,t_0+\epsilon)}(x,t_0)\, d{\mathcal H}^1
\rightarrow \frac32$, in view of \eqref{pa3} and \eqref{pa5}, we obtain from \eqref{pa17}  (using also 
the fact that $\phi_{\rm rad}=1$ on $B_1$ and $=0$ on ${\mathbb R}^{2} \setminus B_{\frac32}$) that
\begin{equation}
\int_{t_1}^{t_0}\int_{B_1}\frac{\rho}{2}\big|h-\frac{(\nabla\rho)^{\perp}}{\rho}\big|^2\, d\|V_t\|dt\leq \int_{t_1}^{t_0}
\int_{B_{\frac32}}{\rho}\frac{|u|^2}{2}\, d\|V_t\|dt+c\beta_*^2+\Cr{c-4}\max\{\mu,\|u\|\}^2
\label{pa18}
\end{equation}
where $\rho=\rho_{(0,t_0)}(x,t)$. Lastly, the term above involving $u$ may be estimated as in \cite[(6.7)-(6.8)]{Kasai-Tonegawa} to get 
\begin{equation}
\int_{t_1}^{t_0}\int_{B_{\frac32}}\rho|u|^2\, d\|V_t\|dt\leq c(p,q)E_1 ^{1-\frac{2}{p}}\|u\|^2.
\label{pa19}
\end{equation}
Since $\frac{\nabla\rho}{\rho}=-\frac{x}{2(t_0-t)}$, the desired estimate follows after redefining $\Cr{c-4}$ depending only on $p,q,\nu,E_1$. 
\hfill{$\Box$}

\begin{prop} Fix $\kappa\in [0,1)$. Under the same assumptions as in  Proposition \ref{pade}, we have
\begin{equation}
\sup_{t\in [\frac54,t_0)}(t_0-t)^{-\kappa}\int_{B_{\frac34}} \rho_{(0,t_0)}(\cdot,t){\rm dist}\, (\cdot,J)^2\, 
d\|V_t\|\leq \Cl[c]{c-5}\max\{\mu,\|u\|\}^2
\label{pa21}
\end{equation}
where $\Cr{c-5}$ depends only on $\kappa,p,q,\nu,E_1$. 
\label{tildest}
\end{prop}
{\it Proof}. 
Define ${\tilde d}:{\mathbb R}^2\rightarrow{\mathbb R}$ such that ${\tilde d}$ 
is positively homogeneous of degree one (i.e.\ ${\tilde d}(\lambda x)=\lambda {\tilde d}(x)$ $\forall \lambda\geq 0$
and $\forall x\in{\mathbb R}^2$), smooth away from $J$, and
\begin{equation}
\begin{array}{ll}
 {\tilde d}(x)={\rm dist}\,(x,J)&\mbox{ $\forall x$ with }{\rm dist}\,(x,J)<\frac{|x|}{5},\\
 \frac12 {\rm dist}\, (x,J)\leq {\tilde d}(x)\leq 2\,{\rm dist}\, (x,J) & \forall x\in {\mathbb R}^2,\\
 |\nabla {\tilde d}(x)|\leq 1 & \forall x\notin J.
\end{array}
\label{pa20}
\end{equation}
By homogeneity, we have $x\cdot \nabla ({\tilde d}^2/|x|^2)=0$, which gives after a little computation that
\begin{equation}
x\cdot\nabla{\tilde d}^2=2{\tilde d}^2.
\label{pa20.5}
\end{equation} 
Let $0\leq \eta\leq 1$ be a non-negative smooth radially symmetric function such that 
\begin{equation}
\eta=0 \quad \forall x\notin B_1 \quad \mbox {and} \quad \eta=1 \quad \forall x\in B_{\frac34}.
\label{pa22}
\end{equation}
Since ${\rm spt} \, \nabla \eta \subset B_{1} \setminus B_{\frac 34}$,  we may assume that
${\rm spt}\, |\nabla\eta| \, \cap \, {\rm spt}\,\|V_t\|$ is contained in $\cup _{j=1}^{3}{\rm graph}\,f_j(\cdot,t)$ for all $t\in [1,3]$, where $f_{j}$ are as in Proposition~\ref{smprop}. Fix $t_1\in [1, \frac54]$ such that
\begin{equation}
\int_{B_2}{\rm dist}\,(\cdot,J)^2\, d\|V_{t_1}\|\leq 4 \mu^2.
\label{pa23}
\end{equation}
Such $t_1$ exists by \eqref{smep1}, the definition of $\mu$.
We next use $(t_0-t)^{-\kappa} g(t)\rho_{(0,t_0)}(x,t){\tilde d}(x)^2\eta(x)$ as a test function in 
\eqref{meq} 
over the time interval $[t_1,t_2]$ with arbitrary $t_2\in [\frac54,t_0)$, where $g$ is a fixed smooth non-negative function with
\begin{equation}
0<g(t)\leq 1
\label{pa23.5}
\end{equation}
which will be chosen later. Denoting, for notational convenience,  $\rho_{(0,t_0)}(x,t)$ by $\rho$
and $(t_0-t)^{-\kappa}g(t)\rho_{(0,t_0)}(x,t)$ by ${\hat \rho}$ respectively, we obtain from \eqref{meq}
that 
\begin{equation}
\begin{split}
&\left.\int_{B_1}\eta{\tilde d}^2{\hat\rho}\, d\|V_t\|\right|_{t=t_1}^{t_2} \leq \int_{t_1}^{t_2}\int_{B_1}\left(-h{\hat \rho}\eta{\tilde d}^2+\nabla ({\hat\rho}\eta{\tilde d}^2)\right)
\cdot(h+u^{\perp})+\eta{\tilde d}^2\frac{\partial{\hat\rho}}{\partial t}\, d\|V_t\|dt.
\end{split}
\label{pa24}
\end{equation}
Since by direct calculation and estimation 
\begin{equation*}
\begin{split}
&\left(-h{\hat \rho}\eta{\tilde d}^2+\nabla ({\hat\rho}\eta{\tilde d}^2)\right)
\cdot(h+u^{\perp})\\  
&=\left(-|h|^2 {\hat \rho}+(\nabla{\hat\rho}\cdot h)\right)\eta{\tilde d}^2+{\hat \rho}\nabla(\eta{\tilde d}^2)\cdot h
+\eta{\tilde d}^2(-h{\hat\rho}+\nabla{\hat\rho})\cdot u^{\perp}  
 +{\hat\rho}\nabla(\eta{\tilde d}^2)\cdot u^{\perp}\\
&\leq -{\hat\rho}\big|h-\frac{(\nabla{\hat\rho})^{\perp}}{\hat\rho}\big|^2\eta{\tilde d}^2-(\nabla{\hat\rho}\cdot h)
\eta{\tilde d}^2+\frac{|(\nabla{\hat\rho})^{\perp}|^2}{\hat\rho}\eta{\tilde d}^2+{\hat\rho}\nabla(\eta{\tilde d}^2)\cdot h \\
& \quad + \frac12{\hat\rho}\big| h-\frac{(\nabla{\hat\rho})^{\perp}}{\hat\rho}\big|^2\eta{\tilde d}^2+\frac12{\hat\rho}\eta{\tilde d}^2|u|^2
+{\hat\rho}\nabla(\eta{\tilde d}^2)\cdot u^{\perp},
\end{split}
\label{pa25}
\end{equation*}
it follows from \eqref{pa24} that
\begin{equation}
\begin{split}
&\left.\int_{B_1}\eta{\tilde d}^2{\hat\rho}\, d\|V_t\|\right|_{t=t_1}^{t_2} \leq
\int_{t_1}^{t_2}\int_{B_1}-(\nabla{\hat\rho}\cdot h)\eta{\tilde d}^2+\frac{|(\nabla{\hat\rho})^{\perp}|^2}{\hat\rho}\eta{\tilde d}^2
+{\hat\rho}\nabla(\eta{\tilde d}^2)\cdot h \\
& \hspace{2in}+\frac12{\hat\rho}\eta{\tilde d}^2|u|^2
+{\hat\rho}\nabla(\eta{\tilde d}^2)\cdot u^{\perp}+\eta{\tilde d}^2\frac{\partial{\hat\rho}}{\partial t}\, d\|V_t\|dt.
\end{split}
\label{pa26}
\end{equation}
Next note that by \eqref{fvf},  
\begin{equation}
\int_{B_1}-(\nabla{\hat\rho}\cdot h)\eta{\tilde d}^2\, d\|V_t\|
=\int_{G_1(B_1)}S\cdot \left(\eta{\tilde d}^2\nabla^2{\hat\rho}+\nabla{\hat\rho}
\otimes\nabla(\eta{\tilde d}^2)\right)\, dV_t(\cdot,S) \quad \mbox{for a.e.} \; t.
\label{pa27}
\end{equation}
Using \eqref{pa27} in \eqref{pa26} and using \eqref{idenbh} (keeping in mind that $\nabla {\hat\rho}=(t_0-t)^{-\kappa}g(t)\nabla\rho$, $\nabla^{2} {\hat\rho}=(t_0-t)^{-\kappa}g(t)\nabla^{2}\rho$), we obtain
\begin{equation}
\begin{split}
&\left.\int_{B_1}{\hat\rho}\eta{\tilde d}^2\, d\|V_t\|\right|_{t=t_1}^{t_2} \leq\int_{t_1}^{t_2}
\int_{G_1(B_1)}S\cdot\left(\nabla{\hat\rho}\otimes\nabla(\eta{\tilde d}^2)\right)+{\hat\rho}\nabla(\eta{\tilde d}^2)
\cdot h\\ & \hspace{1in} +\frac12{\hat\rho}\eta{\tilde d}^2|u|^2
+{\hat\rho}\nabla(\eta{\tilde d}^2)\cdot u^{\perp}
+\eta{\tilde d}^2\rho \frac{d}{d t}\left((t_0-t)^{-\kappa}g\right)\, dV_t(\cdot,S)dt \\
& =I_1+I_2+I_3+I_4+I_5
\end{split}
\label{pa29}
\end{equation}
where $I_{1}, \ldots, I_{5}$ denote the five integrals corresponding to the five summands, in the order listed,  in the integrand on the right hand side of the  above. We analyze these integrals as follows:

\noindent
{\bf Estimation of $I_1+I_2$}.
\newline
Since $S\cdot (v_1\otimes v_2)=v_1\cdot v_2-v_1^{\perp}\cdot v_2$ for any two vectors $v_1,\, v_2$, 
we have for a.e$.$ $t$ (recalling that $v^{\perp}=v-S(v)$), 
\begin{equation}
\begin{split}
S\cdot(\nabla{\hat\rho}&\otimes\nabla(\eta{\tilde d}^2))+{\hat\rho}\nabla(\eta{\tilde d}^2)\cdot h
=\nabla{\hat\rho}\cdot\nabla(\eta{\tilde d}^2)-(\nabla{\hat\rho})^{\perp}\cdot\nabla(\eta{\tilde d}^2)
+{\hat\rho}\nabla(\eta{\tilde d}^2)\cdot h \\
&=\nabla{\hat\rho}\cdot\nabla(\eta{\tilde d}^2)+{\hat\rho}\nabla(\eta{\tilde d}^2)\cdot
\left(h-\frac{(\nabla{\hat\rho})^{\perp}}{\hat\rho}\right) \\
&\leq -\frac{{\hat\rho}x}{2(t_0-t)}\cdot\nabla(\eta{\tilde d}^2)+\frac{\rho}{2}\big|h+\frac{x^{\perp}}{2(t_0-t)}\big|^2
+\frac12 (t_0-t)^{-\kappa}|\nabla(\eta{\tilde d}^2)|^2{\hat\rho}.
\end{split}
\label{pa30}
\end{equation}
Here we also used \eqref{pa23.5}. 
The terms involving $\nabla\eta$ is non-zero only on $B_1\setminus B_{\frac34}$ where
$|{\tilde d}|\leq c\beta_*$ by \eqref{thap4} and $(t_0-t)^{-1}{\hat\rho}$ is uniformly bounded. Thus, using also \eqref{pa20.5}
and \eqref{pa20},
we obtain from \eqref{pa30} that
\begin{equation}
I_1+I_2\leq\int_{t_1}^{t_2}\int_{B_1} -\frac{{\hat\rho}\eta{\tilde d}^2}{t_0-t}+\frac{\rho}{2}\big|h+\frac{x^{\perp}}{2(t_0-t)}\big|^2
+4(t_0-t)^{-\kappa}{\hat\rho}\eta{\tilde d}^2\, d\|V_t\|dt+c\beta_*^2.
\label{pa31}
\end{equation}
{\bf Estimation of $I_3$}.
\newline
We separate integration with respect to the spatial variable $x$ into the region $A_1=\{|x|\leq (t_0-t)^{\frac{\kappa}{2}}\}\cap B_1$ and its complement 
$A_2=B_1\setminus A_1$. On $A_1$, ${\tilde d}(x)\leq 2 {\rm dist}\, (x,J)\leq 2(t_0-t)^{\frac{\kappa}{2}}$
by \eqref{pa20} so recalling that ${\hat\rho}=(t_0-t)^{-\kappa}g\rho$ and \eqref{pa23.5}, we see that
\begin{equation}
\int_{t_1}^{t_2}\int_{A_1}{\hat\rho}\eta{\tilde d}^2 |u|^2\, d\|V_t\|dt\leq \int_{t_1}^{t_2}\int_{B_1}4\rho\eta|u|^2\, d\|V_t\|dt \leq c(p,q,E_1)\|u\|^2
\label{pa32}
\end{equation}
where the last inequality follows from \eqref{pa19}. On $A_2$, we have
\begin{equation}
{\hat\rho}\leq (4\pi)^{-\frac12}(t_0-t)^{-\kappa-\frac12}\exp\left(-\frac{1}{4(t_0-t)^{1-\kappa}}\right)
\label{pa33}
\end{equation}
and hence ${\hat\rho}$ is uniformly bounded for $t\in [1, t_0)$. Thus 
$$\int_{t_{1}}^{t_{2}}\int_{A_{2}} {\hat\rho}\eta{\tilde d}^2 |u|^2\, d\|V_t\|dt \leq c(\kappa,p,q,E_1)\|u\|^2$$ and we conclude that
\begin{equation}
I_3\leq c(\kappa,p,q,E_1)\|u\|^2.
\label{pa34}
\end{equation}
{\bf Estimation of $I_4$}.
\newline
Since 
\begin{equation}
\int_{t_1}^{t_2}\int_{B_1}{\hat\rho}|\nabla(\eta{\tilde d}^2)||u|\, d\|V_t\|dt
\leq\frac12\int_{t_1}^{t_2}\int_{B_1}(t_0-t)^{-\kappa}{\hat\rho}|\nabla(\eta{\tilde d}^2)|^2
+\rho|u|^2\, d\|V_t\|dt,
\label{pa35}
\end{equation}
and the term involving $\nabla\eta$ in \eqref{pa35} can be estimated by $c\beta_*^2$, 
and the last term may be estimated as in \eqref{pa19}, we obtain
\begin{equation}
I_4\leq \int_{t_1}^{t_2}\int_{B_1}
4(t_0-t)^{-\kappa}{\hat\rho}\eta{\tilde d}^2\, d\|V_t\|dt
+c\beta_*^2+c(p,q,E_1)\|u\|^2.
\label{pa36}
\end{equation}

\noindent
{\bf  Computation of $I_5$}.
\newline
Now we make the explicit choice of $g$ given by 
\begin{equation}
g(t)=\exp\left(-8\int_{t_1}^t(t_0-s)^{-\kappa}\, ds\right) \quad \forall t \in [t_{1}, t_{0}),
\label{pa37}
\end{equation}
and note that $g(t) \leq 1$ and also, since $\kappa<1$, that
\begin{equation}
\inf_{t\in [t_1 ,t_0)}g(t)\geq c(\kappa)>0.
\label{pa37.5}
\end{equation}
We emphasize that $c(\kappa)$ here may be chosen independently of $t_1\in [1, \frac54]$ and $t_0\in[\frac32,3]$. 
With this choice of $g$, 
\begin{equation}
\rho\frac{d}{dt}\left((t_0-t)^{-\kappa}g\right)=\frac{\kappa{\hat\rho}}{t_0-t}
-8(t_0-t)^{-\kappa}{\hat\rho}
\label{pa38}
\end{equation}
so we have
\begin{equation}
I_5= \int_{t_1}^{t_2}\int_{B_1}
 \frac{\kappa{\hat\rho} \eta{\tilde d}^2}{t_0-t}
-8(t_0-t)^{-\kappa}{\hat\rho} \eta{\tilde d}^2\, d\|V_t\|dt.
\label{pa39}
\end{equation}
{\bf Conclusion}.
\newline
By combining \eqref{pa31}, \eqref{pa34}, \eqref{pa36} and 
\eqref{pa39}, and using the estimates \eqref{pa2},
\eqref{pa23}, \eqref{pa37.5} and \eqref{pa20}, we deduce the desired estimate \eqref{pa21}.
\hfill{$\Box$}
\section{A priori estimates III: Distance and approximate continuity estimates for the junction point}

Fix $p, q$ as in (A0), $\nu \in (0, 1), E_1 \in [1, \infty)$ and $\kappa\in [0,1).$ Let  $\Cr{c-2}$ be as in Proposition~\ref{hde}, $\Cr{e-2},\Cr{c-4}$ be as in Proposition~\ref{pade},  and 
$\Cr{c-5}$ be as in Proposition~\ref{tildest}. With $U=B_{3}$ and $I = [0, 4]$, suppose $\{V_t\}_{t\in[0,4]}$ and $\{u(\cdot, t)\}_{t\in[0,4]}$ 
satisfy (A1)-(A4). Assume further, with the following new definitions of $\mu$ and $\|u\|$ (in which  spatial integration is over $B_3$ rather than over $B_2$), that 
\begin{equation}
\mu=\left(\int_0^4 \int_{B_3}{\rm dist}\,(\cdot, J)^2\, d\|V_t\|dt\right)^{\frac12}\leq \frac{\Cr{e-2}}{2},
\label{smep1c}
\end{equation}
\begin{equation}
\|u\|=\left(\int_0^4 \big(\int_{B_3}|u|^p\, d\|V_t\|\big)^{\frac{q}{p}}dt\right)^{\frac{1}{q}}\leq \frac{\Cr{e-2}}{2},
\label{smep2c}
\end{equation}
\begin{equation}
\|V_0\|(\phi_{j_1})\leq \frac{3-\nu}{2} \Cr{c-p}, \ \ \|V_4\|(\phi_{j_2})\geq \frac{1+\nu}{2} \Cr{c-p} \quad \mbox {for some} \quad j_{1}, j_{2} \in \{1, 2, 3\} 
\label{smep34c}
\end{equation}
where $\phi_j$ are as in \eqref{phij}.
Note that \eqref{smep1c}-\eqref{smep34c} are more restrictive conditions than hypotheses (\ref{smep1})-(\ref{smep4}) of Proposition \ref{pade}. 
For $\xi\in {\mathbb R}^2$ with $|\xi|<1$, let $V_t^{(\xi)}$ be the translation of $V_t$ by $-\xi$; thus, 
\begin{equation}
V_t^{(\xi)}(\phi)=\int_{G_1(B_3)}\phi(x-\xi,S)\, dV_t(x,S) \quad \mbox{for each} \quad \phi\in C_c(G_1(B_2)).
\label{ya7}
\end{equation} 
Note that if ${\rm spt}\, \|V_t\|$ has a junction point at $\xi$, then 
${\rm spt}\, \|V_t^{(\xi)}\|$ has a junction point at the origin. Now, depending only on
$\nu,E_1,\Cr{e-2}$ (hence ultimately only on  $p,q,\nu,E_1$), 
there exists a small $\Cl[de]{delta-1} \in (0, 1)$ such that, if $|\xi|\leq\Cr{delta-1}$, then $\{V_t^{(\xi)}\}_{t\in[0,4]}$
and $\{u(\cdot+\xi,t)\}_{t\in [0,4]}$ satisfy \eqref{smep1}-\eqref{smep4}  with $\Cr{e-2}$ in place of 
$\Cr{e-1}$.  In particular,  if $|\xi|\leq \Cr{delta-1}$,
$\{V_t^{(\xi)}\}$ and $\{u(\cdot+\xi,t)\}$ satisfy the hypotheses of Proposition \ref{pade} on $B_2\times [0,4].$ Let us fix such $\Cr{delta-1}$. 
\begin{lemma}
For $\xi\in {\mathbb R}^2\setminus\{0\}$ define $J_{\xi}=\{x\in {\mathbb R}^2 : x-\xi\in J\}$. 
Then on one of the three connected components of $\{x\in{\mathbb R}^2 : {\rm dist}\,(x,{\bf R}_{\frac{\pi}{3}}(J))>|\xi|\}$, we
have 
\begin{equation}
\frac{\sqrt{3}}{2}|\xi|\leq {\rm dist}\,(x,J)+{\rm dist}\,(x,J_{\xi}).
\label{xie2}
\end{equation}
\label{little}
\end{lemma}
{\it Proof}. First note that ${\bf R}_{\frac{\pi}{3}}(J)\setminus\{0\}$ is the set of points $x$ such that  the closest point to $J$ from $x$ is not unique.
Given $\xi\in {\mathbb R}^2\setminus \{0\}$, let $A=\{x\in{\mathbb R}^2 : {\rm dist}\,(x,{\bf R}_{\frac{\pi}{3}}(J))>|\xi|\}$. 
Since $x\in A$ is away from ${\bf R}_{\frac{\pi}{3}}(J)$ by at least $|\xi|$, the closest points in $J$ and $J_{\xi}$ to $x$ are
both unique. Let $x_{J}\in J$ be the closest point to $x$ in $J$. 
One checks easily that the closest point in $J_{\xi}$ from $x$ is $x_{J}+\xi^{\perp}$, where $\xi^{\perp}$ is $(T_{x_J} J)^{\perp}(\xi)$.
This implies ${\rm dist}\, (x,J)=|x-x_J|$ and ${\rm dist}\, (x,J_{\xi})=|x-x_J-\xi^{\perp}|$. Then the triangle inequality gives
$|\xi^{\perp}|\leq {\rm dist}\, (x,J)+{\rm dist}\,(x,J_{\xi})$. For $\xi$, there is at least one component of $A$ on which 
$|\xi^{\perp}|\geq \frac{\sqrt{3}}{2}|\xi|$ holds. On this component, we have \eqref{xie2}.
\hfill{$\Box$}

\begin{prop}
There exist $\Cl[eps]{e-4}\in (0,\frac{\Cr{e-2}}{2}]$ and $\Cl[c]{c-6}\in (1,\infty)$ depending only on 
$p,q,\nu,E_1$ such that if $\{V_t\}_{t\in [0,4]}$, $\{u(\cdot,t)\}_{t\in [0,4]}$
satisfy (A1)-(A4) with $U=B_3$, $I = [0, 4]$ and \eqref{smep1c}, \eqref{smep2c}, \eqref{smep34c}
 hold with $\Cr{e-4}$ in place of $\frac{\Cr{e-2}}{2}$
then for any $\xi\in B_1$ and $t_0\in [\frac32,3]$ with
$h(V_{t_0},\cdot)\in L^2(\|V_{t_0}\|)$ and $\Theta(\|V_{t_0}\|,\xi)=\frac32$,
we have
\begin{equation}
|\xi|\leq \Cr{c-6}\max\{\mu,\|u\|\}.
\label{xie4}
\end{equation}
In addition, given $\kappa\in [0,1)$, there exists $\Cl[c]{c-7}\in (1,\infty)$ depending only on $\kappa,p,q,\nu,E_1$ 
such that
\begin{equation}
\sup_{t\in [\frac54,t_0)}(t_0-t)^{-\kappa}\int_{B_{\frac34}(\xi)}\rho_{(\xi,t_0)}(\cdot,t)
{\rm dist}(\cdot,J_{\xi})^2\, d\|V_t\|\leq \Cr{c-7}\max\{\mu,\|u\|\}^2.
\label{xie4.5}
\end{equation}
\label{distrip}
\end{prop}
{\it Proof}. Recall $\Cr{delta-1}$ which is fixed before Lemma \ref{little}.
Corresponding to $\kappa=1/2$, let $\Cr{c-5}$ be chosen using
Proposition \ref{tildest}. Fix $r_0\in (0,\frac14]$ by
\begin{equation}
r_0 = \min\left\{\frac14, \frac12\left(\frac{16}{3}32eE_1 \sqrt{16\pi}\Cr{c-5}\right)^{-1}\right\}.
\label{defr}
\end{equation}
Corresponding to $\tau=\min\{\Cr{delta-1}/2,\, r_0/8\}$, fix $\Cr{e-p}$ 
using Proposition \ref{smprop}. We will choose $\Cr{e-4}\in (0, \Cr{e-p}]$ by restricting further in the following.
For any $\xi\in B_1$ and $t_0$ satisfying the assumptions, due to the
choice of $\Cr{e-p}$, the claim of Proposition \ref{smprop} shows that we have $|\xi|\leq 2\tau\leq \Cr{delta-1}$. Then due to the
choice of $\Cr{delta-1}$ (see the discussion before Lemma \ref{little}),
$\{V_t^{(\xi)}\}_{t\in [0,4]}$ and $\{u(\cdot+\xi,t)\}_{t\in [0,4]}$
satisfy assumptions of Proposition \ref{pade} and \ref{tildest}.
Thus we have
\begin{equation}
\sup_{t\in [\frac54,t_0)}(t_0-t)^{-1/2}\int_{B_{\frac34}}
\rho_{(0,t_0)}(\cdot,t)\, {\rm dist}\,(\cdot, J)^2\, d\|V_t^{(\xi)}\|
\leq \Cr{c-5}\max\{\mu_{\xi},\|u(\cdot+\xi)\|\}^2,
\label{xie5}
\end{equation}
where $\mu_{\xi}$ is the corresponding quantities for $V_{t}^{(\xi)}$ with integration over $B_2$ and $\|u(\cdot+\xi)\|$
is integration over $B_2$.
By the definition of $V_t^{(\xi)}$ and ${\rm dist}\,(x,J_{\xi})^2
\leq 2\, {\rm dist}\,(x,J)^2+2|\xi|^2$, we have
\begin{equation}
\mu_{\xi}^2\leq 2\mu^2+ 32E_1|\xi|^2
\label{xie6}
\end{equation}
and 
\begin{equation}
\|u(\cdot+\xi)\|\leq \|u\|,
\label{xie7}
\end{equation}
where $\|u\|$ is integration over $B_3$.
In the time interval $[t_0-2r_0^2,t_0-r_0^2]\subset[\frac54,t_0)$, we choose $t_1$ so that
\begin{equation}
\int_{B_3}{\rm dist}\,(\cdot, J)^2\, d\|V_{t_1}\|\leq \frac{3}{r_0^2} \mu^2,
\label{xie8}
\end{equation}
\begin{equation}
\int_{B_1}|h(V_{t_1},\cdot)|^2\, d\|V_{t_1}\|\leq \frac{3\Cr{c-2}}{r_0^2}\max\{\mu,\|u\|\}^2.
\label{xie8.5}
\end{equation}
Such $t_1$ exists by the definition of $\mu$ and by \eqref{thap10}.
Using $t_1$ in \eqref{xie5} as well as \eqref{xie6}, \eqref{xie7} and recalling the definition of $V_t^{(\xi)}$, we then obtain
\begin{equation}
\int_{B_{\frac34}(\xi)}\rho_{(\xi,t_0)}(\cdot,t_1){\rm dist}(\cdot,J_{\xi})^2\, d\|V_{t_1}\|
\leq \Cr{c-5} (t_0-t_1)^{1/2} \max\{2\mu^2+32E_1|\xi|^2,\|u\|^2\}.
\label{xie9}
\end{equation}
On $B_{r_0}$, we have (since $|\xi|\leq 2\tau\leq \frac{r_0}{4}$ and $2r_0^2\geq t_0-t_1\geq r_0^2$)
\begin{equation}
\rho_{(\xi,t_0)}(x,t_1)\geq \frac{\exp\left(-\frac{|x|^2+|\xi|^2}{2(t_0-t_1)}\right)}{\sqrt{8\pi r_0^2}}
\geq \frac{e^{-1}}{\sqrt{8\pi}r_0}.
\label{xie10}
\end{equation}
Using $t_0-t_1\leq 2r_0^2$ and noting that $r_0\leq \frac14$ and $|\xi|\leq 
\frac{r_0}{4}$ (implying $B_{r_0}\subset B_{\frac34}(\xi)$), we obtain from \eqref{xie9} and
\eqref{xie10}
\begin{equation}
\frac{1}{r_0}\int_{B_{r_0}} {\rm dist}(\cdot,J_{\xi})^2\, d\|V_{t_1}\|\leq e\sqrt{16\pi}\Cr{c-5} r_0
 \max\{2\mu^2+32E_1 |\xi|^2,\|u\|^2\}.
 \label{xie11}
 \end{equation}
Next, we consider the set $\{x\in {\mathbb R}^2\,:\, {\rm dist}(x,{\bf R}_{\frac{\pi}{3}}(J))>
\frac{r_0}{4}\}$. Denote the three connected components of this set by $W_1,W_2,W_3$. 
By the argument in the proof of Proposition \ref{firprop}, by choosing $\Cr{e-4}$ 
depending only on $r_0$  (as in \eqref{xie8}) and $\Cr{c-2}$ (as in \eqref{xie8.5}) 
(which ultimately depend only on $p,q,\nu,E_1$), 
we can ensure that
\begin{equation}
\|V_{t_1}\|(W_j\cap B_{r_0})\geq \frac{r_0}{2},\,\,j=1,2,3.
\label{xie12}
\end{equation}
By Lemma \ref{little} and the fact that $|\xi|\leq r_0/4$, on one of the components, we have
\eqref{xie2}. Without loss of generality, let this component be $W_1$. Then \eqref{xie12} implies
\begin{equation}
|\xi|^2\leq \frac{2}{r_0}\int_{W_1\cap B_{r_0}}|\xi|^2\,
d\|V_{t_1}\|\leq \frac{16}{3r_0}\int_{W_1\cap B_{r_0}} 
\left({\rm dist}(\cdot,J)^2 +{\rm dist}(\cdot,J_{\xi})^2\right)\, d\|V_{t_1}\|.
\label{xie13}
\end{equation}
The first term of the right-hand side of \eqref{xie13} may be estimated by 
an appropriate constant times $\mu^2$ due to \eqref{xie8}. 
For the second term, we use \eqref{xie11} and \eqref{defr} to deduce
\begin{equation}
\frac{16}{3r_0}\int_{B_{r_0}}{\rm dist}(\cdot,J_{\xi})^2\, d\|V_{t_1}\|
\leq \frac{16}{3}2e\sqrt{16\pi}\Cr{c-5} r_0\max \{\mu,\|u\|\}^2
+\frac{|\xi|^2}{2}.
\label{xie14}
\end{equation}
By relegating the last term of \eqref{xie14} to the left-hand side in \eqref{xie13}
and setting an appropriate $\Cr{c-6}$, we obtain the desired estimate \eqref{xie4}.
By applying Proposition \ref{tildest} to 
$V_{t}^{(\xi)}$ and using \eqref{xie6} and \eqref{xie4}, we obtain \eqref{xie4.5}.
\hfill{$\Box$}
\begin{prop} Corresponding to $\gamma\in (0,\frac12),\ \kappa\in (0,1),\ p,q,\nu,E_1$, there exist $\Cl[eps]{e-5}\in (0, \Cr{e-4}]$ 
depending only on $p,q,\nu,E_1,\gamma$ and $\Cl[c]{c-8}\in (1,\infty)$ depending only on $p,q,\nu,
E_1,\gamma,\kappa$ such that 
the following holds: Suppose $\{V_t\}_{t\in [0,4]}$ and $\{u(\cdot,t)\}_{t\in [0,4]}$
satisfy (A1)-(A4) with $U=B_3$ and $I = [0, 4]$. Assume that we have \eqref{smep1c}-\eqref{smep2c}
with $\Cr{e-5}$ in place of $\frac{\Cr{e-2}}{2}$ and \eqref{smep34c}.
Suppose that we have two points $\xi_1,\xi_2\in B_1$ and two 
times $t_1,t_2\in [\frac32,3]$ such that $h(V_{t_j},\cdot)\in L^2(\|V_{t_j}\|)$
and $\Theta(\|V_{t_j}\|,\xi_j)=\frac32$ for $j=1,2$. 
Then we have
\begin{equation}
|\xi_1-\xi_2|\leq  \Cr{c-8}  \big(\max\{\mu,\|u\|\}^{\gamma}+\sqrt{|t_1-t_2|}\big)^{\kappa} \max\{\mu,\|u\|\}.
\label{hoe0.5}
\end{equation}
\label{hoeld}
\end{prop}
{\it Proof}. Assume $t_1\leq t_2$ without loss of generality. Write for simplicity 
\begin{equation}
{\hat \mu}=\max\{\mu,\|u\|\},\,\, {\hat s}={\hat \mu}^{\gamma}+\sqrt{|t_2-t_1|},\,\, {\hat \xi}=\xi_2-\xi_1.
\label{hoe0}
\end{equation}
From Proposition \ref{distrip}, we already know that $|{\hat \xi}|\leq 2\Cr{c-6} {\hat \mu}$. 
We choose 
${\hat t}\in [t_1-2{\hat s}^2,t_1-{\hat s}^2]$ 
so that 
\begin{equation}
\int_{B_1}|h(V_{\hat t},\cdot)|^2\, d\|V_{\hat t}\|\leq \frac{2\Cr{c-2}}{{\hat s}^2}\hat\mu^2.
\label{hoe1}
\end{equation}
This is possible by \eqref{thap10}. Since ${\hat s}\geq {\hat \mu}^{\gamma}$ with
$\gamma<\frac12$, ${\hat s}$ is relatively larger than ${\hat \mu}$ for all sufficiently small ${\hat \mu}$. 
We utilize this in the following. We restrict ${\hat \mu}$ depending only on $\Cr{c-6}$ and $\gamma$ so that 
\begin{equation}
2\Cr{c-6}{\hat \mu}\leq \frac{\hat s}{10} 
\label{hoe2}
\end{equation}
E.g. $\hat\mu \leq (20\Cr{c-6})^{-\frac{1}{1-\gamma}}$ is sufficient.
Consider the set $B_{\hat s}\cap \{x\in {\mathbb R}^2\,:\,
{\rm dist}(x,{\bf R}_{\frac{\pi}{3}}(J))>2\Cr{c-6} {\hat \mu}\}$. Due to \eqref{hoe2}, this set consists of 
three non-empty connected components, denoted by $W_1,W_2,W_3$. We have 
\begin{equation}
\frac{{\hat \mu}^2}{{\hat s}^2}\leq{\hat \mu}^{2-2\gamma}={\hat\mu}^{2-4\gamma}\cdot{\hat\mu}^{2\gamma}
\label{hoe3}
\end{equation}
with ${\hat \mu}^{2-4\gamma}$ chosen as small as one likes (note $\gamma<1/2$). 
Thus, restricting ${\hat \mu}$ depending only on $\gamma$ and $\Cr{c-2}$, we can make sure
using \eqref{hoe1} and \eqref{hoe3} that ${\rm spt}\,\|V_{\hat t}\|$ 
lies $o({\hat \mu}^{\gamma})$-neighborhood of $J$ (using the argument
in Proposition \ref{firprop}) in $B_1$. The same can be said about
${\rm spt}\,\|V_{\hat t}^{(\xi_1)}\|$, since $|\xi_1|\leq\Cr{c-6}{\hat \mu}$. 
In particular, since $\hat s\geq \hat \mu^{\gamma}$, we have
\begin{equation}
\|V_{\hat t}^{(\xi_1)}\|(W_j)\geq \frac{\hat s}{4}
\label{hoe4}
\end{equation}
for $j=1,2,3$ under this restriction on ${\hat \mu}$. Since
$|{\hat \xi}|\leq 2\Cr{c-6}{\hat \mu}$, by Lemma \ref{little}, on one of 
$W_j$, say on $W_1$, we have \eqref{xie2} with ${\hat \xi}$ in place of $\xi$. Thus 
by \eqref{hoe4} and \eqref{xie2} we obtain
\begin{equation}
\begin{split}
|{\hat \xi}|^2  \leq \frac{4}{\hat s}\int_{W_1}|{\hat\xi}|^2\, d\|V_{\hat t}^{(\xi_1)}\|
&\leq \frac{32}{3\hat s}\int_{W_1}{\rm dist}(\cdot,J)^2+{\rm dist}(\cdot,J_{\hat\xi})^2\,
d\|V_{\hat t}^{(\xi_1)}\| \\
&\leq \frac{32}{3\hat s}\sum_{j=1}^2\int_{B_{2\hat s}(\xi_j)}
{\rm dist}(\cdot, J_{\xi_j})^2\, d\|V_{\hat t}\|.
\end{split}
\label{hoe5}
\end{equation}
For each $j=1,2$ and $x\in B_{2\hat s}(\xi_j)$, we have
\begin{equation}
\rho_{(\xi_j,t_j)}(x,{\hat t})=\frac{\exp\big(-\frac{|x-\xi_j|^2}{4(t_j-{\hat t})}\big)
}{\sqrt{4\pi (t_j-{\hat t})}}\geq \frac{\exp(-1)}{\sqrt{12\pi} {\hat s}},
\label{hoe6}
\end{equation}
where we used ${\hat s}^2\leq t_j-{\hat t}\leq 3{\hat s}^2$ which follows
easily from the definition of ${\hat s}$ and ${\hat t}$. By \eqref{hoe5}
and \eqref{hoe6}, we obtain
\begin{equation}
|{\hat \xi}|^2\leq \frac{32 e \sqrt{12\pi}}{3}
\sum_{j=1}^2\int_{B_{2\hat s}(\xi_j)}\rho_{(\xi_j,t_j)}(\cdot,{\hat t})
{\rm dist}(\cdot,J_{\xi_j})^2\, d\|V_{\hat t}\|.
\label{hoe7}
\end{equation}
For $\kappa\in [0,1)$, by Proposition \ref{distrip}, each of the last two integrals
is bounded by $\Cr{c-7}(t_j-{\hat t})^{\kappa}{\hat\mu}^2$. Since $t_j-{\hat t}\leq 3{\hat s}^2$,
by defining $\Cr{c-8}$ appropriately, we obtain the desired estimate \eqref{hoe0.5}.
\hfill{$\Box$}
\section{Blow-up analysis and improvement of the space-time $L^{2}$-distance}\label{blowup-analysis}
\label{blowupsection}
Throughout this section, we prove a sequence of propositions under the following assumptions. 
Suppose that $\{V_t^{(m)}\}_{t\in[0,4]}$ and $\{u^{(m)}(\cdot,t)\}_{t\in [0,4]}$ ($m\in {\mathbb N}$) are arbitrary sequences
satisfying (A1)-(A4) and \eqref{smep34c} with $U=B_3$, $I = [0, 4]$ and with $V_{t}^{(m)}$, $u^{(m)}$ in place of $V_{t}$, $u$. Suppose that we have sequences $\{\mu^{(m)}\}_{m\in {\mathbb N}}$
and $\{\|u^{(m)}\|\}_{m\in {\mathbb N}}$ which converge to 0 (and which will be defined in the next section) with the property
\begin{equation}
\big(\int_0^4\int_{B_3}{\rm dist}\,(\cdot,J)^2\, d\|V_t^{(m)}\|dt\big)^{\frac12}\leq \mu^{(m)},
\label{blow1}
\end{equation}
\begin{equation}
\big(\int_0^4\big(\int_{B_3}|u^{(m)}|^p\, d\|V^{(m)}_t\|\big)^{\frac{q}{p}}dt\big)^{\frac{1}{q}}\leq \|u^{(m)}\|
\label{blow2}
\end{equation}
and
\begin{equation}
\lim_{m\rightarrow\infty}(\mu^{(m)})^{-1}\|u^{(m)}\|=0.
\label{blow3}
\end{equation}
Fix a decreasing sequence $\{\tau_m\}_{m\in {\mathbb N}}\subset (0,\frac12)$ with $\lim_{m\rightarrow\infty}\tau_m=0$ and use
Proposition \ref{smprop}
with $\tau=\tau_m$ to obtain $\Cr{e-p}(m)$ and $\Cr{c-p-1}(m)$ corresponding to $\tau_m$. We choose a subsequence 
so that, after relabelling, $\max\{\mu^{(m)},\|u^{(m)}\|\}\leq \Cr{e-p}(m)$ for all $m\in {\mathbb N}$. Then we can apply Proposition
\ref{smprop} to $\{V_t^{(m)}\}_{t\in[0,4]}$ and $\{u^{(m)}(\cdot,t)\}_{t\in [0,4]}$ with $\tau=\tau_m$. Let $f_j^{(m)}:[\tau_m,2-\tau_m]
\times [\tau_m,4-\tau_m]\rightarrow{\mathbb R}$, $j=1,2,3$, be the resulting functions, satisfying \eqref{smep6} and \eqref{smep7}. For each fixed $\tau\in (0,\frac12)$ and for all $m\in {\mathbb N}$ with
$\tau_m\leq \tau$, note that $f_j^{(m)}$ satisfies \eqref{smep6} with $Q = Q_{\tau}=[\tau,2-\tau]\times[\tau,4-\tau]$ 
and $\Cr{c-p-1} = \Cr{c-p-1}(\tau)$, i.e., 
\begin{equation}
\|f_j^{(m)}\|_{C^{1,\zeta}(Q_{\tau})}\leq \Cr{c-p-1}(\tau)    
\max\{\mu^{(m)},\|u^{(m)}\|\}.
\label{blow4}
\end{equation}
For each $m\in {\mathbb N}$ and $j=1,2,3$, define
\begin{equation}
\tilde f_j^{(m)}=(\mu^{(m)})^{-1}f_j^{(m)}.
\label{blow5}
\end{equation}
By \eqref{blow4}, \eqref{blow5}, \eqref{blow3} and the Ascoli-Arzel\`{a} compactness theorem, $\{\tilde f_j^{(m)}\}_{m\in {\mathbb N}}$ has a 
subsequence which converges locally uniformly on $(0,2)\times(0,4)$  to some limit function $\tilde f_j$, $j=1,2,3$.
We also have the estimate
\begin{equation}
\|\tilde f_j\|_{C^{1,\zeta}(Q_{\tau})} \leq \Cr{c-p-1}(\tau).
\label{blow6}
\end{equation}
In the following, we denote subsequences by the same index. 
\begin{prop}
The function $\tilde f_j$ belongs to $C^{\infty}((0,1)\times(1,3))$ and satisfies 
the heat equation $\frac{\partial \tilde f_j}{\partial t}=\frac{\partial^2 \tilde f_j}{\partial x^2}$ on 
$(0,1)\times (1,3)$ for $j=1,2,3$. 
\label{blowprop1}
\end{prop}
{\it Proof}. It is enough to prove the claim for $\tilde f_1$ since the proof for the other two is
similarly carried out after suitable rotations. 
Fix $\phi\in C^{\infty}_c((0,1)\times(1,3);{\mathbb R}^+)$, and fix $\tau\in (0,\frac12)$ so that 
${\rm spt}\, \phi\subset Q_{\tau}$. For all sufficiently large $m$, we have $\tau_m<\tau$ and we 
only consider such $m$. Let $\Cr{c-p-1}=\Cr{c-p-1}(\tau)$ be a constant to be fixed depending only on $\tau$. We take 
in \eqref{meq} $\phi^{(m)}(x,t) \equiv (x_2+2\Cr{c-p-1}    \mu^{(m)})\phi(x_1,t)\eta^{(m)}(x_2)$ as a test function,
where, $x=(x_1,x_2)$ and $\eta^{(m)}$ is a $C^{\infty}$ function such that $\eta=1$ for 
$x_2\in [-\Cr{c-p-1}    \mu^{(m)},\Cr{c-p-1}    \mu^{(m)}]$,
$\eta^{(m)}=0$ for $x_2\notin [-2\Cr{c-p-1}    \mu^{(m)},2\Cr{c-p-1}    \mu^{(m)}]$ 
and $0\leq \eta^{(m)}\leq 1$. 
Note that $\phi^{(m)}(\cdot,t)$ has compact support in $B_3$ and is non-negative, so is a valid choice as a test function. Since 
$x_2=f_1^{(m)}(x_1,t)$ for $(x_1,x_2)\in {\rm spt}\,\|V_t^{(m)}\|$, and since $|f^{(m)}|
\leq \Cr{c-p-1}    \mu^{(m)}$ by \eqref{smep6} and \eqref{blow3}, for all sufficiently large $m$, 
we have $\eta^{(m)}=1$ on ${\rm spt}\, \|V_t^{(m)}\|$. Thus in the following computation, even though we need 
$\eta^{(m)}$ for $\phi^{(m)}$ to have non-negativity, we ignore $\eta^{(m)}$. We then have
\begin{equation}
0\leq \int_1^3\int_{B_2} (-h(V_t^{(m)},\cdot)\phi^{(m)}+\nabla\phi^{(m)})\cdot(h(V_t^{(m)},\cdot)+(u^{(m)})^{\perp})
+\frac{\partial\phi^{(m)}}{\partial t}\, d\|V_t^{(m)}\|dt.
\label{blow7}
\end{equation}
By the Cauchy-Schwarz inequality and by dropping a negative $|h|^2$ term, we obtain from \eqref{blow7}
\begin{equation}
\begin{split}
0&\leq \int_1^3\int_{B_2}|u^{(m)}|^2\phi^{(m)}+|u^{(m)}||\nabla\phi^{(m)}|+\frac{\partial\phi^{(m)}}{\partial t}
+\nabla\phi^{(m)}\cdot h(V_t^{(m)},\cdot)\, d\|V_t^{(m)}\|dt\\
&=I_1^{(m)}+I_2^{(m)}+I_3^{(m)}+I_4^{(m)}, \quad {\rm say}.
\end{split}
\label{blow8}
\end{equation}
We next estimate $\lim_{m\rightarrow\infty}(\mu^{(m)})^{-1}I_j^{(m)}$. By the H\"{o}lder inequality, we have
\begin{equation}
\lim_{m\rightarrow\infty}(\mu^{(m)})^{-1}I_1^{(m)}\leq \lim_{m\rightarrow\infty}c(p,q)(\mu^{(m)})^{-1}\|u^{(m)}\|^2=0,
\label{blow9}
\end{equation}
where we used \eqref{blow3}. Similarly, since $|\nabla\phi^{(m)}|\leq c(\phi)$ on ${\rm spt}\, \|V_t^{(m)}\|$
and by \eqref{blow3},
\begin{equation}
\lim_{m\rightarrow\infty} (\mu^{(m)})^{-1}I_2^{(m)}\leq \lim_{m\rightarrow\infty}c(\phi)(\mu^{(m)})^{-1}\|u^{(m)}\|=0.
\label{blow10}
\end{equation}
For $I_3^{(m)}$, we have
\begin{equation}
\begin{split}
\lim_{m\rightarrow\infty}(\mu^{(m)})^{-1}I_3^{(m)}&=\lim_{m\rightarrow\infty}\int_1^3\int_0^1(\tilde f_1^{(m)}+2\Cr{c-p-1})
\frac{\partial\phi}{\partial t}\, \sqrt{1+|\partial_{x_1} f_1^{(m)}|^2}\, dx_1dt\\
&=\int_1^3\int_0^1 (\tilde f_1+2\Cr{c-p-1})
\frac{\partial\phi}{\partial t}\, dx_1dt,
\end{split}
\label{blow11}
\end{equation}
where we used the uniform convergence $\tilde f_1^{(m)}\rightarrow\tilde f_1$ and \eqref{blow4}.
For $I_4^{(m)}$, since $\nabla\phi^{(m)}=(0,1)\phi+(x_2+2\Cr{c-p-1}    \mu^{(m)})\nabla\phi$, writing
$h(V_t^{(m)},x)=h=(h_1,h_2)$, we have
\begin{equation}
(\mu^{(m)})^{-1}I_4^{(m)}=(\mu^{(m)})^{-1}\int_1^3\int_{B_2} \{\phi h_2+(x_2+2\Cr{c-p-1}    \mu^{(m)})\nabla\phi\cdot h\}\, d\|V_t^{(m)}\|dt.
\label{blow12}
\end{equation}
Since $|x_2+2\Cr{c-p-1}    \mu^{(m)}|\leq 3\Cr{c-p-1}    \mu^{(m)}$ and using the estimate \eqref{thap10}
which is valid here, the second term of the integral converges to $0$. For the first term, by the first variation formula,
\begin{equation}
\int_{B_2}\phi h_2\, d\|V_t^{(m)}\|=-\int_{B_2} S\cdot (\nabla\phi\otimes (0,1))\, dV_t^{(m)}(\cdot,S).
\label{blow13}
\end{equation}
Since $S=(1+|\partial_{x_1} f^{(m)}_1|^2)^{-1}(1,\partial_{x_1} f^{(m)}_1)\otimes(1,\partial_{x_1} f^{(m)}_1)$ and $\nabla\phi
=(\partial_{x_1}\phi,0)$, we have
\begin{equation}
-\int_{B_2} S\cdot (\nabla\phi\otimes (0,1))\, dV_t^{(m)}(\cdot,S)=
-\int_0^1 (1+|\partial_{x_1}f_1^{(m)}|^2)^{-\frac12}\partial_{x_1}f_1^{(m)}\partial_{x_1}\phi\,dx_1.
\label{blow14}
\end{equation}
Since $\nabla \tilde f_1^{(m)}\rightarrow \nabla\tilde f_1$ uniformly, \eqref{blow12}-\eqref{blow14} show that
\begin{equation}
\lim_{m\rightarrow\infty} (\mu^{(m)})^{-1}I_4^{(m)}=-\int_1^3 \int_0^1 \partial_{x_1}\tilde f_1\partial_{x_1}\phi\, dx_1.
\label{blow15}
\end{equation}
Combining \eqref{blow8}-\eqref{blow11} and \eqref{blow15}, we obtain (writing $x_1$ as $x$)
\begin{equation}
0\leq \int_1^3\int_0^1 (\tilde f_1+2\Cr{c-p-1})\frac{\partial\phi}{\partial t}-
\frac{\partial \tilde f_1}{\partial x}\frac{\partial\phi}{\partial x}\, dxdt.
\label{blow16}
\end{equation}
We carry out the same argument with $\phi^{(m)}=(2\Cr{c-p-1}\mu^{(m)}-x_2)\phi(x_1,t)\eta^{(m)}(x_2)$, which
is again non-negative with compact support. The limit in this case produces
\begin{equation}
0\leq \int_1^3\int_0^1 (2\Cr{c-p-1}-\tilde f_1)\frac{\partial\phi}{\partial t}+\frac{\partial \tilde f_1}{\partial x}\frac{\partial \phi}{\partial x}\,
dxdt.
\label{blow17}
\end{equation}
Since $\phi$ has a compact support in $(0,1)\times (1,3)$, the term involving 
$\Cr{c-p-1}$ is $0$. Thus \eqref{blow16} and \eqref{blow17} give
\begin{equation}
0=\int_1^3\int_0^1 \tilde f_1\frac{\partial\phi}{\partial t}-\frac{\partial \tilde f_1}{\partial x}\frac{\partial \phi}{\partial x}\,
dxdt.
\label{blow18}
\end{equation}
We have proved that \eqref{blow18} holds for arbitrary $\phi\in C^{\infty}_c((0,1)\times(1,3);{\mathbb R}^+)$. One can 
then prove that \eqref{blow18} holds for $\phi\in C^{\infty}_c((0,1)\times(1,3))$ which is not necessarily non-negative.
By the standard regularity theory of parabolic equation, $\tilde f_1$ is smooth and satisfies the heat equation.
\hfill{$\Box$}

For the sequence $\{V_t^{(m)}\}_{t\in [0,4]}$ ($m\in{\mathbb N}$) under consideration, we define the following.
\begin{define}
\begin{equation}
\begin{split}
T_g=\cap_{m\in{\mathbb N}}\Big\{t\in [\frac32,3]\, :&\,  V_t^{(m)} \mbox{is a unit density 1-varifold},\ \ 
h(V_t^{(m)},\cdot)\in L_{loc}^2(\|V_t^{(m)}\|)\\
&\mbox{ and } \Theta(\|V_t^{(m)}\|,x)=1\mbox{ or }\frac32,\, \ \forall x\in {\rm spt}\,\|V_t^{(m)}\|\Big\},
\end{split}
\label{trisin1}
\end{equation}
where all the conditions are required  to be satisfied in $B_3$.
\end{define}
Since above conditions are satisfied for a.e$.$ $t\in [\frac12,4]$ for each $\{V_t^{(m)}\}_{t\in [0,4]}$, $T_g$ is a 
full measure set in $[\frac32,3]$, i.e., ${\mathcal L}^1([\frac32,3]\setminus T_g)=0$. 
We next define the following sets.
\begin{define}
For $m\in{\mathbb N}$ and $t\in T_g$, define 
\begin{equation}
\xi^{(m)}(t)=\Big\{x\in B_1\, :\, \Theta(\|V_t^{(m)}\|,x)=\frac32\Big\},\ \ \tilde\xi^{(m)}(t)=\Big\{\frac{x}{\mu^{(m)}}\ : \ x\in \xi^{(m)}(t)
\Big\}.
\label{trisin2}
\end{equation}
\end{define}
Since ${\rm spt}\, \|V_t^{(m)}\|$ away from the origin consists of three $C^{1,\zeta}$ curves, there has to be at least one point $x$
in $B_1$ with $\Theta(\|V_t^{(m)}\|,x)=\frac32$. Otherwise $\Theta(\|V_t\|.x)=1$ $\forall x\in {\rm spt}\, \|V_t^{(m)}\|\cap B_1$
and ${\rm spt}\, \|V_t^{(m)}\|$ has to be a union of regular embedded curves, a contradiction. Thus $\xi^{(m)}(t)\neq \emptyset$
for all $t\in T_g$.
We now apply Proposition \ref{distrip} and \ref{hoeld}. Fix $\gamma\in (0,\frac12)$ and $\kappa\in (0,1)$. Then for all sufficiently
large $m$, \eqref{xie4} and \eqref{hoe0.5} combined with \eqref{blow3} show that for any $a\in \tilde\xi^{(m)}(t)$, 
\begin{equation}
|a|\leq \Cr{c-6},
\label{trisin3}
\end{equation}
and for any $a\in \tilde\xi^{(m)}(t_1)$ and $b\in \tilde\xi^{(m)}(t_2)$, 
\begin{equation}
|a-b|\leq \Cr{c-8}((\mu^{(m)})^{\gamma}+\sqrt{|t_1-t_2|})^{\kappa}.
\label{trisin4}
\end{equation}
\begin{prop}
There exists a $\frac{\kappa}{2}$-H\"{o}lder continuous 
function $\tilde\xi:[\frac32,3]\rightarrow{\mathbb R}^2$ such that 
\begin{equation}
\sup_{t\in [\frac32,3]}|\tilde\xi(t)|\leq \Cr{c-6},
\label{trisin5}
\end{equation}
\begin{equation}
\sup_{t_1,t_2\in [\frac32,3],\, t_1\neq t_2}\frac{ |\tilde\xi(t_1)-\tilde\xi(t_2)|}{|t_1-t_2|^{\frac{\kappa}{2}}}\leq \Cr{c-8}
\label{trisin6}
\end{equation}
and $\tilde\xi^{(m)}(t)$ converges uniformly on $T_g$ to $\tilde\xi(t)$ in the Hausdorff distance.
\label{blowprop2}
\end{prop}
{\it Proof}. Choose a countable dense set $\{t_i\}_{i\in {\mathbb N}}\subset T_g$. For all sufficiently large $m$, 
$\tilde\xi^{(m)}(t)$ is bounded uniformly by \eqref{trisin3}. Also by \eqref{trisin4}, one notes
that the diameter of $\tilde\xi^{(m)}(t)$ is $\leq \Cr{c-8}(\mu^{(m)})^{\gamma}$ and 
converges to $0$ as $m\rightarrow\infty$. Thus for each fixed $i\in {\mathbb N}$,
$\{\tilde\xi^{(m)}(t_i)\}_{m\in {\mathbb N}}$ has a converging subsequence whose limit is a single point. 
By the diagonal argument, we can extract
a subsequence (denoted by the same index) so that $\{\tilde\xi^{(m)}(t_i)\}_{m\in {\mathbb N}}$ converges to a limit point
denoted by $\tilde\xi(t_i)$. By \eqref{trisin4}, $\tilde\xi$ is $\frac{\kappa}{2}$-H\"{o}lder continuous on this countable set,
and one can extend the definition of $\tilde\xi$ uniquely to the whole $[\frac32,3]$ with the same H\"{o}lder constant.
For $t\in T_g\setminus \{t_i\}_{i\in \mathbb N}$, by using \eqref{trisin4}, one
can prove that $\{\tilde\xi^{(m)}(t)\}_{m\in{\mathbb N}}$ also converges to $\tilde\xi(t)$ and that the convergence is uniform. 
\hfill{$\Box$}
\begin{define}
For each $t\in [\frac32,3]$ and $j=1,2,3$, let $\tilde\xi_j^{\perp}(t)\in {\mathbb R}$ be obtained as follows. 
For $j=1$, set $\tilde\xi_1^{\perp}(t)$ be the second coordinate of 
$\tilde\xi(t)$. For $j=2,3$, rotate $\tilde\xi(t)$ by $\frac{2\pi(j-1)}{3}$ 
clockwise, and take its second coordinate to be $\tilde\xi_j^{\perp}(t)$. 
\label{defcp}
\end{define}
Since ${\bf R}_0+{\bf R}_{-\frac{2\pi}{3}}+{\bf R}_{-\frac{4\pi}{3}}={\bf 0}$, we have
\begin{equation}
\tilde\xi_1^{\perp}(t)+\tilde\xi_2^{\perp}(t)+\tilde\xi_3^{\perp}(t)=0
\label{trisin7}
\end{equation}
for all $t\in [\frac32,3]$. 
\begin{prop}
\label{blowprop3}
For any $t_0\in [\frac32,3]$, we have
\begin{equation}
\sup_{t\in [\frac54,t_0)} (t_0-t)^{-\kappa-\frac12}\int_{0}^{\frac12} e^{-\frac{x^2}{4(t_0-t)}}
\sum_{j=1}^3 |\tilde f_j(x,t)-\tilde\xi_j^{\perp}(t_0)|^2\, dx\leq \sqrt{4\pi}\,\Cr{c-7}.
\label{trisin8}
\end{equation}
\end{prop}
{\it Proof}. If we prove 
\begin{equation}
(t_0-t)^{-\kappa-\frac12}\int_{\tau}^{\frac12}e^{-\frac{x^2}{4(t_0-t)}} \sum_{j=1}^3 |\tilde f_j(x,t)-\tilde\xi_j^{\perp}(t_0)|^2\, dx
\leq \sqrt{4\pi}\, \Cr{c-7},
\label{trisin9}
\end{equation}
for arbitrary $t_0\in T_g$, $t\in [\frac54,t_0)$ and $\tau\in (0,\frac12)$,
then by the continuity of $\tilde\xi_j^{\perp}$, \eqref{trisin9} is true for all 
$t_0\in [\frac32,3]$ and we will end the proof of \eqref{trisin8}. Thus we 
fix arbitrary such $t_0,t,\tau$. By \eqref{trisin1}, there exists a sequence of non-empty sets
$\{\xi^{(m)}(t_0)\}_{m\in {\mathbb N}}$ as in \eqref{trisin2}. From each $\xi^{(m)}(t_0)$,
choose one point $\xi_*^{(m)}(t_0)\in \xi^{(m)}(t_0)$. Now, for all sufficiently large $m$,
we may apply Proposition \ref{distrip} with $\xi$ there replaced by $\xi_*^{(m)}(t_0)$. 
Thus for all sufficiently large $m$, we have
\begin{equation}
(t_0-t)^{-\kappa}\int_{B_{\frac58}\setminus B_{\frac{\tau}{2}}} \rho_{(\xi_*^{(m)}(t_0),t_0)}(\cdot,t){\rm dist}\,(\cdot, J_{\xi_*^{(m)}(t_0)})^2
\, d\|V_t^{(m)}\|\leq \Cr{c-7}(\mu^{(m)})^2.
\label{trisin10}
\end{equation}
As we have seen already, we may represent ${\rm spt}\, \|V_t^{(m)}\|$ by $f_j^{(m)}$ in the relevant domain
after suitable rotations. At a point in $x\in {\rm spt}\, \|V_t^{(m)}\|$ represented by $(s,f_j^{(m)}(s,t))$ after a rotation, 
\begin{equation}
{\rm dist}\, (x,J_{\xi_*^{(m)}(t_0)})=|f_j^{(m)}(s,t) - \xi_j^{(m)\perp}(t_0)|,
\label{trisin11}
\end{equation}
where $\xi_j^{(m)\perp}(t_0)=\mbox{second coordinate of }{\bf R}_{-\frac{2(j-1)\pi}{3}}(\xi_*^{(m)}(t_0))$. Since
$(\mu^{(m)})^{-1}\xi^{(m)}_*(t_0)\rightarrow\tilde \xi(t_0)$, we have $(\mu^{(m)})^{-1} \xi_j^{(m)\perp}(t_0)
\rightarrow \tilde\xi_j^{\perp}(t_0)$. Since $(\mu^{(m)})^{-1}f_j^{(m)}\rightarrow \tilde f_j$ uniformly away from
the origin, with \eqref{trisin11}, we have
\begin{equation}
\lim_{m\rightarrow\infty}(\mu^{(m)})^{-2} \int_{B_{\frac58}\setminus B_{\frac{\tau}{2}}}\rho_{(\xi_*^{(m)}(t_0),t_0)}{\rm dist}\,(\cdot,J_{\xi_*^{(m)}(t_0)})^2\, d\|V_t^{(m)}\|= \sum_{j=1}^3 \int_{\frac{\tau}{2}}^{\frac58} \rho_{(0,t_0)}|\tilde f_j -\tilde\xi_j^{\perp}(t_0)|^2\, dx.
\label{trisin12}
\end{equation}
Recalling the definition of $\rho_{(0,t_0)}$, \eqref{trisin10} and \eqref{trisin12} prove \eqref{trisin9} and we end the proof.
\hfill{$\Box$}
\begin{lemma}
\label{blowlem1}
There exists  $\Cl[c]{c-9}\in (1,\infty)$ depending only on $\kappa,p,q,\nu,E_1$ with the property that
\begin{equation}
\sup_{t\in [\frac54,\frac52]}\left( \int_0^{\frac12}|\tilde f_j(x,t)|^2\, dx\right)^{\frac12}\leq \Cr{c-9}
\label{trisin13}
\end{equation}
for $j=1,2,3$.
\end{lemma}
{\it Proof}. We simply choose $t_0=3$ in Proposition \ref{blowprop3}. 
Then, for any $t\in[\frac54,\frac52]$ and $x\in(0,\frac12)$, $t_0-t\leq \frac74$ and
$\frac{x^2}{4(t_0-t)}\leq \frac18$. Moreover, $|\tilde\xi_j^{\perp}(3)|\leq\Cr{c-6}$ by \eqref{trisin5}. Combining these 
facts and with a suitable constant depending only on $\Cr{c-6}$ and $\Cr{c-7}$, and thus ultimately depending
only on $\kappa,p,q,\nu,E_1$, we obtain \eqref{trisin13} from \eqref{trisin8}.
\hfill{$\Box$}
\begin{define}
We define $\tilde f_{AV}\in C^{\infty}((0,1)\times(\frac54,\frac52))$ by
\begin{equation}
\label{defavef}
\tilde f_{AV}(x,t)=\frac13 (\tilde f_1(x,t)+\tilde f_2(x,t)+\tilde f_3(x,t)).
\end{equation}
\end{define}
\begin{prop}
\label{blowprop4}
The odd extension of $\tilde f_{AV}$ with respect to $x$ satisfies the heat equation on $(-1,1)\times(\frac54,\frac52)$
and is $C^{\infty}$ there. In particular, $\tilde f_{AV}(0,t)=0$ for $t\in (\frac54,\frac52)$. 
\end{prop}
{\it Proof}. Fix $\tau\in (0,\frac14)$. In \eqref{trisin8}, we use $t\in (\frac54,\frac52)$  and  $t_0=\frac{\tau^2}{4}+t$. Then
we obtain
\begin{equation}
\int_0^{\tau} e^{-1} \sum_{j=1}^3 |\tilde f_j(x,t)-\tilde\xi_j^{\perp}(t_0)|^2\, dx\leq \sqrt{4\pi}\, \Cr{c-7}
\left(\frac{\tau^2}{4}\right)^{\kappa+\frac12}.
\label{trisin14}
\end{equation}
Using \eqref{trisin7}, \eqref{defavef} and \eqref{trisin14}, we obtain
\begin{equation}
\int_0^{\tau}|\tilde f_{AV}(x,t)|^2\, dx\leq e\sqrt{4\pi} 4^{-\kappa-\frac12}\Cr{c-7} \tau^{2\kappa+1}.
\label{trisin15}
\end{equation}
We continue to denote the odd extension of $\tilde f_{AV}(x,t)$ for $x\in (-1,0)$ by the same notation. 
For $\phi\in C^{\infty}_c((-1,1)\times(\frac54,\frac52))$, we need to prove
\begin{equation}
\int_{Q} \tilde f_{AV}\frac{\partial \phi}{\partial t}+\tilde f_{AV}\frac{\partial^2 \phi}{\partial x^2}\, dxdt=0
\label{trisin16}
\end{equation}
with $Q=(-1,1)\times(\frac54,\frac52)$. Since $\tilde f_{AV}$ is odd with respect to $x$, we only 
need to prove \eqref{trisin16} for odd $\phi$. Let $\eta_{\tau}\in C^{\infty}({\mathbb R})$ be a function such 
that $\eta_{\tau}(x)=1$ for $|x|\geq \tau$, $\eta_{\tau}(x)=0$ for $|x|\leq\frac{\tau}{2}$, $0\leq \eta_{\tau}\leq 1$
and $|\eta_{\tau}'|\leq 4\tau^{-1}$ and $|\eta_{\tau}''|\leq 16\tau^{-2}$. By integration by parts and the fact that $\tilde f_{AV}$
satisfies the heat equation away from $\{x=0\}$, we have
\begin{equation}
\begin{split}
\int_{Q} \tilde f_{AV}\frac{\partial \phi}{\partial t}+\tilde f_{AV}\frac{\partial^2 \phi}{\partial x^2}\, dxdt&=\lim_{\tau\rightarrow 0+}
\int_{Q}\tilde f_{AV}\eta_{\tau}\frac{\partial \phi}{\partial t}+\tilde f_{AV}\eta_{\tau}\frac{\partial^2 \phi}{\partial x^2}\, dxdt\\
&=-\lim_{\tau\rightarrow 0+} \int_{Q\cap\{|x|\leq \tau\}} \big(2\eta_{\tau}'\frac{\partial \phi}{\partial x}+\phi\eta_{\tau}''\big)
\tilde f_{AV}\, dxdt.
\end{split}
\label{trisin17}
\end{equation}
Since $|\phi(x,t)|\leq c(|\nabla\phi|)|x|$ by the oddness of $\phi$, we have $|\phi \eta_{\tau}''|\leq c\tau$. Then by using 
\eqref{trisin15}, we may prove that \eqref{trisin17} is $=0$, which proves \eqref{trisin16}. The standard regularity theory
shows that $\tilde f_{AV}$ is $C^{\infty}$ on $Q$.
\hfill{$\Box$}
\begin{prop}
\label{blowprop5}
For $j,j'\in \{1,2,3\}$, $t\in (\frac54,\frac52)$ and $x\in (-1,0)$, define $(\tilde f_j-\tilde f_{j'})(x,t)=(\tilde f_j-\tilde f_{j'})(-x,t)$.
Then $\tilde f_j-\tilde f_{j'}$ satisfies the heat equation on $(-1,1)\times\left(\frac54,\frac52\right)$ and is $C^{\infty}$ there. 
In particular, we have
$\frac{\partial (\tilde f_j-\tilde f_{j'})}{\partial x}(0,t)=0$ for $t\in \left(\frac54,\frac52\right)$. 
\end{prop}
{\it Proof}. We consider $\tilde f_1-\tilde f_2$ since others can be similarly proved. 
Given $\tau\in (0,\frac14)$, fix an arbitrary $\tau'\in (0,\tau)$. With respect to $\tau'$, we obtain
$\Cr{c-p-1}(\tau')$ by Proposition \ref{smprop}. For all sufficiently large $m$, define
\begin{equation*}
T_g^{(m)}=\left\{t\in T_g\cap \left[\frac54,\frac52\right]\,:\, \int_{B_2}|h(V_t^{(m)},\cdot)|^2\phi_{\rm rad}^2\, d\|V_t^{(m)}\|\leq \Cr{alpha-1}^2\right\},
\end{equation*}
where $\Cr{alpha-1}$ is from Proposition \ref{firprop}. By \eqref{thap10}, note that we have
\begin{equation}
{\mathcal L}^1\left([\frac54,\frac52]\setminus T_g^{(m)}\right)\leq \Cr{alpha-1}^{-2}(\mu^{(m)})^2 \Cr{c-2} .
\label{trisin18}
\end{equation}
For any $t\in T_g^{(m)}$, by Proposition \ref{firprop}, ${\rm spt}\,\|V_t^{(m)}\|$ is represented by 
$f^{(m)}_j$ ($j=1,2,3$) inside $B_1$, and there is only one junction point where
three curves are joined with angle $\frac{2\pi}{3}$. Using the similar notation in the proof of Proposition \ref{firprop}, 
each curves are represented as in \eqref{fir7} with $f_j^{(m)}(\cdot,t)$ and $s_j^{(m)}$ in place of $f_j$ and $s_j$ there. 
We have \eqref{fir8}-\eqref{fir10} satisfied similarly. 
Suppose without loss of generality that $s_1^{(m)}\leq s_2^{(m)}$. For all sufficiently large $m$, we have 
$|s_2^{(m)}|<\tau'$ since $|s_j^{(m)}|\leq \Cr{c-6}\mu^{(m)}$ by \eqref{xie4}. For each $s\in (\tau',\tau)$, 
we have (omitting $t$ dependence)
\begin{equation}
\begin{split}
\big|\frac{\partial( f_1^{(m)}-f_2^{(m)})}{\partial x}(s)\big|&\leq 
\big|\frac{\partial( f_1^{(m)}-f_2^{(m)})}{\partial x}(s_2^{(m)})\big|
+\int_{s_2^{(m)}}^s\sum_{j=1,2} \big|\frac{\partial^2 f_j^{(m)}}{\partial x^2}\big|\, dx \\
&\leq \big|\frac{\partial f_1^{(m)}}{\partial x}(s_2^{(m)})-\frac{\partial f_1^{(m)}}{\partial x}(s_1^{(m)})\big|
+2\int_{B_{2\tau}}|h(V_t^{(m)},\cdot)|\, d\|V_t^{(m)}\|,
\end{split}
\label{trisin19}
\end{equation}
where we have used $\frac{\partial f_1^{(m)}}{\partial x}(s_1^{(m)})=\frac{\partial f_2^{(m)}}{\partial x}(s_2^{(m)})$
which follows from \eqref{fir10}. The first term of \eqref{trisin19} may be bounded by the second term, so we
obtain from \eqref{trisin19}
\begin{equation}
\sup_{s\in (\tau',\tau)}\big|\frac{\partial (f_1^{(m)}-f_2^{(m)})}{\partial x}(s,t)\big|
\leq 4\int_{B_{2\tau}}|h(V_t^{(m)},\cdot)|\, d\|V_t^{(m)}\|.
\label{trisin20}
\end{equation}
For any $t\in [\frac54,\frac52]\setminus T_g^{(m)}$, we have
\begin{equation}
\sup_{s\in (\tau',\tau)}\big|\frac{\partial( f_1^{(m)}-f_2^{(m)})}{\partial x}(s,t)\big|
\leq 2 \Cr{c-p-1}(\tau') \mu^{(m)}.
\label{trisin21}
\end{equation}
Combining \eqref{trisin18}, \eqref{trisin20} and \eqref{trisin21}, we obtain
\begin{equation}
\int_{\frac54}^{\frac52}\sup_{s\in (\tau',\tau)}\big|\frac{\partial (f_1^{(m)}-f_2^{(m)})}{\partial x}(s,t)\big|
\, dt\leq 4\int_{\frac54}^{\frac52}\int_{B_{2\tau}}|h(V_t^{(m)},\cdot)|\, d\|V_t^{(m)}\|+2\Cr{c-p-1}\Cr{alpha-1}^{-2}\Cr{c-2}(\mu^{(m)})^3.
\label{trisin22}
\end{equation}
We may estimate
\begin{equation}
\begin{split}
\int_{\frac54}^{\frac52}\int_{B_{2\tau}}|h(V_t^{(m)},\cdot)|\, d\|V_t^{(m)}\|&\leq (5 \tau E_1)^{\frac12} \big(\int_1^3\int_{B_2}
|h(V_t^{(m)},\cdot)|^2\phi_{\rm rad}^2\, d\|V_t^{(m)}\|\big)^{\frac12}\\
&\leq (5\tau E_1)^{\frac12}\Cr{c-2}^{\frac12}\mu^{(m)},
\end{split}
\label{trisin23}
\end{equation}
by \eqref{thap10}. Thus dividing \eqref{trisin22} by $\mu^{(m)}$ and letting $m\rightarrow\infty$, \eqref{trisin23} and 
the uniform convergence of $\frac{\partial \tilde f_j^{(m)}}{\partial x}$ to $\frac{\partial \tilde f_j}{\partial x}$ show that
\begin{equation}
\int_{\frac54}^{\frac52}\sup_{s\in (\tau',\tau)}\big|\frac{\partial(\tilde f_1-\tilde f_2)}{\partial x}(s,t)\big|\, dt
\leq (5\tau E_1)^{\frac12}\Cr{c-2}^{\frac12}.
\label{trisin24}
\end{equation}
Since the right-hand side of \eqref{trisin24} does not depend on $\tau'$, the same inequality holds with $\tau'=0$. 
Arguing as in the proof for $\tilde f_{AV}$, we may prove that the even extension of $\tilde f_1-\tilde f_2$ with
respect to $x$ now satisfies the heat equation weakly by using \eqref{trisin24} (with $\tau'=0$). 
This time, it is sufficient to use even test functions
$\phi$, and use also $|\frac{\partial \phi}{\partial x}|\leq c(|\nabla^2\phi|)|x|$ in the proof to estimate the truncation error due to $\eta_{\tau}$. 
We omit the detail since it is 
similar to the previous case. In particular, we have $\tilde f_1-\tilde f_2$ is $C^{\infty}$ on $(-1,1)\times(\frac54,\frac52)$.
Since it is evenly extended, the $x$-derivative vanishes on $x=0$. 
\hfill{$\Box$}
\begin{prop}
We have $\tilde f_j\in C^{\infty}([0,1)\times(\frac54,\frac52))$ and $\tilde f_j(0,t)=\tilde\xi_j^{\perp}(t)$
for $t\in (\frac54,\frac52)$ and for $j=1,2,3$. Moreover, for any $k,l\in {\mathbb N}\cup\{0\}$, there exists
$\Cl[c]{c-10}\in (1,\infty)$ depending only on $\kappa,p,q,\nu,E_1,k,l$ such that, for $j=1,2,3$, 
\begin{equation}
\sup_{[0,\frac14)\times(\frac32,\frac52)}\Big|\frac{\partial^{k+l}\tilde f_j}{\partial x^{k}\partial t^{l}}\Big|
\leq \Cr{c-10}.
\label{trisin27}
\end{equation}
\label{blowprop6}
\end{prop}
{\it Proof}. We prove this  for the case $j=1$; the other cases follow by similar reasoning.
Since we may write $\tilde f_1=\tilde f_{AV}+\frac13(\tilde f_1-\tilde f_2)+\frac13(\tilde f_1-\tilde f_3)$,
Proposition \ref{blowprop4} and \ref{blowprop5} show that $\tilde f_1$ may be smoothly extended for $\{x\leq 0\}$.
More precisely, for $x\in (-1,0)$ and $t\in (\frac54,\frac52)$, define 
\begin{equation}
\tilde f_1 (x,t)=\big\{-\tilde f_{AV}+
\frac13(\tilde f_1-\tilde f_2)+\frac13(\tilde f_1-\tilde f_3)\big\}\Big|_{(-x,t)}=\frac13(\tilde f_1-2\tilde f_2-2\tilde f_3)(-x,t).
\label{trisin25}
\end{equation}
Then $\tilde f_1$ is in $C^{\infty}((-1,1)\times(\frac54,\frac42))$ and satisfies the heat equation on its domain. Moreover, by
\eqref{trisin13} and \eqref{trisin25}, we have
\begin{equation}
\sup_{t\in [\frac54,\frac52]}\left(\int_{-\frac12}^{\frac12}|\tilde f_1(x,t)|^2\, dx\right)^{\frac12}\leq \frac83\Cr{c-9},
\label{trisin26}
\end{equation}
and by \eqref{trisin8}, $\tilde f_1(0,t)=\tilde \xi_1^{\perp}(t)$ holds for $t\in (\frac54,\frac52)$. 
Since $\tilde f_1$ satisfies the heat equation with the estimate \eqref{trisin26}, the standard regularity theory 
(\cite{ladyzhenskaja}) shows that any partial derivatives of $\tilde f_1$ on $(-\frac14,\frac14)\times (\frac32,\frac52)$
can be bounded by a constant depending only on $\Cr{c-9}$ and the order of differentiation. Since $\Cr{c-9}$ depends only on $\kappa,p,q,\nu,E_1$,
we have \eqref{trisin27} with a suitable constant $\Cr{c-10}$.
\hfill{$\Box$}
\begin{define}
Recall the definition of $\tilde\xi(t)$ (in Definition \ref{defcp}). For each $m\in {\mathbb N}$, define $\tilde J^{(m)}\subset {\mathbb R}^2$ to be the set obtained by first rotating $J$ counterclockwise by $\arctan(\mu^{(m)}\frac{\partial\tilde f_1}{\partial x}(0,2))$ and then translating by
$\mu^{(m)}\tilde \xi(2)$. In the similarity class of $J$, $\tilde J^{(m)}$ is the element characterized by the properties that it has junction point at $\mu^{(m)}\tilde \xi(2)$, and the slope of  its ray close to the positive $x$-axis is equal to $\mu^{(m)}\frac{\partial\tilde f_1}{\partial x}(0,2)$.
Denote the junction point of ${\tilde J}^{(m)}$ by $a^{(m)}$ (thus $a^{(m)} =\mu^{(m)}{\tilde \xi}(2)$). 
\label{deftp}
\end{define}
To clarify the property of $\tilde J^{(m)}$ concerning the slope of its ray close the the $x$-axis, we recall that the second coordinate of ${\bf R}_{-\frac{(j-1)2\pi}{3}}(\tilde \xi(t))$ has been denoted by 
$\tilde\xi_j^{\perp}(t)$ and is equal to $\tilde f_j(0,t)$ for $j=1,2,3$. The ray of $\tilde J^{(m)}$ close to the
$x$-axis can be expressed as
\begin{equation}
\Big\{\Big(x,\mu^{(m)}\tilde f_1(0,2)+\mu^{(m)}(x-\mu^{(m)}\tilde \xi_1(2))\frac{\partial \tilde f_1}{\partial x}(0,2)\Big)\in {\mathbb R}^2\, :\,
x\in (\mu^{(m)}\tilde \xi_1(2),\infty)\Big\}.
\label{trisin28}
\end{equation}
More generally, for $j=1,2,3$, the half line of ${\bf R}_{-\frac{(j-1)2\pi}{3}}(\tilde J^{(m)})$ close to the $x$-axis is
\begin{equation}
\Big\{\Big(x,\mu^{(m)}\tilde f_j(0,2)+\mu^{(m)}(x-\mu^{(m)}v_j)\frac{\partial \tilde f_j}{\partial x}(0,2)\Big)\in {\mathbb R}^2\, :\,
x\in (\mu^{(m)}v_j,\infty)\Big\},
\label{trisin29}
\end{equation}
where $v_j$ is the first coordinate 
of ${\bf R}_{-\frac{(j-1)2\pi}{3}}(\tilde\xi(2))$. It is important to note that we used
$\frac{\partial \tilde f_1}{\partial x}(0,2)=\frac{\partial \tilde f_j}{\partial x}(0,2)$ ($j=2,3$) which follows from Proposition \ref{blowprop5}
and \ref{blowprop6}. 
\begin{prop}
\label{blowprop7}
There exists $\Cl[c]{c-11}$ depending only on $p,q,\nu,E_1$ such that, for all $\theta\in (0,\frac14)$, we have
\begin{equation}
\limsup_{m\rightarrow\infty} \frac{1}{(\mu^{(m)})^2\theta^5}\int_{2-\theta^2}^{2+\theta^2}\int_{B_{\theta}
(a^{(m)})} {\rm dist}\,(\cdot,\tilde J^{(m)})^2\, d\|V_t^{(m)}\|
dt\leq \Cr{c-11}\theta^{2}
\label{trisin30}
\end{equation}
and
\begin{equation}
d(\tilde J^{(m)},J)\leq \Cr{c-11}\mu^{(m)}.
\label{trisin30.5}
\end{equation}
\end{prop}
{\it Proof}. Fix $\theta\in (0,\frac14)$ and $\tau\in (0,\theta)$. For any $t\in (2-\theta^2,2+\theta^2)$ with $t+\tau^2\in T_g$, choose a point
$\xi_*^{(m)}\in \xi^{(m)}(t+\tau^2)$ (recall \eqref{trisin2}). Then by Proposition \ref{distrip} with
$\kappa=\frac12$ fixed, for all sufficiently large $m$, we have
\begin{equation}
\tau^{-2\kappa}\int_{B_{\frac34}(\xi_*^{(m)})} \rho_{(\xi_*^{(m)},t+\tau^2)}(\cdot,t){\rm dist}\,(\cdot,J_{\xi_*^{(m)}})^2\, d\|V_t^{(m)}\|
\leq \Cr{c-7}(\mu^{(m)})^2.
\label{trisin31}
\end{equation}
Since $\rho_{(\xi_*^{(m)},t+\tau^2)}(x,t)\geq (4\pi\tau^2)^{-\frac12}e^{-1}$ for $|x-\xi_*^{(m)}|\leq 2\tau$, we have from \eqref{trisin31}
\begin{equation}
\int_{B_{\tau}} {\rm dist}(\cdot,J_{\xi_*^{(m)}})^2\, d\|V_t^{(m)}\|\leq \Cr{c-7}(4\pi)^\frac12 e \tau^{1+2\kappa}
(\mu^{(m)})^2.
\label{trisin32}
\end{equation}
Here, the integration should be over $B_{2\tau}(\xi_*^{(m)})$, but since $\xi_*^{(m)}\rightarrow 0$ (uniformly in $t$) as $m\rightarrow \infty$, 
we have $B_{\tau}\subset B_{2\tau}(\xi_*^{(m)})$ for sufficiently large $m$ and we obtain \eqref{trisin32}. We next wish to replace
$J_{\xi_*^{(m)}}$ by $\tilde J^{(m)}$ in \eqref{trisin32}. The (Hausdorff) distance between $J_{\xi_*^{(m)}}$ and $J_{\mu^{(m)}\tilde\xi(2)}$ 
in $B_1$ may be estimated by $3\Cr{c-6}\mu^{(m)}$ for all sufficiently large $m$ due to
$|(\mu^{(m)})^{-1}\xi_*^{(m)}-\tilde \xi(t+\tau^2)|\leq o(1)$  and $|\tilde \xi(t+\tau^2)-\tilde \xi(2)|\leq
2\Cr{c-6}$  by Proposition \ref{blowprop2}. The distance between $J_{\mu^{(m)}\tilde\xi(2)}$ and $\tilde J^{(m)}$ is 
bounded by $\Cr{c-10}\mu^{(m)}$ due to the estimate \eqref{trisin27} (with $k=1$ and $l=0$) 
for the angle of rotation. Thus we have for any $x\in B_{\tau}$
\begin{equation}
{\rm dist}(x,\tilde J^{(m)})^2\leq 
2 \,{\rm dist}(x,J_{\xi_*^{(m)}})^2
+2(\mu^{(m)})^2 (3\Cr{c-6}+\Cr{c-10})^2.
\label{trisin33}
\end{equation}
Now combining \eqref{trisin32} and \eqref{trisin33}, we obtain
\begin{equation}
\limsup_{m\rightarrow\infty}\frac{1}{(\mu^{(m)})^{2}}\int_{B_{\tau}}{\rm dist}\, (\cdot,\tilde J^{(m)})^2\, d\|V_t^{(m)}\|
\leq 2\Cr{c-7}(4\pi)^\frac12 e \tau^{1+2\kappa}+6\tau(3\Cr{c-6}+\Cr{c-10})^2.
\label{trisin34}
\end{equation}
We next estimate the integration over $B_{2\theta}\setminus B_{\tau}$. Fix any $t\in (2-\theta^2,2+\theta^2)$. For all 
sufficiently large $m$ depending on $\tau$, ${\rm spt}\, \|V_t^{(m)}\|\cap B_{2\theta}
\setminus B_{\tau}$ is represented as a union of graphs using $f_j^{(m)}(\cdot,t)=\mu^{(m)}\tilde f_j^{(m)}(\cdot,t)$. 
Recalling \eqref{trisin29}, we have (denoting $\frac{\partial}{\partial x}$ by sub-index $x$ for simplicity)
\begin{equation}
\begin{split}
&\frac{1}{(\mu^{(m)})^2}\int_{B_{2\theta}\setminus B_{\tau}}{\rm dist}(\cdot,\tilde J^{(m)})^2\, d\|V_t^{(m)}\| \\
&\leq \sum_{j=1}^3 \int_{\frac{\tau}{2}} ^{2\theta}
\big|\tilde f^{(m)}_j(x,t)-\tilde f_j(0,2)-(x-\mu^{(m)}v_j) (\tilde f_j)_x(0,2)\big|^2\,
 \sqrt{1+(\mu^{(m)})^2 (\tilde f_j^{(m)})_x^2}\, dx \\
 &\leq \sum_{j=1}^3 \int_{\frac{\tau}{2}} ^{2\theta}
2\big|\tilde f^{(m)}_j(x,t)-\tilde f_j(0,2)-x (\tilde f_j)_x(0,2)\big|^2\, dx+ c(\Cr{c-p-1}(\tau),\Cr{c-10})(\mu^{(m)})^2.
 \end{split}
 \label{trisin35}
 \end{equation}
 We know already that $\tilde f_j^{(m)}$ converges to $\tilde f_j$ on $[\frac{\tau}{2},2\theta]$, and
 \begin{equation}
 |\tilde f_j(x,t)-\tilde f_j(0,2)-x(\tilde f_j)_x (0,2)|\leq \Cr{c-10}(|x|^2+|t-2|)
 \label{trisin36}
 \end{equation}
 by Taylor's theorem and \eqref{trisin27}. Since $|x|\leq2 \theta$ and $|t-2|\leq \theta^2$, \eqref{trisin35}
 and \eqref{trisin36} prove
 \begin{equation}
 \limsup_{m\rightarrow\infty}\frac{1}{(\mu^{(m)})^2} \int_{B_{2\theta}\setminus B_{\tau}}{\rm dist}\,(\cdot, \tilde J^{(m)})^2\, d\|V_t^{(m)}\|
 \leq 24\Cr{c-10}^2 \theta^5.
 \label{trisin37}
 \end{equation}
 Since $\tau$ is arbitrary, combining \eqref{trisin34} and \eqref{trisin37} and setting $\Cr{c-11}\geq 48\Cr{c-10}^2$, we obtain the
 desired estimate \eqref{trisin30}, also by observing $B_{\theta}(a^{(m)})\subset
 B_{2\theta}$ for all sufficiently large $m$. 
 By the definition of $\tilde J^{(m)}$, \eqref{trisin30.5} follows as well with a suitable choice of 
 constant.
 \hfill{$\Box$}
\section{Pointwise estimates: Proof of Theorem~\ref{mainreg}}
With Proposition \ref{blowprop7} established, a standard iteration argument establishes the desired 
estimates as well as the expected geometry of the flow as a regular triple junction moving by curvature. For completeness, we present the detailed argument. 
\begin{prop}
Corresponding to $p,q,\nu,E_1$, there exist $\Cl[eps]{e-6}\in (0,1)$, $\theta_*\in (0,\frac14)$ and $\Cl[c]{c-12}\in (1,\infty)$ such that the following holds: For $R\in (0,\infty)$ and $U=B_{4R}$, suppose $\{V_t\}_{t\in [-2R^2,2R^2]}$ and
$\{u(\cdot,t)\}_{t\in [-2R^2,2R^2]}$ satisfy (A1)-(A4). Assume 
\begin{equation}
\mu=\left( R^{-5}\int_{-2R^2}^{2R^2} \int_{B_{4R}} {\rm dist}\,(\cdot,J)^2 \, d\|V_t\|dt\right)^{\frac12}<\Cr{e-6},
\label{dec3}
\end{equation}
\begin{equation}
\exists j_1,j_2\in \{1,2,3\} : 
R^{-1}\|V_{-2R^2}\|(\phi_{j_1,J,R})\leq \frac{3-\nu}{2}\Cr{c-p}, \ \ R^{-1}\|V_{2R^2}\|(\phi_{j_2,J,R})\geq \frac{1+\nu}{2}\Cr{c-p}, 
\label{dec5}
\end{equation}
and denote
\begin{equation}
\|u\|=R^{\zeta} \left(\int_{-2R^2}^{2R^2}\left(\int_{B_{4R}} |u|^p\, d\|V_t\|\right)^{\frac{q}{p}}\right)^{\frac{1}{q}}.
\label{dec4}
\end{equation}
Then there exists $ J'={\bf R}_{\theta}(J)+\xi\in {\mathcal J}$ such that
\begin{equation}
d_R(J',J)\leq \Cr{c-12}\mu \quad and
\label{dec7}
\end{equation}
\begin{equation}
\left((\theta_* R)^{-5}\int_{-2(\theta_*R)^2}^{2(\theta_*R)^2} \int_{B_{4\theta_* R}(\xi)} {\rm dist}\, (\cdot, J')^2\, d\|V_t\|dt\right)^{\frac12}
\leq \theta_*^{\zeta} \max\{\mu,\Cr{c-12}\|u\|\}.
\label{dec8}
\end{equation}
Moreover, if we additionally assume that $\|u\|\leq \Cr{e-6}$, then we have
\begin{equation}
(\theta_* R)^{-1}\|V_{-2(\theta_*R)^2}\|(\phi_{j, J',\theta_*R})\leq \frac{3-\nu}{2}\Cr{c-p}, \ \ (\theta_* R)^{-1}\|V_{2(\theta_* R)^2}\|
(\phi_{j,J',\theta_*R})\geq \frac{1+\nu}{2}\Cr{c-p}, \ \ j=1,2,3.
\label{dec8.5}
\end{equation}
\label{disjun2}
\end{prop}
{\it Proof}. We may assume that $R=1$ after a parabolic change of variables. 
We prove the claim by contradiction. For all $m\in {\mathbb N}$, consider a set of sequences
$\{V_t^{(m)}\}_{t\in [-2,2]}$ and $\{u^{(m)}(\cdot, t)\}_{t\in [-2,2]}$ satisfying (A1)-(A4) with $U=B_4$
such that \eqref{dec3} and \eqref{dec5} are satisfied with $V_t^{(m)},\frac{1}{m}$, that is, for all $m\in {\mathbb N}$,
\begin{equation}
\mu^{(m)} \equiv \left(\int_{-2}^{2}\int_{B_4}{\rm dist}\,(\cdot, J)^2\, d\|V_t^{(m)}\|dt\right)^{\frac12}<\frac{1}{m},
\label{dec10}
\end{equation}
\begin{equation}
\|V_{-2}^{(m)}\|(\phi_{j})\leq \frac{3-\nu}{2} \Cr{c-p}, \ \ \|V_{2}^{(m)}\|(\phi_{j})\geq \frac{1+\nu}{2}\Cr{c-p}, \ \ j=1,2,3.
\label{dec11}
\end{equation}
Define $\|u^{(m)}\|$ as in \eqref{dec4} with $R=1$ and $u^{(m)}$ in place of $u$. 
The negation then implies that for any $J'={\bf R}_{\theta}(J)+\xi\in {\mathcal J}$ with
\begin{equation}
d(J',J)\leq m\mu^{(m)},
\label{dec13}
\end{equation}
we have
\begin{equation}
\left(\theta_*^{-5}\int_{-2\theta_*^2}^{2\theta_*^2} \int_{B_{4\theta_* }(\xi)} {\rm dist}\, (\cdot,J')^2\, d\|V^{(m)}_t\|dt\right)^{\frac12}
> \theta_*^{\zeta} \max\{\mu^{(m)}, m\|u^{(m)}\|\}.
\label{dec14}
\end{equation}
Here $\theta_*\in (0,\frac14)$ will be chosen depending only on $p,q,E_1$. 

We next proceed to use the argument in the previous section. First,
use $J'=J$ in \eqref{dec14}. 
Then we have
\begin{equation}
\theta_*^{\zeta}\max\{\mu^{(m)},m\|u^{(m)}\|\}<\left(\theta_*^{-5}\int_{-2}^2\int_{B_4}{\rm dist}\, (\cdot,J)^2\, d\|V_t^{(m)}\|dt\right)^{\frac12}
\leq \theta_*^{-\frac52} \mu^{(m)}.
\label{dec15}
\end{equation}
Thus, \eqref{dec15} shows \eqref{blow3} is satisfied. We also have $\lim_{m\rightarrow\infty}\mu^{(m)}=0$ by \eqref{dec10}. 
We shift $t$ by $-2$ so that the time interval will be $[0,4]$ from $[-2,2]$. 
We have \eqref{blow1}, \eqref{blow2} and \eqref{smep34c} satisfied
thus all the assumptions in the previous section are satisfied.  In the argument, we may fix $\kappa=\frac12$ from the beginning. The 
conclusion of Proposition \ref{blowprop7} shows that for all $m>\Cr{c-11}$, $\tilde J^{(m)}$ satisfies \eqref{dec13} due to \eqref{trisin30.5}
while we have \eqref{trisin30}. On the other hand, \eqref{dec14} shows
\begin{equation}
\theta_*^{2\zeta}\leq \limsup_{m\rightarrow\infty} \frac{1}{(\mu^{(m)})^2\theta_*^5} \int_{2-2\theta_*^2}^{2+2\theta_*^2}
\int_{B_{4\theta_*}(a^{(m)})}{\rm dist}\,(\cdot,\tilde J^{(m)})^2\, d\|V_{t-2}^{(m)}\|dt.
\label{dec16}
\end{equation}
If we let $\theta=4\theta_*$ in \eqref{trisin30} and compare \eqref{dec16}, we obtain
\begin{equation}
\theta_*^{2\zeta}\leq 4^7 \Cr{c-11}\theta_*^2.
\label{dec17}
\end{equation}
Note that $\Cr{c-11}$ and $\zeta$ depend only on $p,q,\nu,E_1$. Since $\zeta\in (0,1)$, we obtain a contradiction for 
suitably small $\theta_*$ depending only on $\Cr{c-11},\zeta$ and thus ultimately only on $p,q,\nu,E_1$. Once $\theta_*$
is fixed, then we may use Proposition \ref{smprop} with $\tau=\frac{\theta_*}{2}$.  For suitably small $\Cr{e-6}$, 
we can make sure that ${\rm spt}\|V_t\|$ on $B_{4\theta_*}(\xi)\setminus B_{\frac{\theta_*}{2}}(\xi)$ is close to $J$ (and thus to $J'$)
in $C^{1,\zeta}$ and \eqref{dec8.5} can be guaranteed. 
\hfill{$\Box$}
\begin{prop}
Corresponding to $\nu,E_1,p,q$ there exist $\Cl[eps]{e-7}\in (0,1)$ and $\Cl[c]{c-13}\in (1,\infty)$
with the following property. Under the assumptions of Proposition \ref{disjun2} where
$\Cr{e-6}$ is replaced by $\Cr{e-7}$ and with $\|u\|\leq \Cr{e-7}$, 
\newline
(1) there exists $J_0\in {\mathcal J}$ with the junction point at $\hat a$ such that 
\begin{equation}
d_R(J_{0},J)\leq \Cr{c-13}\max\{\mu,\Cr{c-12}\|u\|\},
\label{del9}
\end{equation}
(2) for $0<s<R$, there exists $J_s\in {\mathcal J}$ with the junction point at $a_s$ such that
\begin{equation}
d_s(J_s,J_0)+\left(s^{-5}\int_{-2s^2}^{2s^2}\int_{B_{4s}(a_s)} {\rm dist}\, (\cdot,J_s)^2\,d\|V_t\|dt\right)^{\frac12}
\leq \Cr{c-13}(s/R)^{\zeta}\max\{\mu,\Cr{c-12}\|u\|\},
\label{del10}
\end{equation}
(3) ${\rm spt}\,\|V_0\|\cap B_{2R}$ consists of three curves meeting at $\hat a$ with 120 degrees. 

\label{disjun3}
\end{prop}
{\it Proof}. We may assume $R=1$ after a change of variables. We choose $\Cr{e-7}\in (0,1)$ and $\Cr{c-13}\in (1,\infty)$ so that
\begin{equation}
\Cr{c-12}\Cr{e-7}<\Cr{e-6},
\label{del4}
\end{equation}
\begin{equation}
(1+\theta_*^{-1})\Cr{c-12}\sum_{j=0}^{\infty}\theta_*^{(j-1)\zeta}\leq \Cr{c-13},
\label{del6}
\end{equation}
\begin{equation}
2\theta_*^{-\frac52 -\zeta}\leq \Cr{c-13}.
\label{del7}
\end{equation}
Set $J^{(0)}=J$ and we inductively prove that for $m=1,2,\cdots$, we have
$J^{(m)}={\bf R}_{\tau^{(m)}}(J)+a^{(m)}\in {\mathcal J}$ such that 
\begin{equation}
d_{\theta_*^{m-1}}(J^{(m)},J^{(m-1)})\leq \Cr{c-12}\theta_*^{(m-1)\zeta}\max\{\mu,\Cr{c-12}\|u\|\},
\label{del1}
\end{equation}
\begin{equation}
\mu^{(m)}=\left(\theta_*^{-5m}\int_{-2\theta_*^{2m}}^{2\theta_*^{2m}} \int_{B_{4\theta_*^m}(a^{(m)})}{\rm dist}\, (\cdot, J^{(m)})^2\, d\|V_t\|dt\right)^{\frac12}
\leq \theta_*^{m\zeta}\max\{\mu,\Cr{c-12}\|u\|\},
\label{del2}
\end{equation}
\begin{equation}
\theta_*^{-m}\|V_{-2\theta_*^{2m}}\|(\phi_{j,J^{(m)},\theta_*^m})\leq \frac{3-\nu}{2} \Cr{c-p}, \ \ \theta_*^{-m}\|V_{2\theta_*^{2m}}\|(\phi_{j,J^{(m)},\theta_*^m})\geq
\frac{1+\nu}{2}\Cr{c-p},\ j=1,2,3.
\label{del3}
\end{equation}
Since $\Cr{e-7}<\Cr{e-6}$ by \eqref{del4}, Proposition \ref{disjun2} gives the proof for the existence of
$J^{(1)}=J'$ satisfying \eqref{del1}-\eqref{del3} for $m=1$ case. Under the inductive 
assumption up to $m$, we have
\begin{equation}
\mu^{(m)}\leq \theta_*^{m\zeta}\max\{\mu,\Cr{c-12}\|u\|\}<\Cr{e-6}
\label{del8}
\end{equation}
by \eqref{dec3}, $\|u\|\leq \Cr{e-7}$ and \eqref{del4}. Then with \eqref{del3} and $J$ 
replaced by $J^{(m)}$ and $R=\theta_*^m$, we have the assumptions \eqref{dec3}
and \eqref{dec5} satisfied. Hence we have a $J^{(m+1)}\in {\mathcal J}$ such that
\begin{equation}
d_{\theta_*^m}(J^{(m+1)},J^{(m)})\leq \Cr{c-12}\mu^{(m)}\leq \Cr{c-12}\theta_*^{m\zeta}
\max\{\mu,\Cr{c-12}\|u\|\}
\label{del12}
\end{equation}
by \eqref{dec7} and \eqref{del2} and
\begin{equation}
\mu^{(m+1)}\leq \theta_*^{\zeta}\max\{\mu^{(m)},\Cr{c-12}\theta_*^{m\zeta}\|u\|\}
\leq \theta_*^{(m+1)\zeta}\max\{\mu,\Cr{c-12}\|u\|\}
\label{del13}
\end{equation}
by \eqref{dec4}, \eqref{dec8} and \eqref{del2}. \eqref{del3} for $m+1$ is also satisfied
due to \eqref{dec8.5}. Thus \eqref{del1}-\eqref{del3} for $m+1$ in place of $m$ are
verified. We next check that $J^{(m)}$ converges. By \eqref{del1} and \eqref{del6},
$\{\tau^{(m)}\}$ is a Cauchy sequence and 
\begin{equation}
|\tau^{(m)}|\leq \sum_{j=1}^m d_{\theta_*^{j-1}}(J^{(j)},J^{(j-1)})
\leq \frac{\Cr{c-13}}{2}
\max\{\mu,\Cr{c-12}\|u\|\}
\label{del14}
\end{equation}
and $\tau^{(m)}$ converges to some $\hat\tau$. We also have
\begin{equation}
|a^{(m)}|\leq  \sum_{j=1}^m |a^{(j)}-a^{(j-1)}|\leq \sum_{j=1}^{m}
\theta_*^{j-1} d_{\theta_*^{j-1}}(J^{(j)},J^{(j-1)})\leq  \frac{\Cr{c-13}}{2}
\max\{\mu,\Cr{c-12}\|u\|\}
\label{del15}
\end{equation}
and $a^{(m)}$ converges to some $\hat a$. Setting $J_0={\bf R}_{\hat\tau}(J)+\hat a$, 
we prove \eqref{del9} by \eqref{del14} and \eqref{del15}. Next, for $\theta_*^{m+1}
\leq s<\theta_*^m$, let $J_s=J^{(m)}$. Then,
\begin{equation}
\begin{split}
d_s(J_s,J_0)&\leq \sum_{j=m+1}^{\infty}(s^{-1}|a^{(j)}-a^{(j-1)}|+|\tau^{(j)}-\tau^{(j-1)}|)
\leq \sum_{j=m+1}^{\infty} (s^{-1}\theta_*^{j-1}+1)d_{\theta_*^{j-1}}(J^{(j)},J^{(j-1)}) \\
&\leq \sum_{j=m+1}^{\infty} (\theta_*^{-1}+1)\Cr{c-12}\theta_*^{(j-1)\zeta}\max\{\mu,
\Cr{c-12}\|u\|\}\leq \frac{s^{\zeta}}{2}\Cr{c-13}\max\{\mu,\Cr{c-12}\|u\|\}
\end{split}
\label{del16}
\end{equation}
by \eqref{del1} and \eqref{del6},
\begin{equation}
\left(s^{-5}\int_{-2s^2}^{2s^2}\int_{B_{4s}(a^{(m)})} {\rm dist}(\cdot,J_s)^2\, d\|V_t\|dt\right)^{\frac12}
\leq \theta_*^{-\frac52}\mu^{(m)}\leq \frac{s^{\zeta}}{2}\Cr{c-13}\max\{\mu,\Cr{c-12}\|u\|\}
\label{del17}
\end{equation}
by \eqref{del2} and \eqref{del7}. Thus \eqref{del16} and \eqref{del17} prove \eqref{del10}.
Now we see that the junction point $a_s$ of $J_s$ satisfies $|\hat a -a_s|=O(s^{1+\zeta})$,
and the decay of \eqref{del17} via the graphical representation shows that the dilation of 
${\rm spt}\,\|V_0\|$ centered at $\hat a$ consists of three $C^{1,\zeta}$ curves approaching to $J_0$. 
Also there cannot be any other junction point, thus the claim of (3) follows. 
\hfill{$\Box$}

\noindent
{\it Proof of Theorem \ref{mainreg}}. 
We may assume that $R=1$. For each $s\in [-1,1]$, we use
Proposition \ref{disjun3} with $R=1/2$ and the time variable 
shifted by $s$. For this purpose, 
choose $\Cr{e-8}$ so that $2^\frac52 \Cr{e-8}<\Cr{e-7}$.
Then for all $s\in [-1,1]$, we have
\begin{equation}
\label{mr8}
\left((1/2)^{-5}\int_{s-1/2}^{s+1/2}\int_{B_2}{\rm dist}(\cdot,J)^2\,
d\|V_t\|dt\right)^{\frac12}<\Cr{e-7}
\end{equation}
by \eqref{mr1}. Since $\Cr{e-8}<\Cr{e-7}$, we also have
\begin{equation}
(1/2)^{\zeta}\left(\int_{s-1/2}^{s+1/2}\left(\int_{B_2}|u|^p
\, d\|V_t\|\right)^{\frac{q}{p}}dt\right)^{\frac{1}{q}}<\Cr{e-7}
\label{mr9}
\end{equation}
by \eqref{mr3}. By Proposition \ref{smprop} and restricting
$\Cr{e-8}$ appropriately, we may guarantee that
\begin{equation}
2\|V_{s-1/2}\|(\phi_{j,J,1/2})\leq \frac{3-\nu}{2} \Cr{c-p}, \ \
2\|V_{s+1/2}\|(\phi_{j,J,1/2})\geq \frac{1+\nu}{2} \Cr{c-p}, \ \
j=1,2,3
\label{mr10}
\end{equation}
for all $s\in [-1,1]$. Since the assumptions of Proposition \ref{disjun3}
are satisfied by \eqref{mr8}-\eqref{mr10}, there exist 
$J_0(s)
\in {\mathcal J}$ with the junction point at $\hat a(s)$ and
satisfying \eqref{del9}, and for each
$\lambda\in (0,\frac12)$, 
$J_{\lambda}(s)\in {\mathcal J}$ with the junction point at
$\hat a_{\lambda}(s)$ satisfying \eqref{del10}.
The point $\hat a(s)$ is the unique junction point of 
${\rm spt}\,\|V_s\|$ in $B_1$.
We next prove $\hat a\in C^{\frac{1+\zeta}{2}}$ and
the corresponding norm estimate. 
Since \eqref{del9} gives $\sup_{s\in [-1,1]}|\hat a(s)|$
bound (with $2^{\frac52}\Cr{c-13}$ in place of $\Cr{c-13}$), we only 
need to check the $\frac{1+\zeta}{2}$-H\"{o}lder norm bound of $\hat a$.
For any $-1\leq s_1<s_2\leq 1$, set $\lambda=\min\{\frac14,\sqrt{s_2
-s_1}\}$. Then by \eqref{del10}, we have
\begin{equation}
\begin{split}
d_{\lambda}(J_{\lambda}(s_1),J_0(s_1))&+\big(\lambda^{-5}
\int_{s_1-2\lambda^2}^{s_1+2\lambda^2}\int_{B_{4\lambda}
(\hat a_{\lambda}(s_1))}{\rm dist}(\cdot,J_{\lambda}(s_1))^2\,
d\|V_t\|dt\big)^{\frac12} \\ &
\leq 2^{\frac52}\Cr{c-13}\lambda^{\zeta}
\max\{\mu,\Cr{c-12}\|u\|\}.
\end{split}
\label{mr11}
\end{equation}
Since $[s_2-\frac12\lambda^2,s_2+\frac12\lambda^2]
\subset [s_1-2\lambda^2,s_1+2\lambda^2]$, by restricting
$\Cr{e-8}$ so that
$2^5 \Cr{c-13}\lambda^{\zeta}\max\{\mu,\Cr{c-12}\|u\|\}$ $<
\Cr{e-7}$ holds, we may use Proposition \ref{disjun3} again
with $R=\frac{\lambda}{2}$ and $J$ replaced by $J_{\lambda}(s_1)$.
Since we know that $J_0(s_2)$ is the unique triple junction of
${\rm spt}\|V_{s_2}\|$, \eqref{del9} gives
\begin{equation}
d_{\frac{\lambda}{2}}(J_{0}(s_2),J_{\lambda}(s_1))\leq
\Cr{c-13}\max\{2^5 \Cr{c-13}\lambda^{\zeta}\max\{\mu,\Cr{c-12}\|u\|\},
\Cr{c-12}\big(\frac{\lambda}{2}\big)^{\zeta}\|u\|\}.
\label{mr12}
\end{equation}
Now \eqref{mr11} and \eqref{mr12} give the desired $C^{\frac{1+\zeta}{2}}$
estimate for $\hat a$ with an appropriate choice of $\Cr{c-14}$. 
From the argument up to this point, it is clear that we have \eqref{mr6} 
and \eqref{mr6.5}. The estimates \eqref{mr11} and \eqref{mr12} also 
prove 
\begin{equation}
\Big\|\frac{\partial f_j}{\partial x}(l_j(\cdot),\cdot)\Big\|_{C^{\frac{\zeta}{2}}([-1,1])}
\leq \Cr{c-14}\max\{\mu,\|u\|\}.
\label{mr13}
\end{equation}
With a similar argument, one can prove (using $\lambda=x-l_j(s)$ and 
Proposition \ref{disjun3})
\begin{equation}
\Big|\frac{\partial f_j}{\partial x}(x,s)-\frac{\partial f_j}{\partial x}(l_j(s),s)\Big|
\leq \Cr{c-14} (x-l_j(s))^{\zeta}\max\{\mu,\|u\|\}.
\label{mr14}
\end{equation}
To obtain $C^{1,\zeta}$ estimate for $f_j$, fix
$-1\leq s_1<s_2\leq 1$ and $l_j(s_i)\leq x_i\leq 1$ ($i=1,2$) 
and set $\lambda=\min\{1/4,\max\{|x_1-x_2|,\sqrt{s_2-s_1}\}\}$.
If $\max\{x_1-l_j(s_1), x_2-l_j(s_2)\}\leq \lambda$, one can check
that
\begin{equation}
\Big|\frac{\partial f_j}{\partial x}(x_1,s_1)-\frac{\partial f_j}{\partial x}(x_2,s_2)\Big|\leq \Cr{c-14} \lambda^{\zeta}\max\{\mu,\|u\|\}
\label{mr15}
\end{equation}
by using \eqref{mr13}, \eqref{mr14} and the triangle inequalities. 
The similar estimate for $f_j$ (with $\lambda^{1+\zeta}$ in place of 
$\lambda^{\zeta}$ on the right-hand side) can be obtained. If
$\max\{x_1-l_j(s_1), x_2-l_j(s_2)\}> \lambda$, and assuming
$x_1-l_j(s_1)>\lambda$ without loss of generality, we may 
use Proposition \ref{disjun3} with $\lambda=x_1-l_j(s_1)$ and 
Proposition \ref{smprop} to obtain a $C^{1,\zeta}$ estimate.
With an appropriate choice of $\Cr{c-14}$, we may finish the proof of
\eqref{mr7}. 
\hfill{$\Box$}
\section{Partial regularity: Proof of Theorem~\ref{mainpa}} \label{secpart}
In this section, we combine Theorem \ref{mainreg} with a 
stratification theorem of singular sets. The general idea of
stratification using tangent cone goes back to Federer \cite{Fed}
and it has been adapted in a number of variational problems. 
Here we use a non-trivial adaptation to Brakke flows due to White \cite{White0}. 
For the proof of Theorem \ref{mainpa}, we first recall some definitions and results from 
\cite{Ilmanenp,White0}.  

(a) {\it Existence of tangent flow}.
For any fixed $(y,s)\in U\times (0,\Lambda)$ and $\lambda>0$, define
\begin{equation}
V^{(y,s),\lambda}_t(\phi)=\lambda^{-1}\int_{G_1(\lambda^{-1}(U-y))} \phi(y+\lambda x, S)\, d V_{s+\lambda^2 t}(x,S)\label{pat1}
\end{equation}
 for $\phi\in C_c(G_1(\lambda^{-1}(U-y)))$ and 
 $t\in (-\lambda^{-2} s,\lambda^{-2}(\Lambda-s))$.
$V^{(y,s),\lambda}_t$ is a parabolically rescaled flow at $(y,s)$. For any
positive sequence $\{\lambda_i\}_{i\in {\mathbb N}}$ converging to 0, 
there exist a subsequence $\{\lambda_{i_j}\}_{j\in {\mathbb N}}$
and a family of varifolds $\{\tilde V_t\}_{t\in {\mathbb R}}$ with the following 
property: we have
$\tilde V_t\in {\bf IV}_1
({\mathbb R}^2)$ for a.e$.$ $t\in {\mathbb R}$, $\{\tilde V_t\}_{t\in
{\mathbb R}}$ satisfies \eqref{meq} with $u=0$ on ${\mathbb R}^2\times{\mathbb R}$,
$\tilde V_t=\tilde V_{\lambda^2 t}^{(0,0),\lambda}$ for all $t<0$ and $\lambda>0$, and
$\lim_{j\rightarrow\infty}\|V^{(y,s),\lambda_{i_j}}_t\|=\|\tilde V_t\|$ for all $t\in 
{\mathbb R}$. The proof of the existence of such flow is in 
\cite[Lem. 8]{Ilmanenp} and also see \cite[Sec. 7]{White0}. 
It is for $u=0$, 
but the proof goes through
even with non-zero $u$ since Huisken's monotonicity formula \cite{Huisken} 
holds with
a minor error term due to $u$ (see \cite[Sec. 6]{Kasai-Tonegawa} for the
detail) and it
vanishes as $\lambda\rightarrow 0$. We call $\{\tilde V_t\}_{t\in {\mathbb R}}$
a {\it tangent flow} at $(y,s)$. Note that $\{\tilde V_t\}_{t\in {\mathbb R}}$ inherits the property of  
(A2) with the same constant $E_1$.

(b) {\it backwards-cone-like functions}.
For any 
tangent flow $\{\tilde V_t\}_{t\in 
{\mathbb R}}$ at $(y,s)\in U\times(0,\Lambda)$ and for $(x,t)\in {\mathbb R}^2\times {\mathbb R}$, define
\begin{equation}
g(x,t)=\lim_{\tau\rightarrow 0+}
\int_{{\mathbb R}^2}\rho_{(x,t)}(x',t-\tau)\, d\|\tilde V_{t-\tau}\|(x').
\label{pat2}
\end{equation}
By Huisken's monotonicity formula, the limit in \eqref{pat2} always exists. 
The set of all such $g$ obtained from a tangent flow at $(y,s)$ is denoted
by ${\mathcal G}(y,s)$.
The function $g$ has the following property which is called {\it backwards-cone-like} \cite[Sec. 8]{White0}: 
\begin{equation}
\label{pat2.1}
g(x,t)\leq g(0,0) \ \ \forall (x,t)\in {\mathbb R}^2\times {\mathbb R},
\end{equation}
\begin{equation}
g(x,t)=g(0,0) \Longrightarrow g(x+x', t+t')=g(x+\lambda x',t+\lambda^2 t')
\  \forall t'\leq 0, \ \forall x'\in {\mathbb R}^2, \ \forall \lambda>0.
\label{pat2.2}
\end{equation}
Define
\begin{equation}
{\mathcal V}(g)=\{x\in {\mathbb R}^2 : g(x,0)=g(0,0)\}, \ \
{\mathcal S}(g)=\{(x,t)\in {\mathbb R}^2\times {\mathbb R} : g(x,t)=g(0,0)\}.
\label{pat3}
\end{equation}
We note that ${\mathcal V}(g)$ and ${\mathcal S}(g)$ are
denoted by $V(g)$ and $S(g)$ in \cite{White0}, respectively, but we changed the 
notation here to avoid possible confusion. 
Then \cite[Th. 8.1]{White0} proves that 
\begin{equation}
\label{pat3.5}
g(x,t)=g(x+x',t) \ \ \forall x\in {\mathbb R}^2, \, \forall x' \in {\mathcal V}(g),\,
\forall t\leq 0,
\end{equation}
${\mathcal V}(g)$ is a vector subspace of ${\mathbb R}^2$, and
${\mathcal S}(g)$ is either ${\mathcal V}(g)\times \{0\}$ or ${\mathcal V}(g)
\times (-\infty,a]$ for some $a\in [0,\infty]$. In the latter case, $g$ is
time-independent up to time $t=a$: that is, $g(x,t)=g(x,t')$ for all $t\leq t'<a$
and $x\in {\mathbb R}^2$. 
Depending on whether ${\mathcal S}(g)$ is equal to ${\mathcal V}(g)\times {\mathbb R}$,
${\mathcal V}(g)\times
(-\infty,a]$ for some $a\in [0,\infty)$ or ${\mathcal V}(g)\times \{0\}$, $g$
is called {\it static}, {\it quasi-static} or {\it shrinking}, 
respectively. ${\mathcal V}(g)$
is called the {\it spatial spine} of $g$. In the present situation of ${\mathbb R}^2$, dimension
of ${\mathcal V}(g)$ denoted by ${\rm dim}{\mathcal V}(g)$ can be either
2, 1 or 0.  

(c) {\it Stratification}.
Define for $g\in {\mathcal G}(y,s)$ 
\begin{equation}
{\mathcal D}(g)=\left\{\begin{array}{ll}2+{\rm dim}{\mathcal V}(g) &
\mbox{ if $g$ is static}, \\  
{\rm dim}{\mathcal V}(g) &\mbox{ if $g$ is quasi-static or
shrinking,}
\end{array}
\right.
\label{pat4}
\end{equation}
and define for $k\in \{0\}\cup{\mathbb N}$ 
\begin{equation}
\Sigma_k=\{(y,s)\in U\times(0,\Lambda): {\mathcal D}(g)\leq k\ \ \forall
g\in {\mathcal G}(y,s)\}.
\label{pat5}
\end{equation}
Then \cite[Th. 8.2]{White0} proves that ${\rm dim}\,\Sigma_k\leq k$ and 
$\Sigma_0$ is a discrete set. Here, ${\rm dim}$ is the Hausdorff 
dimension with respect to the parabolic metric. 

{\it Proof of Theorem \ref{mainpa}}. Define $\Sigma_1$ as in \eqref{pat5} 
whose parabolic Hausdorff dimension is at most 1. Let $(y,s)\in U\times(0,\Lambda)
\setminus \Sigma_1$. By the definition of \eqref{pat5}, there exists $g\in {\mathcal G}(y,s)$ with ${\mathcal D}(g)\geq 2$. In the following, fix such $g$.
\newline
{\it Claim 1}.
There exists a constant $g_0>0$ depending only on $E_1$ 
such that either
$g(x,t)\geq g_0$ or $g(x,t)=0$. In the latter case, there exists a 
space-time neighborhood $U_{x,t}$ of $(x,t)$ with $U_{x,t}\cap
\cup_{t'}({\rm spt}\,\|\tilde V_{t'}\|\times\{t'\})=\emptyset$. 
\newline
{\it Proof of claim 1}. This is a well-known fact but we include the proof
for the convenience of the reader. By Brakke's clearing out lemma \cite{Brakke} or 
\cite[Cor.\,6.3]{Kasai-Tonegawa}, there exist constants $\tilde g_0>0$ 
and $L>1$ depending on $E_1$
such that for any $\tau>0$, $\|\tilde V_{t-2\tau}\|(B_{\sqrt{\tau}L}(x))<\tilde g_0 \sqrt{\tau}$ implies $\|\tilde V_{t'}\|
(B_{\sqrt{\tau}}(x))=0$ for all $t'\in [t-\tau,t+\tau]$. Assume 
$g(x,t)>0$. By the monotone property of \eqref{pat2}, for sufficiently small $\tau$, we have
\begin{equation*}
2 g(x,t)>\int_{B_{\sqrt{\tau} L}(x)} \rho_{(x,t)}(x',t-2\tau)\, d\|\tilde V_{t-2\tau}\|(x')
\geq \frac{e^{-\frac{L^2}{8}}}{\sqrt{8\pi\tau}} \|\tilde V_{t-2\tau}\|(B_{\sqrt{\tau}
L}(x)). \label{par5.1}
\end{equation*}
Thus, for sufficiently small $g_0$ depending on $\tilde g_0$ and $L$, 
if $g(x,t)<g_0$, we have $\|\tilde V_{t'}\|(B_{\sqrt{\tau}}(x))=0$ as above
and this implies $g(x,t)=0$, leading to a contradiction. Thus we have $g(x,t)
\geq g_0$ if $g(x,t)>0$. If $g(x,t)=0$, then the same argument shows
the last statement, concluding the proof of claim 1. 
\newline
{\it Claim 2}. ${\rm dim}{\mathcal V}(g)=2$ implies that there exists
a space-time neighborhood $U_{y,s}\subset U\times(0,\Lambda)$ such that
$U_{y,s}\cap \cup_{t}({\rm spt}\,\|V_t\|\times\{t\})=\emptyset$. 
\newline
{\it Proof of claim 2}. By \eqref{pat3}, $g(\cdot,0)$ is a constant function on 
${\mathbb R}^2$. Suppose $g(0,0)\geq g_0$. By the monotone property of
\eqref{pat2} and using (A2), we have for any $x\in {\mathbb R}^2$,
$\tau>0$ and $R>0$
\begin{equation}
g_0\leq g(x,0)\leq \int_{B_{\sqrt{\tau} R}(x)}\rho_{(x,0)}(x',-\tau)\, d\|\tilde V_{-\tau}\|(x')
+E_1 o(1)
\label{pat6}
\end{equation}
where $o(1)$ here means $\lim_{R\rightarrow\infty}o(1)=0$. Hence, 
fixing large $R$ depending only on $E_1$ so the last term is less than 
$\frac{g_0}{2}$, and then set $\delta=\frac{g_0}{2} \sqrt{4\pi}$. 
With this choice and \eqref{pat6} show
\begin{equation}
\delta\sqrt{\tau}\leq \|\tilde V_{-\tau}\|(B_{\sqrt{\tau} R}(x))
\label{pat7}
\end{equation}
for all $\tau>0$ and $x\in {\mathbb R}^2$. But then, 
since $B_R$ may contain $O(\tau^{-1})$ number of disjoint balls
of radius $\sqrt{\tau}R$, we may prove that
$\|\tilde V_{-\tau}\|(B_R)\geq C \delta/\sqrt{\tau}$
which goes to infinity as $\tau\rightarrow 0$. This is a contradiction
to (A2). Thus $g(0,0)=0$, and by \eqref{pat2.1}, $g$ is identically 0 on ${\mathbb R}^2
\times {\mathbb R}$. Then claim 1 shows 
$\|\tilde V_t\|=0$ for all $t\in {\mathbb R}$. 
Let $\{\lambda_j\}_{j=1}^{\infty}$ be a sequence such that 
$\|V_t^{(y,s),\lambda_j}\|\rightarrow \|\tilde V_t\|(=0)$ as was
described in (a). Then one has $\lim_{j\rightarrow \infty}
\|V_{t'-2}^{(y,s),\lambda_j}\|(B_{L})=0$ for $t'\in [-1,1]$, $L$ as in claim 1.
Then by \cite[Cor.6.3]{Kasai-Tonegawa}, 
for sufficiently large $j$, we have
$\|V_{t'}^{(y,s),\lambda_j}\|(B_{1})=0$ for all $t'\in [-1,1]$.
This shows that there exists a space-time neighborhood of $(y,s)$ 
on which $\|V_t\|$ has measure zero, completing the proof of claim 2. 
\newline
{\it Claim 3}. ${\rm dim}{\mathcal V}(g)=1$ implies that there exists a space-time
neighborhood $U_{y,s}\subset U\times(0,\Lambda)$ such that
$U_{y,s}\cap \cup_{t} ({\rm spt}\,\|V_t\|\times \{t\})$ is represented 
as a $C^{1,\zeta}$ graph. 
\newline
{\it Proof of claim 3}. Since ${\mathcal D}(g)\geq 2$, \eqref{pat4} shows
that $g$ is static, i.e., ${\mathcal S}(g)={\mathcal V}(g)\times{\mathbb R}$.
We have $g(0,0)\geq g_0$, or else, $g$ is identically 0 and ${\rm dim}{\mathcal V}(g)=2$.
As described in (b), $g$ is independent of $t$, and by \eqref{pat3.5}, 
invariant in ${\mathcal V}(g)$ direction. Thus, if $g(x,t)(=g(x,0))>0$ for some
$x\in {\mathbb R}^2\setminus {\mathcal V}(g)$, then $g(\lambda x+x',t)
=g(x,t)$ for all $x'\in {\mathcal V}(g)$ and $\lambda>0$, letting
$g$ having a positive constant on the half-space of ${\mathbb R}^2$
with the boundary ${\mathcal V}(g)$.
This leads to a contradiction by the similar argument in the proof of claim 2.
Thus we have $g=0$ outside of ${\mathcal S}(g)$ and positive constant
on ${\mathcal S}(g)$. Similarly, we may also prove that $\|\tilde V_t\|(
{\mathbb R}^2\setminus {\mathcal V}(g))=0$ for all $t\in {\mathbb R}$.
For a.e$.$ $t\in {\mathbb R}$, we have 
$\tilde V_t\in {\bf IV}_1({\mathbb R}^2)$, thus 
$\tilde V_t=\theta(x,t)|{\mathcal V}(g)|$ for some ${\mathcal H}^1$ a.e$.$ 
integer-valued function $\theta(\cdot,t)$. Since 
$h({\tilde V_t},\cdot)\in L^2_{loc}$ for a.e$.$ $t$, going back to the
definition of the first variation, one can check that $\theta(\cdot,t)$ has to 
be a constant function. 
Since $g$ is constant on ${\mathcal S}(g)$,
one easily sees that $\theta$ is independent of $t$. 
Thus with some integer $\theta_0$, $\tilde V_t=\theta_0 |{\mathcal V}(g)|$ 
for all $t$. Now by (A5), we necessarily have $\theta_0=1$. Let $V_t^{(y,s),\lambda_j}$ be a 
sequence converging to $\tilde V_t$. 
Then by \cite[Th.\,8.7 or Prop.\,9.1]{Kasai-Tonegawa}, 
for sufficiently large $j$, we may conclude that $\cup_t ({\rm spt}\,\|V_t^{(y,s),
\lambda_j}\|\times \{t\})$ is represented as a $C^{1,\zeta}$ graph in 
$B_1\times (-1,1)$. This concludes the proof of claim 3. 
\newline
{\it Claim 4}. ${\rm dim}{\mathcal V}(g)=0$ implies that there exists a space-time
neighborhood $U_{y,s}\subset U\times(0,\Lambda)$ such that
$U_{y,s}\cap \cup_{t} ({\rm spt}\,\|V_t\|\times \{t\})$ is represented 
as a $C^{1,\zeta}$ triple junction as in Theorem \ref{mainreg}. 
\newline
{\it Proof of claim 4}. Since ${\mathcal D}(g)\geq 2$, by \eqref{pat4}, 
$g$ is static. $g$ is independent of $t$ and 
$g(x,t)=g(\lambda x,0)$ for all $x\in{\mathbb R}^2$ and $\lambda>0$. 
We shall write $g(x)$ instead of $g(x,t)$. 
Define $W=\{x\in {\mathbb R}^2 : |x|=1,\, g(x)\geq g_0\}$. Then 
following the similar argument as in the proof of claim 2, we may prove
that the number of element of $W$ is finite. More precisely, if we pick
$W'\subset W$ consisted of $N$ elements, 
we may choose $O(N/\sqrt{\tau})$ number of
disjoint balls of radius $\sqrt{\tau} R$ centered at $\cup_{\lambda>0}
\lambda W'$ inside of $B_{R}$, each satisfying \eqref{pat7}. 
Then we would have $\|\tilde V_{-\tau}(B_R)\|\geq O(N)$, thus 
$N$ cannot go to $\infty$. Thus $\{g\geq g_0\}$ consists of 
a finite number of half rays denoted by $\{l_j \subset {\mathbb R}^2: j=1,\cdot, N\}$
emanating from 0, and $g$ is constant
on each half ray. As in the proof of claim 3, one can argue that 
there exist some positive integers $\theta_j$ such that $\tilde V_t
=\sum_{j=1}^{N} \theta_j |l_j|$. Then again using (A5) which says $\Theta(\|V_s\|,y)<2$
and arguing as before, we have $N\leq 3$ and $\sum_{j=1}^N\theta_j
\leq 3$. The conditions $h(\tilde V_t,\cdot)
\in L^2_{loc}$ and ${\rm dim}{\mathcal V}(g)=0$ limit the possibility to 
$\tilde V_t=|{\bf R}_{\theta}(J)|$ with some $\theta\in [0,2\pi)$. 
We are now ready to apply Theorem \ref{mainreg} to 
$V_t^{(y,s),\lambda_j}$. Note that (after a rotation by $\theta$) 
that \eqref{mr1}-\eqref{mr3} are satisfied for all sufficiently large $j$.
Thus this concludes the proof of claim 4. Since $U_{y,s}$ in all three
cases do not intersect with $\Sigma_1$, $\Sigma_1$ is a closed set 
and this ends the proof of Theorem \ref{mainpa}.
\hfill{$\Box$}

\section{The top dimensional part of the genuine singular set}
\label{difpat}
Under the hypotheses (A1)-(A5) of Theorem \ref{mainpa}, let $\Sigma_1$ and $\Sigma_0$
be defined as in \eqref{pat5}. We know that $\Sigma_1$ is closed (by Theorem \ref{mainpa}), 
and that ${\rm dim} \, \Sigma_{1} \leq 1$ where ${\rm dim}$ is the parabolic Hausdorff dimension. Moreover,  
$\Sigma_0$ is discrete by \cite[Th.\,8.2]{White0}. 

We may further characterize the ``top dimensional part'' of the singular set, i.e.\  $\Sigma_1
\setminus \Sigma_0,$ in terms of tangent flows as follows: 
\begin{thm}
For any $(y,s)\in \Sigma_1\setminus \Sigma_0$, there exists a 
quasi-static tangent flow $\{\tilde V_t\}_{t\in {\mathbb R}}$ such that, ${\rm spt} \, \|\tilde V_{t}\| = \emptyset$ or 
${\rm spt}\,\|\tilde V_t\|=S$ for some $S\in {\bf G}(2,1)$. Moreover, there exists a set of integers $\theta_1>\cdots
>\theta_N\geq 0$ ($N\geq 2$) and real numbers $0\leq a_1<\cdots<a_{N-1}$ such that 
$\tilde V_t=\theta_j |S|$ for $t\in (a_{j-1},a_j)$, $j=1,\cdots,N$, where $a_0=-\infty$ and $a_N=\infty$.
If $\theta_1=1$, then $a_1=0$ and $\theta_2=0$. 
\label{difpatreg}
\end{thm}
A heuristic meaning of the above 
is that $\Sigma_1\setminus \Sigma_0$ is the set of points where some curve instantaneously
disappears. 

{\it Proof}. Take any $(y,s)\in \Sigma_1\setminus \Sigma_0$.
Then there exists $g\in {\mathcal G}(y,s)$ such that ${\mathcal D}(g)= 1$. 
By \eqref{pat4}, it has to be quasi-static or shrinking, and ${\rm dim}{\mathcal V}(g)=1$. 
By \eqref{pat3.5}, $g$ is invariant in ${\mathcal V}(g)$
direction while having the backwards-cone-like property \eqref{pat2.2} with respect to $(0,0)$.
The set $\{x\in {\mathbb R}^2: g(x,-1)>0\}$ then consists of a finite number of lines 
parallel to ${\mathcal V}(g)$ and by the argument in the proof of Theorem \ref{mainpa}, 
one can prove that $\tilde V_{-1}$ is a sum of varifolds with integer multiples supported on such lines. Due to the
backwards-cone-like property and also the fact that $\tilde V_t$ is a curvature flow, one
can prove that $\tilde V_t=\theta_1 |{\mathcal V}(g)|$ ($\theta_1\in {\mathbb N}$) for $t<0$ (otherwise it has to move
at non-zero speed even if it is a line). This shows that $g$ has to be quasi-static. 
By \cite{White0}, we know that ${\mathcal S}(g)={\mathcal V}(g)\times(-\infty,a]$ for some $a\in [0,\infty)$.
This in particular shows $\tilde V_t=\theta_1|{\mathcal V}(g)|$ for $t\in (-\infty,a)$. 
One can then prove, for example using the clearing-out lemma \cite{Brakke}, that 
${\rm spt}\,\|\tilde V_t\|\subset {\mathcal V}(g)$ for all $t>0$ and by $h(\tilde V_t,\cdot)\in L^2_{loc}$, 
that $\tilde V_t=\theta(t)|{\mathcal V}(g)|$ for some $\theta(t)\in {\mathbb N}$. The fact that
they are time-discretely decreasing can be easily seen from the curvature flow inequality.
If $\theta_1=1$, and if $a_1>0$, then this would mean that ${\tilde V}_t=|{\mathcal V}(g)|$
in a neighborhood of $(0,0)$. Since $V_t^{(y,s),\lambda_j}$ is approaching to ${\tilde V}_t$,
for sufficiently large $j$, we may apply \cite[Th.\,8.7]{Kasai-Tonegawa} to $V_t^{(y,s),\lambda_j}$ 
in some small neighborhood of $(0,0)$ and conclude that $(y,s)$ is a $C^{1,\zeta}$ regular point of $V_t$
(see the definition in \cite{Kasai-Tonegawa}). But then 
the tangent flow at $(y,s)$ should be static, a contradiction. Thus if $\theta_1=1$, then 
$a_1=0$. This completes the proof.
\hfill{$\Box$}

\noindent
{\bf Remark:} If we assume further that there exists no quasi-static tangent flow with ${\rm dim}\,{\mathcal V}(g)=1$, 
Theorem \ref{difpatreg} 
shows that $\Sigma_1\setminus\Sigma_0=\emptyset$. If this is satisfied,
the picture is akin to that of the motion of grain boundaries
where networks of curves joined by triple junctions move continuously with occasional
collisions of junctions only at discrete points in $\Sigma_0$. 

\end{document}